\documentclass{amsart}

\usepackage{amsmath} 
\usepackage{amssymb}

\newtheorem{theorem}{Theorem}[section] 
\newtheorem{claim}{Claim}[theorem]
 
\newtheorem{proposition}[theorem]{Proposition} 
\newtheorem{corollary}[theorem]{Corollary} 

\theoremstyle{definition}
\newtheorem{definition}[theorem]{Definition}

\newtheorem{problem}[theorem]{Problem}

\theoremstyle{remark}
\newtheorem{remark}[theorem]{Remark}
\newtheorem{notation}[theorem]{Notation}

\newtheorem{context}[theorem]{Context}

\numberwithin{equation}{section}
\newcommand{\forces}{\Vdash} 
\newcommand{\bV}{{\bf V}} 
 
\newcommand{\conc}{{}^\frown\!}
\newcommand{\lh}{{\rm lh}} 
\newcommand{\rest}{{\restriction}}
\newcommand{\mrot}{{\rm root}\/} 
\newcommand{\suc}{{\rm succ}} 

\newcommand{\cf}{{\rm cf}}
\newcommand{\nor}{{\bf nor}} 
\newcommand{\pos}{{\bf pos}}
\newcommand{\dis}{{\bf dis}}
\newcommand{\val}{{\bf val}}
\newcommand{\ad}{\alpha_{\rm dn}}
\newcommand{\au}{\alpha_{\rm up}}
\newcommand{\rht}{{\rm ht}}
\newcommand{\otp}{{\rm otp}}
\newcommand{\TCR}{{\rm TCR}^\lambda}
\newcommand{\CR}{{\rm CR}^\lambda}
\newcommand{\tree}{{\rm tree}}
\newcommand{\vtl}{\vartriangleleft}
\newcommand{\llim}{{\lim}_\lambda}

\newcommand{\bairel}{{}^{\lambda}\lambda}
\newcommand{\exl}{{{\rm ex}(\bar{\lambda})}}
\newcommand{\cov}{{\rm cov}}

\newcommand{\Dom}{{\rm Dom}}
\newcommand{\Rang}{{\rm Rng}}
\newcommand{\dom}{{\rm Dom}} 

\newcommand{\Gsg}{{\Game^{\rm Sacks}_{\bar{\lambda}}}}
\newcommand{\Gbd}{{\Game^{\rm bd}_\lambda}}
\newcommand{\Gfl}{{\Game^{\rm fuzzy}_\lambda}}
\newcommand{\st}{{\rm st}}
\newcommand{\nst}{\name{\rm st}}
\newcommand{\rk}{{\rm rk}}
\newcommand{\proj}{{\rm proj}}
\newcommand{\vare}{\varepsilon}
\newcommand{\cl}{{\rm cl}}
\newcommand{\cA}{{\mathcal A}}

\newcommand{\bbD}{{\mathbb D}}
\newcommand{\cD}{{\mathcal D}}

\newcommand{\bH}{{\bf H}}
\newcommand{\cH}{{\mathcal H}}

\newcommand{\cG}{{\mathcal G}}
\newcommand{\cI}{{\mathcal I}}
\newcommand{\bcI}{{\bar{\cI}}}
\newcommand{\cJ}{{\mathcal J}}

\newcommand{\bbP}{{\mathbb P}}
\newcommand{\cP}{{\mathcal P}}

\newcommand{\bbQ}{{\mathbb Q}}

\newcommand{\gR}{{\mathfrak R}}

\newcommand{\cS}{{\mathcal S}}
\newcommand{\cT}{{\mathcal T}}

\newcommand{\cX}{{\mathcal X}}
\newcommand{\cY}{{\mathcal Y}}
\newcommand{\gY}{{\mathfrak Y}}

\newcount\skewfactor
\def\mathunderaccent#1#2 {\let\theaccent#1\skewfactor#2
\mathpalette\putaccentunder}
\def\putaccentunder#1#2{\oalign{$#1#2$\crcr\hidewidth
\vbox to.2ex{\hbox{$#1\skew\skewfactor\theaccent{}$}\vss}\hidewidth}}
\def\name{\mathunderaccent\tilde-3 }

\newcommand{\PART}[2]{
	\newpage    
	\begin{center} \LARGE  Part #1 

                                    #2  
        \end{center}

	\bigskip\medskip

	\renewcommand{\thesection}{#1.\arabic{section}}
}

\begin{document}

\setcounter{page}{0}

\title{Sheva--Sheva--Sheva: Large Creatures}

\author{Andrzej Ros{\l}anowski}
\address{Department of Mathematics\\
 University of Nebraska at Omaha\\
 Omaha, NE 68182-0243, USA}
\email{roslanow@member.ams.org}
\urladdr{http://www.unomaha.edu/$\sim$aroslano}
\thanks{The first author thanks the Hebrew University of Jerusalem for
 support during his visit to Jerusalem in Summer'2001 and he also
 acknowledges partial support from {\em University Committee on Research} of
 the University of Nebraska at Omaha. He also thanks his  wife,
 Ma{\l}gorzata Jankowiak--Ros{\l}anowska for supporting him when he  was
 preparing the final version of this paper.}  

\author{Saharon Shelah}
\address{Institute of Mathematics\\
 The Hebrew University of Jerusalem\\
 91904 Jerusalem, Israel\\
 and  Department of Mathematics\\
 Rutgers University\\
 New Brunswick, NJ 08854, USA}
\email{shelah@math.huji.ac.il}
\urladdr{http://www.math.rutgers.edu/$\sim$shelah}
\thanks{The research of the second author was partially supported by the
United States-Israel Binational Science Foundation. Publication 777.}  

\subjclass{Primary 03E35, 03E40; Secondary: 03E05, 03E55}
\date{February 2003}

\begin{abstract}
We develop the theory of the forcing with trees and creatures for an
inaccessible $\lambda$ continuing Ros{\l}anowski and Shelah \cite{RoSh:470},
\cite{RoSh:672}. To make a real use of these forcing notions (that is to
iterate them without collapsing cardinals) we need suitable iteration
theorems, and those are proved as well. (In this aspect we continue
Ros{\l}anowski and Shelah \cite{RoSh:655} and Shelah \cite{Sh:587},
\cite{Sh:667}.) 
\end{abstract}

\maketitle

\tableofcontents

\newpage
\section{Introduction}
The present paper has two themes. 

The first is related to the quest for the right generalization of properness
to higher cardinals (that is, for a property of forcing notions that would
play in iterations with uncountable supports similar role to that of
standard properness in CS iterations). The evidence that there is no
straightforward generalization of properness to larger cardinals was given
already in Shelah \cite{Sh:b} (see \cite[Appendix 3.6(2)]{Sh:f}).
Substantial progress has been achieved in Shelah \cite{Sh:587},
\cite{Sh:667}, but the properties there were tailored for generalizing the
case {\em no new reals\/} of \cite[Ch. V]{Sh:f}. Then Ros{\l}anowski and
Shelah \cite{RoSh:655} gave an iterable condition for not collapsing
$\lambda^+$ in $\lambda$--support iterations of $\lambda$--complete forcing
notions (with possibly adding subsets of $\lambda$). Very recently Eisworth
\cite{Ei0x} has given another property preserved in $\lambda$--support
iterations (and implying that $\lambda^+$ is not collapsed). At the moment
it is not clear if the two properties (the one of \cite{RoSh:655} and that
of \cite{Ei0x}) are equivalent, though they have similar flavour. However,
the existing iterable properties still do not cover many examples of natural
forcing notions, specially those which come naturally in the context of
$\lambda$--reals. This brings us to the second theme: developing the forcing
for $\lambda$--reals.

A number of cardinal characteristics related to the Baire space
${}^\omega\omega$, the Cantor space ${}^\omega2$ and/or the combinatorial
structure of $[\omega]^{\textstyle \omega}$ can be extended to the spaces
$\bairel,{}^\lambda 2$ and $[\lambda]^{ \textstyle\lambda}$ for any infinite
cardinal $\lambda$. Following the tradition of Set Theory of the Reals we
may call cardinal numbers defined this way for $\bairel$ (and related 
spaces) {\em cardinal characteristics of $\lambda$--reals}. The menagerie of
those characteristics seems to be much larger than the one for the
continuum. But to decide if the various definitions lead to different (and
interesting) cardinals we need a well developed forcing technology. 

There has been a serious interest in cardinal characteristics of the
$\lambda$--reals in literature. For example, Cummings and Shelah
\cite{CuSh:541} investigates the natural generalizations ${\mathfrak
b}_\lambda$ of the unbounded number and the dominating number ${\mathfrak
d}_\lambda$, giving simple constraints on the triple of cardinals
$({\mathfrak b}_\lambda,{\mathfrak d}_\lambda, 2^\lambda)$ and  proving that
any triple of cardinals obeying these constraints can be realized. In a
somewhat parallel work \cite{ShSj:643}, Shelah and Spasojevi{\v c} study
${\mathfrak b}_\lambda$ and the generalization ${\mathfrak t}_\lambda$ of
the tower number. Zapletal \cite{Za97} investigated the splitting number
${\mathfrak s}_\lambda$ -- here the situation is really complicated as the
inequality ${\mathfrak s}_\lambda>\lambda^+$ needs large cardinals. One of
the sources of interest in characteristics of the $\lambda$--reals is their
relevance for our understanding of the club filter on $\lambda$ (or the dual
ideal on non-stationary subsets of $\lambda$) -- see, e.g., 
Balcar and Simon \cite[\S5]{BS89}, Landver \cite{La90}, Matet and
Pawlikowski \cite{MPxx}, Matet, Ros{\l}anowski and Shelah
\cite{MRSh:799}. First steps toward developing forcing for $\lambda$--reals
has been done long time ago: in 1980 Kanamori \cite{Ka80} presented a
systematic treatment of the $\lambda$--perfect--set forcing in products and 
iterations. Recently,  Brown \cite{Br02}, \cite{Br0x} discussed the
$\lambda$--superperfect forcing and other tree--like forcing notions. 

Our aim in this paper is to provide tools for building forcing notions
relevant for $\lambda$--reals (continuing in this Ros{\l}anowski and Shelah
\cite{RoSh:470}, \cite{RoSh:672}) and give suitable iteration theorems (thus
continuing Ros{\l}anowski and Shelah \cite{RoSh:655}). However, we restrict
our attention to the case when $\lambda$ is a strongly inaccessible
uncountable cardinal (after all, $\aleph_0$ is inaccessible), see
\ref{incon} below.  
 
The structure of the paper is as follows. It is divided into two parts,
first one presents iteration theorems, the second one gives examples and
applications. In Section A.1 we present some basic notions and methods
relevant for iterating $\lambda$--complete forcing notions. The next
section, A.2, presents preservation of $\lambda$--analogue of the Sacks
property (in Theorem \ref{itSacks}) as well as preservation of being
$\bairel$--bounding (in Theorem \ref{itbound}). Section A.3 introduces {\em
fuzzy properness}, a more complicated variant of {\em properness over 
semi-diamonds} from \cite{RoSh:655}. Of course, we prove a suitable
iteration theorem (see Theorem \ref{1.6}). Then we give examples for the
properties discussed in Part A. We start with showing that a forcing notion
useful for uniformization is fuzzy proper (in Section B.4), and then we turn
to developing forcing notions built with the use of {\em trees and 
creatures}. In Section B.5 we set the terminology and notation, and in the
next section we discuss when the resulting forcing notions have the two
bounding properties discussed in \S A.2. Section B.7 shows how our methods
result in suitably proper forcing notions, and the last section introduces
some new characteristics of the $\lambda$--reals.  
\medskip

\noindent{\bf Notation}\qquad
Our notation is rather standard and compatible with that of classical
textbooks (like Jech \cite{J}). In forcing we keep the older convention that
{\em a stronger condition is the larger one}. Our main conventions are
listed below.  

\begin{notation}
\begin{enumerate}
\item For a forcing notion $\bbP$, $\Gamma_\bbP$ stands for the canonical
$\bbP$--name for the generic filter in $\bbP$. With this one exception, all
$\bbP$--names for objects in the extension via $\bbP$ will be denoted with a
tilde below (e.g., $\name{\tau}$, $\name{X}$). The weakest element of $\bbP$
will be denoted by $\emptyset_\bbP$ (and we will always assume that there is
one, and that there is no other condition equivalent to it). We will also
assume that all forcing notions under considerations are atomless. 

By ``$\lambda$--support iterations'' we mean iterations in which domains of
conditions are of size $\leq\lambda$. However, we will pretend that
conditions in a $\lambda$--support iteration $\bar{\bbQ}=\langle\bbP_\zeta,
\name{\bbQ}_\zeta:\zeta<\zeta^*\rangle$ are total functions on $\zeta^*$
and for $p\in\lim(\bar{\bbQ})$ and $\alpha\in\zeta^*\setminus\Dom(p)$ we
will let $p(\alpha)=\name{\emptyset}_{\name{\bbQ}_\alpha}$. 

\item For a filter $D$ on $\lambda$, the family of all $D$--positive subsets
of $\lambda$ is called $D^+$. (So $A\in D^+$ if and only if $A\subseteq
\lambda$ and $A\cap B\neq\emptyset$ for all $B\in D$.)

The club filter of $\lambda$ is denoted by $\cD_\lambda$.

\item Ordinal numbers will be denoted be the lower case initial letters of
the Greek alphabet ($\alpha,\beta,\gamma,\delta\ldots$) and also by $i,j$
(with possible sub- and superscripts). 

Cardinal numbers will be called $\theta,\kappa,\lambda,\mu$ (with
possible sub- and superscripts); $\lambda$ is a fixed inaccessible cardinal
(see \ref{incon}).

\item By $\chi$ we will denote a {\em sufficiently large} regular cardinal; 
$\cH(\chi)$ is the family of all sets hereditarily of size less than
$\chi$. Moreover, we fix a well ordering $<^*_\chi$ of $\cH(\chi)$. 

\item A bar above a letter denotes that the object considered is a sequence;
usually $\bar{X}$ will be $\langle X_i:i<\zeta\rangle$, where $\zeta$
is the length $\lh(\bar{X})$ of $\bar{X}$. Sometimes our sequences will be
indexed by a set of ordinals, say $S\subseteq\lambda$, and then $\bar{X}$
will typically be $\langle X_\delta:\delta\in S\rangle$. 

But also, $\eta,\nu$ and $\rho$ (with possible sub- and superscripts) will
denote sequences (nodes in quasi trees). 

For two sequences $\eta,\nu$ we write $\nu\vtl\eta$ whenever
$\nu$ is a proper initial segment of $\eta$, and $\nu\trianglelefteq\eta$
when either $\nu\vtl\eta$ or $\nu=\eta$. 

\item We will consider several games of two players. One player will be
called {\em Generic\/} or {\em Complete\/} or just {\em I player}, and we
will refer to this player as ``she''. Her opponent will be called {\em
Antigeneric\/} or {\em Incomplete} or just {\em II player\/} and will be
referred to as ``he''. 
\end{enumerate}
\end{notation}

\begin{definition}
\label{quasitree}
\begin{enumerate}
\item {\em A $\lambda$--quasi tree} is a set $T$ of sequences of length
${<}\lambda$ with the $\vtl$--smallest element denoted by $\mrot(T)$.  
\item A $\lambda$--quasi tree  $T$ is {\em a $\lambda$--tree} if it is
closed under initial segments longer then $\lh(\mrot(T))$. 
\item A $\lambda$--quasi tree is complete if the union of any
$\vtl$--increasing sequence of length less than $\lambda$ of members of $T$
is in $T$.  
\item For a $\lambda$--quasi tree $T$ and $\eta\in T$ we define {\em the
successors of $\eta$ in $T$}, {\em maximal points of $T$}, {\em the
restriction of $T$ to $\eta$},  and {\em the height of $T$} by: 
\[\suc_T(\eta)=\{\nu\in T: \eta\vtl\nu\ \&\ \neg(\exists\rho\in
T)(\eta\vtl\rho\vtl\nu)\},\]
\[\max(T)=\{\nu\in T:\mbox{ there is no }\rho\in T\mbox{ such that }
\nu\vtl\rho\},\]
\[T^{[\eta]}=\{\nu\in T: \eta\trianglelefteq \nu\}, \ \mbox{ and }\
\rht(T)=\sup\{\lh(\eta):\eta\in T\}.\]
   
We put $\hat{T}= T\setminus\max(T)$.

\item For $\delta<\lambda$ and a $\lambda$--quasi tree $T$ we let 
\[(T)_\delta=\{\eta\in T:\lh(\eta)=\delta\}\quad\mbox{ and }\quad (T)_{<
\delta}=\{\eta\in T:\lh(\eta)<\delta\}.\]
The set of all limit $\lambda$--branches through $T$ is 
\[\llim(T)\stackrel{\rm def}{=}\{\eta:\eta\mbox{ is a $\lambda$--sequence }\
\mbox{ and }\ (\forall \beta<\lambda)(\exists\alpha>\beta)(\eta\rest
\alpha\in T)\}.\]   

\item A subset $F$ of a $\lambda$--quasi tree $T$ is {\em a front} of $T$ if
no two distinct members of $F$ are $\vtl$--comparable and 
\[(\forall\eta\in\llim(T)\cup\max(T))(\exists \alpha<\lambda)(\eta\rest
\alpha\in F).\] 
\end{enumerate}
\end{definition}

Note that if $T$ is a complete $\lambda$--quasi tree of height $<\lambda$,
then $\max(T)$ is a front of $T$ and every $\vtl$--increasing sequence of
members of $T$ has a $\vtl$--upper bound in $\max(T)$.

In the present paper we assume the following. 
 
\begin{context}
\label{incon}
\begin{enumerate}
\item[(a)] $\lambda$ is a strongly inaccessible cardinal,
\item[(b)] $\bar{\lambda}=\langle\lambda_\alpha:\alpha<\lambda\rangle$ is a
strictly increasing sequence of uncountable regular cardinals,
$\sup\limits_{\alpha<\lambda}\lambda_\alpha=\lambda$,
\item[(c)] for each $\alpha<\lambda$,
\[\prod_{\beta<\alpha}\lambda_\beta <\lambda_\alpha\quad\mbox{ and }\quad
(\forall\xi<\lambda_\alpha)(|\xi|^\alpha<\lambda_\alpha).\] 
\end{enumerate}
\end{context}

\PART{A}{Iteration theorems for $\lambda$--support iterations}

\section{Iterations of complete forcing notions and trees of conditions}
In this section we recall some basic definitions and facts concerning
complete forcing notions and $\lambda$--support iterations.  

\begin{definition}
\label{strcom}
Let $\bbP$ be a forcing notion.
\begin{enumerate}
\item For a condition $r\in\bbP$ and a set $S\subseteq \lambda$, let
$\Game_0^\lambda(\bbP,S,r)$ be the following game of two players, {\em
Complete} and  {\em Incomplete}:   
\begin{quotation}
\noindent the game lasts $\lambda$ moves and during a play the players
construct a sequence $\langle (p_i,q_i): i<\lambda\rangle$ of pairs of
conditions from $\bbP$ in such a way that $(\forall j<i<\lambda)(r\leq p_j
\leq q_j\leq p_i)$ and at the stage $i<\lambda$ of the game:\\
if $i\in S$, then Complete chooses $p_i$ and Incomplete chooses $q_i$, and\\ 
if $i\notin S$, then Incomplete chooses $p_i$ and Complete chooses $q_i$.
\end{quotation}
Complete wins if and only if for every $i<\lambda$ there are legal moves for
both players. 
\item We say that the forcing notion $\bbP$ is {\em $(\lambda,
S)$--strategically complete} if Complete has a winning strategy in the game
$\Game_0^\lambda(\bbP,S,r)$ for each condition $r\in\bbP$. We say that 
$\bbP$ is {\em strategically $({<}\lambda)$--complete} if it is $(\lambda,
\emptyset)$--strategically complete.
\item We say that $\bbP$ is {\em $\lambda$--complete} if every
$\leq_{\bbP}$--increasing chain of length less than $\lambda$ has an upper  
bound in $\bbP$. 
\item Let $N\prec (\cH(\chi),\in,<^*_\chi)$ be a model such that
  ${}^{<\lambda} N\subseteq N$, $|N|=\lambda$ and $\bbP\in N$. We say that a
  condition $p\in\bbP$ is {\em $(N,\bbP)$--generic in the standard
  sense\/} (or just: {\em $(N,\bbP)$--generic\/}) if for every
  $\bbP$--name $\name{\tau}\in N$ for an ordinal we have $p\forces$``
  $\name{\tau}\in N$ ''. 
\item $\bbP$ is {\em $\lambda$--proper in the standard sense\/} (or just:
  {\em $\lambda$--proper\/}) if there is $x\in \cH(\chi)$  such that for
  every model $N\prec (\cH(\chi),\in,<^*_\chi)$ satisfying  
\[{}^{<\lambda} N\subseteq N,\quad |N|=\lambda\quad\mbox{ and }\quad\bbP,x
  \in N, \]
and every condition $q\in N\cap\bbP$ there is an $(N,\bbP)$--generic
condition $p\in\bbP$ stronger than $q$.
\end{enumerate}
\end{definition}

\begin{remark}
\begin{enumerate}
\item On strategic completeness (and variants) see \cite[\S A.1]{Sh:587}. 
Plainly, $\lambda$--completeness implies strategic
$({<}\lambda)$--completeness.  
\item Note that if ${\mathbb P}$ is strategically $({<}\lambda)$--complete
and $D$ is a normal filter on $\lambda$, then in $\bV^{\bbP}$ the filter $D$ 
generates a proper normal filter on $\lambda$. [Abusing notation, we may
call this filter also by $D$.]   
\end{enumerate}
\end{remark}

\begin{proposition}
\label{obsA.4}
Suppose that $\bbP$ is a $({<}\lambda)$--strategically complete (atomless)
forcing notion, $\alpha^*<\lambda$ and $q_\alpha\in\bbP$ (for $\alpha<
\alpha^*$). Then there are conditions $p_\alpha\in\bbP$ (for $\alpha<
\alpha^*$) such that $q_\alpha\leq p_\alpha$ and for distinct
$\alpha,\alpha'<\alpha^*$ the conditions $p_\alpha, p_{\alpha'}$ are
incompatible.
\end{proposition}

\begin{proof}
For $\alpha<\alpha^*$ let $\st_\alpha$ be the winning strategy of Complete
in the game $\Game_0^\lambda(\bbP,\emptyset,q_\alpha)$. By induction on
$i<\alpha^*$ we define conditions $q^i_\alpha,p^i_\alpha$ as follows:
\medskip

$p^0_0=q_0$, $q^0_0$ is the answer of Complete to $\langle p^0_0\rangle$
according to $\st_0$, $q^0_\alpha=p^0_\alpha=q_\alpha$ for $\alpha>0$.

Suppose that conditions $p^j_\alpha,q^j_\alpha$ have been defined for
$j<i$, $\alpha<\alpha^*$ (where $i<\alpha^*$) so that 
\begin{enumerate}
\item[$(\alpha)$] $(\forall \alpha<\alpha'<i)(q^{\alpha'}_\alpha,
q^{\alpha'}_{\alpha'}\mbox{ are incompatible })$,
\item[$(\beta)$]  for each $\alpha<i$, $\langle (p^j_\alpha,q^j_\alpha):
\alpha\leq j< i\rangle$ is a play of $\Game_0^\lambda(\bbP,\emptyset,
q_\alpha)$ in which Complete uses the strategy $\st_\alpha$, and
\item[$(\gamma)$] $p^j_\alpha=q^j_\alpha=q_\alpha$ for $\alpha\geq i>j$.
\end{enumerate}
For $\alpha<i$ let $r_\alpha$ be a condition stronger than all $q^j_\alpha$
for $j<i$ (there is one by $(\beta)$). If every $r_\alpha$ (for $\alpha<i$)
is incompatible with $q_i$, then we let $p^i_\alpha=r_\alpha$ for
$\alpha<i$, $p^i_\alpha=q_\alpha$ for $\alpha\geq i$. Otherwise, let
$\alpha_0<i$ be the first such that $r_{\alpha_0},q_i$ are compatible. Then
we may pick two incompatible conditions $p^i_{\alpha_0},p^i_i$ above both
$r_{\alpha_0}$ and $q_i$. Next we let $p^i_\alpha=r_\alpha$ for $\alpha<i$,
$\alpha\neq\alpha_0$ and $p^i_\alpha=q_\alpha$ for $\alpha>i$. Finally, for
$\alpha\leq i$, $q^i_\alpha$ is defined as the answer of Complete according
to $\st_\alpha$ to $\langle (p^j_\alpha,q^j_\alpha):j<i\rangle\conc\langle
p^i_\alpha\rangle$, and $q^i_\alpha=q_\alpha$ for $\alpha>i$. 
\medskip

After the inductive definition is carried out we may pick upper bounds
$p_\alpha$ to $\langle q^j_\alpha:j<\alpha^*\rangle$ (for $\alpha<\alpha^*$;
exist by $(\beta)$). The conditions $p_\alpha$ are pairwise incompatible by
$(\alpha)$, so we are done. 
\end{proof}

Both completeness and strategic completeness are preserved in iterations: 

\begin{proposition}
\label{first}
Suppose that $\langle \bbP_\alpha,\name{\bbQ}_\alpha:\alpha<\zeta^*\rangle$
is a $\lambda$--support iteration such that for each $\alpha<\zeta^*$ 
\[\forces_{\bbP_\alpha}\mbox{`` $\name{\bbQ}_\alpha$ is
$\lambda$--complete. ''}\]
Then the forcing $\bbP_{\zeta^*}$ is $\lambda$--complete. 
\end{proposition}

\begin{proposition}
\label{pA.6}
Suppose $\bar{\bbQ}=\langle\bbP_i,\name{\bbQ}_i: i<\gamma\rangle$ is a
$\lambda$--support iteration and for each $i<\gamma$
\[\forces_{\bbP_i}\mbox{`` $\name{\bbQ}_i$ is strategically
$({<}\lambda)$--complete ''.}\]
Then:
\begin{enumerate}
\item[(a)] $\bbP_\gamma$ is strategically $({<}\lambda)$--complete.
\item[(b)] Moreover, for each $\vare\leq\gamma$ and $r\in\bbP_\vare$ there
is a winning strategy $\st(\vare,r)$ of Complete in the game
$\Game_0^\lambda(\bbP_\vare,\emptyset,r)$ such that, whenever $\vare_0<
\vare_1\leq\gamma$ and $r\in\bbP_{\vare_1}$, we have:
\begin{enumerate}
\item[(i)]  if $\langle (p_i,q_i):i<\lambda\rangle$ is a play of
$\Game_0^\lambda(\bbP_{\vare_0},\emptyset,r\rest\vare_0)$ in
which Complete follows the strategy
$\st(\vare_0,r\rest\vare_0)$,  

then $\langle (p_i\conc r\rest [\vare_0,\vare_1),q_i\conc r\rest 
[\vare_0,\vare_1)):i<\lambda\rangle$ is a play of
$\Game_0^\lambda(\bbP_{\vare_1},\emptyset,r)$ in which Complete uses
$\st(\vare_1,r)$; 
\item[(ii)] if $\langle (p_i,q_i):i<\lambda\rangle$ is a play of
$\Game_0^\lambda(\bbP_{\vare_1},\emptyset,r)$ in which Complete plays
according to the strategy $\st(\vare_1,r)$, 

then $\langle (p_i\rest\vare_0,q_i\rest\vare_0):i<\lambda
\rangle$ is a play of $\Game_0^\lambda(\bbP_{\vare_0},\emptyset,r\rest
\vare_0)$ in which Complete uses $\st(\vare_0,r)$;
\item[(iii)] if $\langle (p_i,q_i):i<i^*\rangle$ is a partial play of
$\Game_0^\lambda(\bbP_{\vare_1},\emptyset,r)$ in which Complete uses
$\st(\vare_1,r)$ and $p'\in\bbP_{\vare_0}$ is stronger than all
$p_i\rest\vare_0$ (for $i<i^*$), then there is $p^*\in
\bbP_{\vare_1}$ such that $p'=p^*\rest\vare_0$ and $p^*\geq p_i$
for $i<i^*$. 
\end{enumerate}
\end{enumerate}
\end{proposition}

\begin{definition}
[Compare {\cite[A.3.3, A.3.2]{Sh:587}}]
\label{dA.5}
\begin{enumerate}
\item Let $\alpha,\gamma$ be ordinals, $\emptyset\neq w\subseteq \gamma$.
{\em A standard $(w,\alpha)^\gamma$--tree\/} is a pair $\cT=(T,\rk)$ such
that  
\begin{itemize}
\item $\rk:T\longrightarrow w\cup\{\gamma\}$, 
\item if $t\in T$ and $\rk(t)=\vare$, then $t$ is a sequence $\langle
(t)_\zeta: \zeta\in w\cap\vare\rangle$, where each $(t)_\zeta$ is a
sequence of length $\alpha$,
\item $(T,\vtl)$ is a tree with root $\langle\rangle$ and such that every
chain in $T$ has a $\vartriangleleft$--upper bound it $T$.
\end{itemize}
We will keep the convention that $\cT^x_y$ is $(T^x_y,\rk^x_y)$.
\item Suppose that $w_0\subseteq w_1\subseteq\gamma$, $\alpha_0\leq
\alpha_1$, and $\cT_1=(T_1,\rk_1)$ is a standard $(w_1,\alpha_1
)^\gamma$--tree. {\em The projection $\proj^{(w_1,\alpha_1)}_{(w_0,
\alpha_0)}(\cT_1)$ of $\cT_1$ onto $(w_0,\alpha_0)$} is defined as a
standard $(w_0,\alpha_0)^\gamma$--tree $\cT_0=(T_0,\rk_0)$ such that  
\[T_0=\{\langle (t)_\zeta\rest\alpha_0:\zeta\in w_0\cap\rk_1(t)\rangle:
t=\langle (t)_\zeta:\zeta\in w_1\cap\rk_1(t)\rangle\in T_1\}.\]
The mapping 
\[T_1\ni\langle (t)_\zeta:\zeta\in w_1\cap\rk_1(t)\rangle\longmapsto\langle
(t)_\zeta\rest\alpha_0: \zeta\in w_0\cap\rk_1(t)\rangle\in T_0\]
will be denoted $\proj^{(w_1,\alpha_1)}_{(w_0,\alpha_0)}$ too.
\item We say that $\bar{\cT}=\langle\cT_\alpha:\alpha<\alpha^*\rangle$ is
{\em a legal sequence of $\gamma$--trees\/} if for some increasing
continuous sequence $\bar{w}=\langle w_\alpha:\alpha<\alpha^*\rangle$ of
subsets of $\gamma$ we have
\begin{enumerate}
\item[(i)]  $\cT_\alpha$ is a standard $(w_\alpha,\alpha)^\gamma$--tree (for
$\alpha<\alpha^*$), 
\item[(ii)] if $\alpha<\beta<\alpha^*$, then $\cT_\alpha=\proj^{(w_\beta,
\beta)}_{(w_\alpha,\alpha)}(\cT_\beta)$. 
\end{enumerate}
\item Suppose that $\bar{\cT}=\langle\cT_\alpha:\alpha<\alpha^*\rangle$ is a
legal sequence of $\gamma$--trees and $\alpha^*$ is a limit ordinal. Let
$w_\alpha\subseteq\gamma$ be such that $\cT_\alpha$ is a standard
$(w_\alpha,\alpha)^\gamma$--tree (for $\alpha<\alpha^*$) and let
$w=\bigcup\limits_{\alpha<\alpha^*}w_\alpha$. {\em The inverse limit
$\lim\limits^{\leftarrow}(\bar{\cT})$ of $\bar{\cT}$} is a standard
$(w,\alpha^*)^\gamma$--tree $(T^{\lim},\rk^{\lim})$ such that  
\begin{enumerate}
\item[$(\otimes)$] $T^{\lim}$ consists of all sequences $t$ satisfying
\begin{enumerate}
\item[(i)]   $\Dom(t)$ is an initial segment of $w$ (not necessarily
proper),  
\item[(ii)]  if $\zeta\in\Dom(t)$, then $(t)_\zeta$ is a sequence of length 
$\alpha^*$, 
\item[(iii)] $\langle (t)_\zeta\rest\alpha:\zeta\in w_\alpha\cap\Dom(t)
\rangle\in T_\alpha$ for each $\alpha<\alpha^*$.
\end{enumerate}
\end{enumerate}
\item A legal sequence $\bar{\cT}=\langle\cT_\alpha:\alpha<\alpha^*\rangle$
is {\em continuous\/} if for each limit ordinal $\beta<\alpha^*$, 
$\cT_\beta=\lim\limits^{\leftarrow}(\bar{\cT}\rest\beta)$.
\item Let $\bar{\bbQ}=\langle\bbP_i,\name{\bbQ}_i:i<\gamma\rangle$ be a
$\lambda$--support iteration. {\em A standard tree of conditions in
$\bar{\bbQ}$} is a system $\bar{p}=\langle p_t:t\in T\rangle$ such that 
\begin{itemize}
\item $(T,\rk)$ is a standard $(w,\alpha)^\gamma$--tree for some $w\subseteq
\gamma$ and an ordinal $\alpha$, 
\item $p_t\in\bbP_{\rk(t)}$ for $t\in T$, and
\item if $s,t\in T$, $s\vtl t$, then $p_s=p_t\rest\rk(s)$. 
\end{itemize}
\item Let $\bar{p}^0,\bar{p}^1$ be standard trees of conditions in
$\bar{\bbQ}$, $\bar{p}^i=\langle p^i_t:t\in T_i\rangle$, where $\cT_0=
\proj^{(w_1,\alpha_1)}_{(w_0,\alpha_0)}(\cT_1)$, $w_0\subseteq w_1\subseteq
\gamma$, $\alpha_0<\alpha_1$. We will write $\bar{p}^0\leq^{w_1,
\alpha_1}_{w_0,\alpha_0} \bar{p}^1$ (or just $\bar{p}^0\leq \bar{p}^1$)
whenever for each $t\in T_1$, letting $t'=\proj^{(w_1,\alpha_1)}_{(w_0,
\alpha_0)}(t)\in T_0$, we have $p^0_{t'}\rest\rk_1(t)\leq p^1_t$. 
\end{enumerate}
\end{definition}

\begin{remark}
Concerning Definition \ref{dA.5}(4), note that even though $T^{\lim}$ could
be empty, it does satisfy requirements of \ref{dA.5}(1) (so
$\lim\limits^{\leftarrow}(\bar{\cT})$ is indeed a standard
$(w,\alpha^*)^\gamma$--tree). Also, if the sequence $\bar{\cT}$ is
continuous (and $T_\alpha$'s are not empty), then $T^{\lim}\neq\emptyset$. 
\end{remark}

\begin{proposition}
\label{pA.7}
Assume that $\bar{\bbQ}=\langle\bbP_i,\name{\bbQ}_i:i<\gamma\rangle$ is a
$\lambda$--support iteration such that for all $i<\gamma$ we have 
\[\forces_{\bbP_i}\mbox{`` $\name{\bbQ}_i$ is strategically
$({<}\lambda)$--complete ''.}\]
Suppose that $\bar{p}=\langle p_t:t\in T\rangle$ is a standard tree of
conditions in $\bar{\bbQ}$, $|T|<\lambda$, and $\cI\subseteq\bbP_\gamma$
is open dense. Then there is a standard tree of conditions $\bar{q}=\langle
q_t:t\in T\rangle$ such that $\bar{p}\leq \bar{q}$ and $(\forall t\in T)(
\rk(t)=\gamma\ \Rightarrow\ q_t\in\cI)$.
\end{proposition}

\begin{proof}
For $\vare\leq\gamma$ and $r\in\bbP_\vare$, let $\st(
\vare,r)$ be a winning strategy of Complete in $\Game_0^\lambda(
\bbP_\vare,\emptyset,r)$ as in \ref{pA.6}(b). Let 
\[T^{\max}\stackrel{\rm def}{=}\{t\in T:\neg (\exists t'\in T)(t\vtl t')\}
=\{t_\zeta:\zeta<\kappa\}\]  
(where $\kappa<\lambda$ is a cardinal). We construct partial plays $\langle( 
p^\zeta_i,q^\zeta_i):i\leq\kappa\rangle$ of $\Game_0^\lambda(
\bbP_{\rk(p_{t_\zeta})},\emptyset,p_{t_\zeta})$ (for $\zeta<\kappa$) in
which Complete uses strategy $\st(\rk(p_{t_\zeta}),p_{t_\zeta})$ and such
that   
\begin{enumerate}
\item[$(\alpha)$]  if $\zeta<\kappa$ and $\rk(p_{t_\zeta})=\gamma$, then
$p^\zeta_\zeta\in\cI$, 
\item[$(\beta)$]   if $t\vtl t_\zeta$, $t\vtl t_\xi$, $t\in T$, $\zeta,\xi
<\kappa$, $i\leq\kappa$,\\ 
then $p^\zeta_i\rest\rk(t)=p^\xi_i\rest \rk(t)$ and $q^\zeta_i\rest\rk(t)=
q^\xi_i\rest \rk(t)$. 
\end{enumerate}
So suppose we have defined $p^\zeta_j,q^\zeta_j$ for $\zeta<\kappa$, $j<i<
\kappa$. First we look at $\langle (p^i_j,q^i_j):j<i\rangle$ -- it is a play
of $\Game_0^\lambda(\bbP_{\rk(p_{t_i})},\emptyset,p_{t_i})$ in which
Complete uses $\st(\rk(p_{t_i}),p_{t_i})$, so we may find a condition
$p^i_i\in\bbP_{\rk(p_{t_i})}$ stronger than all $p^i_j,q^i_j$ for $j<i$, and
such that $\rk(p_{t_i})=\gamma\ \Rightarrow\ p^i_i\in\cI$. Next, for
$\zeta<\kappa$, $\zeta\neq i$, we define $p^\zeta_i$ as follows. Let $t\in
T$ be such that $t\vtl t_\zeta$, $t\vtl t_i$ and $\rk(t)$ is the largest
possible. We declare that  
\[\Dom(p^\zeta_i)=(\Dom(p^i_i)\cap\rk(t))\cup\bigcup\limits_{j<i}\Dom(
q_j^\zeta)\cup\Dom(p_{t_\zeta}),\]
and $p^\zeta_i\rest\rk(t)=p^i_i\rest\rk(t)$, and for $\vare\in
[\rk(t),\gamma)$ we have that $p^\zeta_i(\vare)$ is the
$<^*_\chi$--first $\bbP_\vare$--name for a member of
$\name{\bbQ}_\vare$ such that 
\[p^\zeta_i\rest\vare\forces_{\bbP_\vare}\mbox{`` $p^\zeta_i(
\vare)$ is an upper bound to }\{p_{t_\zeta}(\vare)\}\cup
\{q^\zeta_j(\vare):j<i\}\mbox{ ''.}\]
The definition of $p^\zeta_i$'s is correct by \ref{pA.6}(b)(iii+ii). Also,
by the choice of ``the $<^*_\chi$--first'' names and clause $(\beta)$ at 
earlier stages we get clause $(\beta)$ for $p^\zeta_i$'s. 

Finally we define $q^\zeta_i$ (for $\zeta<\kappa$) as the condition given 
to Complete by $\st(\rk(t_\zeta),p_{t_\zeta})$ in answer to
$\langle(p^\zeta_j,q^\zeta_j):j<i\rangle\conc\langle p^\zeta_i\rangle$. 
(Again, one easily verifies $(\alpha),(\beta)$.) 

The conditions $p^\zeta_\kappa,q^\zeta_\kappa$ are chosen in a similar
manner except that we do not have to worry about entering $\cI$ anymore, so
we may take $p^0_\kappa$ to be any bound to the previously defined
conditions $p^0_i,q^0_i$, and other $p^\zeta_\kappa,q^\zeta_\kappa$ are
defined as earlier. 

After the above construction is carried out, for $t\in T$ we let 
\begin{quote}
$q_t=p^\zeta_\kappa\rest\rk(t)$ for some (equivalently: all) $\zeta<\kappa$ 
such that $t\trianglelefteq t_\zeta$.
\end{quote}
It should be clear that $\bar{q}=\langle q_t:t\in T\rangle$ is as required.
\end{proof}

Let us close this section with recalling an important result on
easy ensuring that $\lambda$---support iteration satisfies the
$\lambda^{++}$--cc. Its proof is a fairly straightforward modification of
the proof of the respective result for CS iterations; see \cite[Ch. III,
Thm 4.1]{Sh:f}, Abraham \cite[\S 2]{Ab} for the CS case, Eisworth 
\cite[\S 3]{Ei0x} for the general case of $\lambda$--support iterations.  

\begin{theorem}
\label{lppcc}
Assume $2^\lambda=\lambda^+$, $\lambda^{<\lambda}=\lambda$. Let $\bar{\bbQ}=
\langle\bbP_i,\name{\bbQ}_i:i<\lambda^{++}\rangle$ be $\lambda$--support
iteration such that for all $i<\lambda^{++}$ we have  
\begin{itemize}
\item $\bbP_i$ is $\lambda$--proper,
\item $\forces_{\bbP_i}$ ``$|\name{\bbQ}_i|\leq\lambda^+$''. 
\end{itemize}
Then the limit $\bbP_{\lambda^{++}}$ satisfies the $\lambda^{++}$--cc.
\end{theorem}

\section{Bounding properties}
The results on preservation in CS iterations of properties like the Sacks
property and ${}^\omega\omega$--bounding property were among the earliest in
the theory of proper forcing. Here we introduce relatives of these two
properties for $\lambda$--reals and we show suitable iteration theorems. For 
both properties, the properness is ``built into the property''.    

Recall that $\lambda,\bar{\lambda}$ are assumed to be as specified in
Context \ref{incon}.

\begin{definition}
\label{da2}
Let $\bbP$ be a forcing notion.
\begin{enumerate}
\item For a condition $p\in\bbP$ and an ordinal $i_0<\lambda$ we define a
game $\Gsg(i_0,p,\bbP)$ of two players, {\em the Generic player} and {\em
the Antigeneric player}. A play lasts $\lambda$ moves indexed by ordinals
from the interval $[i_0,\lambda)$, and during it the players construct
a sequence $\langle (s_i,\bar{q}^i,\bar{p}^i):i_0\leq i<\lambda\rangle$ as
follows. At stage $i$ of the play (where $i_0\leq i<\lambda$), first 
Generic chooses $s_i\subseteq {}^{\leq i+1}\lambda$ and a system
$\bar{q}^i=\langle q^i_\eta:\eta\in s_i\cap {}^{ i+1}\lambda
\rangle$ such that  
\begin{enumerate}
\item[$(\alpha)$] $s_i$ is a complete $\lambda$--tree of height $i+1$, and 
\[(\forall\eta\in s_i)(\exists\nu\in s_i)(\eta\trianglelefteq\nu\ \&\
\lh(\nu)=i+1),\]
and $\lh(\mrot(s_i))=i_0$, 
\item[$(\beta)$]  for all $j$ such that $i_0\leq j<i$ we have $s_j=s_i\cap
{}^{\leq j+1}\lambda$,
\item[$(\gamma)$] $q^i_\eta\in\bbP$ for all $\eta\in s_i\cap
{}^{ i+1}\lambda$, and 
\item[$(\delta)$] if $i_0\leq j<i$, $\nu\in s_i\cap {}^{ j+1}
\lambda$ and $\nu\vtl\eta\in s_i\cap {}^{ i+1}\lambda$, then
$p^j_\nu\leq q^i_\eta$ and $p\leq q^i_\eta$, 
\item[$(\vare)$] $|s_i\cap {}^{ i+1}\lambda|<\lambda_i$.
\end{enumerate}
Then Antigeneric answers choosing a system $\bar{p}^i=\langle p^i_\eta:
\eta\in s_i\cap {}^{ i+1}\lambda\rangle$ of conditions in $\bbP$
such that $q^i_\eta\leq p^i_\eta$ for each $\eta\in s_i\cap {}^{
i+1}\lambda$.  

The Generic player wins a play if she always has legal moves (so the play
really lasts $\lambda$ steps) and there are a condition $q\geq p$ and a
$\bbP$--name $\name{\rho}$ such that 
\begin{enumerate}
\item[$(\circledast)$]\qquad $q\forces_{\bbP}\mbox{`` }\name{\rho}\in
\bairel\ \&\ \big(\forall i\in [i_0,\lambda\big))\big(\name{\rho}\rest (i+1)
\in s_i\ \&\ q^i_{\name{\rho}\rest (i+1)}\in\Gamma_{\bbP}\big)\mbox{ ''.}$
\end{enumerate}
\item We say that $\bbP$ has the {\em strong $\bar{\lambda}$--Sacks
property\/} whenever the Generic player has a winning strategy in the game
$\Gsg(i_0,p,\bbP)$ for any $i_0<\lambda$ and $p\in\bbP$.
\item We say that $\bbP$ has the {\em $\bar{\lambda}$--Sacks property\/} if
for every $p\in\bbP$ and a $\bbP$--name $\name{\tau}$ such that
$p\forces\name{\tau}:\lambda\longrightarrow\bV$, there are a condition
$q\geq p$ and a sequence $\langle a_\alpha:\alpha<\lambda\rangle$ such that
$|a_\alpha|<\lambda_\alpha$ (for $\alpha<\lambda$) and $q\forces\mbox{``
}(\forall \alpha<\lambda)(\name{\tau}(\alpha)\in a_\alpha)\mbox{ ''}$.
\end{enumerate}
\end{definition}

\begin{remark}
\label{remSacks}
\begin{enumerate}
\item At a stage $i<\lambda$ of a play of $\Gsg(i_0,p,\bbP)$, the
Antigeneric player may play stronger conditions, and using \ref{obsA.4} we
may require that if $\bar{p}^i=\langle p^i_\eta:\eta\in s_i\cap
{}^{ i+1}\lambda\rangle$ is his move, then the conditions
$p^i_\eta$ are pairwise incompatible.Thus the winning criterion
$(\circledast)$ could be replaced by 
\begin{enumerate}
\item[$(\circledast)^-$]\qquad $q\forces_{\bbP}\mbox{`` }\big(\forall i\in
[i_0,\lambda\big))\big(\exists\eta\in s_i\cap {}^{ i+1}\lambda
\big)\big(q^i_\eta\in\Gamma_{\bbP}\big)\mbox{ ''}$
\end{enumerate}
(thus eliminating the use of $\name{\rho}$). However, the $\lambda$--branch
along which the conditions are from the generic filter will be new (so we
cannot replace the name $\name{\rho}$ by an object $\rho\in\bairel$). 
\item Note that if Generic has a winning strategy in $\Gsg(0,p,\bbP)$,
then she has one in $\Gsg(i_0,p,\bbP)$ for all $i_0<\lambda$. (Remember: the
sequence $\bar{\lambda}$ is increasing.) The reason why we have $i_0$ as a
parameter is a notational convenience.
\item Plainly, if Generic has a winning strategy in $\Gsg(i_0,p,\bbP)$, then
she has one with the following property:
\begin{enumerate}
\item[$(\boxtimes_{\rm nice})$] if $s_i,\bar{q}_i$ are given to Generic as a
move at a stage $i\in [i_0,\lambda)$, then for every $\eta\in s_i\cap
{}^{ i}\lambda$, the set $\{\alpha<\lambda:\eta\conc\langle\alpha
\rangle\in s_i\}$ is an initial segment of $\lambda_i$ and $\eta(j)=0$ for
all $j<i_0$. 
\end{enumerate}
Strategies satisfying the condition $(\boxtimes_{\rm nice})$ will be called
{\em nice}.
\item Easily, if $\bbP$ has the strong $\bar{\lambda}$--Sacks property, then
it is strategically $({<}\lambda)$--complete and has the
$\bar{\lambda}$--Sacks property.
\end{enumerate}
\end{remark}

\begin{theorem}
\label{itSacks}
Suppose that $\bar{\bbQ}=\langle\bbP_\alpha,\name{\bbQ}_\alpha:\alpha<\gamma 
\rangle$ is a $\lambda$--support iteration such that for all
$\alpha<\gamma$: 
\[\forces_{\bbP_\alpha}\mbox{`` $\name{\bbQ}_\alpha$ has the strong
$\bar{\lambda}$--Sacks property ''.}\]
Then:
\begin{enumerate}
\item[(a)] $\bbP_\gamma$ has the $\bar{\lambda}$--Sacks property.
\item[(b)] If $N\prec (\cH(\chi),\in,<^*_\chi)$, $|N|=\lambda$,
${}^{ <\lambda}N\subseteq N$ and $\bar{\lambda},\lambda,p,
\bar{\bbQ},\bbP_\gamma,\ldots\in N$, $p\in \bbP_\gamma$, then there is an
$(N,\bbP_\gamma)$--generic condition $r\in\bbP_\gamma$ stronger than $p$.
\end{enumerate}
\end{theorem}

\begin{proof}
(a)\qquad First note that each $\bbP_\alpha$ is strategically
$({<}\lambda)$--complete (by \ref{pA.6}; remember \ref{remSacks}(4)), so our
assumptions on $\lambda,\bar{\lambda}$ hold in intermediate universes
$\bV^{\bbP_\alpha}$. 

For $\alpha<\gamma$ and $i_0<\lambda$ and a $\bbP_\alpha$--name $\name{q}$
for a condition in $\name{\bbQ}_\alpha$, let $\nst_\alpha(i_0,\name{q})$ be
the $<^*_\chi$--first $\bbP_\alpha$--name for a nice (see \ref{remSacks}(3))
winning strategy of the Generic player in the game $\Gsg(i_0,\name{q},
\name{\bbQ}_\alpha)$. 

Let $\name{\tau}$ be a $\bbP_\gamma$--name for a function from $\lambda$ to
$\bV$, $p\in \bbP_\gamma$. Pick a model $N\prec (\cH(\chi),\in,<^*_\chi)$
such that 
\[\bar{\lambda},\lambda,\name{\tau},p,\bar{\bbQ},\bbP_\gamma,\ldots \in
N,\quad\mbox{ and }\quad |N|=\lambda\quad\mbox{ and }\quad {}^{<\lambda}N
\subseteq N.\]
Note that if $i_0<\lambda$, $\alpha\in N\cap\gamma$, and $\name{q}\in N$ is
a $\bbP_\alpha$--name for a member of $\name{\bbQ}_\alpha$, then
$\nst_\alpha(i_0,\name{q})\in N$. Also, as $\bar{\bbQ}$ is a
$\lambda$--support iteration of $({<}\lambda)$--strategically complete
forcing notions, we may use \ref{pA.6} inside $N$ and for each $\vare\in
N\cap\gamma$ and $r\in\bbP_\vare\cap N$ fix a winning strategy
$\st^*(\vare,r)\in N$ of Complete in the game $\Game_0^\lambda(
\bbP_\vare,\emptyset,r)$ so that conditions (i)--(iii) of
\ref{pA.6}(b) hold true.

Fix a list $\bar{\cI}=\langle\cI_\xi:\xi<\lambda\rangle$ of all open dense
subsets of $\bbP_\gamma$ from $N$ and a one-to-one mapping $\pi:N\cap\gamma
\longrightarrow\lambda$. For $i<\lambda$ let $w_i=\pi^{-1}[i]$ (thus
$\bar{w}=\langle w_i:i<\lambda\rangle$ is an increasing continuous sequence
of subsets of $N\cap\gamma$, each of size $<\lambda$, and $\bigcup\limits_{
i<\lambda} w_i=N\cap\gamma$). 

By induction on $i<\lambda$ we define sequences $\langle\cT_i:i<\lambda
\rangle$ and $\langle\bar{p}^i,\bar{p}^i_*:i<\lambda$ is not a limit ordinal
$\rangle$ such that the following requirements are satisfied.  
\begin{enumerate}
\item[$(\alpha)$] $\langle\cT_i:i<\lambda\rangle$ is a continuous legal
sequence of $\gamma$--trees; $\cT_i\in N$ is a standard
$(w_i,i)^\gamma$--tree, $|T_{i+1}|<\lambda_i$, and $(\forall t\in
T_i)(\exists t'\in T_i)(t\trianglelefteq t'\ \&\ \rk_i(t')=\gamma)$. 
\item[$(\beta)$]  For $i<\lambda$ and $t\in T_i$ such that $\rk_i(t)<\gamma$
let $\psi_i(t)=\{(s)_{\rk_i(t)}:t\vtl s\in T_i\}$.\\
Then (for each $i,t$ as above) $\emptyset\neq \psi_i(t)\subseteq
\prod\limits_{j<i}\lambda_j$ and for each $\eta\in \psi_i(t)$ and $i'<
\pi(\rk_i(t))$ we have $\eta(i')=0$. 
\item[$(\gamma)$] If $\xi\in N\cap\gamma$, $\pi(\xi)<i<j<\lambda$, $t\in
T_j$, $\rk_j(t)=\xi$ and $t'=\proj^{w_j,j}_{w_i,i}(t)\in T_i$ (so
$\rk_i(t')=\xi$), then $\psi_i(t')=\{\eta\rest i:\eta\in\psi_i(t)\}$. 
\item[$(\delta)$] $T_0=\{\langle\rangle\}$, $\bar{p}^0=\langle p^0_{\langle
\rangle}\rangle$, $p^0_{\langle\rangle}=p$, and for $i<\lambda$,
$\bar{p}^{i+1}=\langle p^{i+1}_t:t\in T_{i+1}\rangle$ and $\bar{p}^{i+1}_*=
\langle p^{i+1}_{*,t}:t\in T_{i+1}\rangle$ are standard trees of conditions
in $\bar{\bbQ}$, both belonging to $N$ and such that $\bar{p}^{i+1}_*\leq
\bar{p}^{i+1}$. 
\item[$(\vare)$] If $i<j<\lambda$, then $\bar{p}^{i+1}\leq^{w_{j+1},
j+1}_{w_{i+1},i+1}\bar{p}^{j+1}$. 
\item[$(\zeta)$] If $t_{i+1}\in T_{i+1}$ (for $i<\lambda$) are such that
$\rk_{i+1}(t_{i+1})=\gamma$ and $t_{i+1}=\proj^{w_{j+1},j+1}_{w_{i+1},
i+1}(t_{j+1})$ (for $i<j$), then $\langle p^{i+1}_{*,t_{i+1}},
p^{i+1}_{t_{i+1}}:i<\lambda\rangle$ is a play of the game $\Game_0^\lambda(
\bbP_\gamma,\emptyset,p)$ in which Complete uses the strategy $\st^*(\gamma,
p)$. 
\item[$(\eta)$] If $t\in T_{i+1}$, $\rk_{i+1}(t)=\gamma$, then $p^{i+1}_t\in
\cI_\xi$ for all $\xi\leq i$ and $p^{i+1}_t $ forces a value to
$\name{\tau}(i)$. 
\item[$(\theta)$] Assume that $\xi\in N\cap\gamma$, $\pi(\xi)=i_0\leq i$ and
$t\in T_{i+1}$ is such that $\rk_{i+1}(t)=\xi$. Let, for $j\leq i$,
$t_j=\proj^{w_{i+1},i+1}_{w_j,j}(t)$ and let $\name{r}$ be the
$<^*_\chi$--first $\bbP_\xi$--name for a member of $\name{\bbQ}_\xi$ such
that 
\[\begin{array}{ll}
\forces_{\bbP_\xi}&\mbox{`` if there is a common upper bound to }\{p^j_{t_j}
(\xi):j\leq i_0\mbox{ is non-limit }\},\\
&\ \mbox{ then $\name{r}$ is such an upper bound, else $\name{r}=p(\xi)$
''.}
  \end{array}\]
Furthermore, for $i_0\leq j\leq i$ and $\eta\in\psi_{j+1}(t_{j+1})$, fix
$s^{j+1}_\eta\in T_{j+1}$ such that $\rk_{j+1}(s^{j+1}_\eta)>\xi$ and
$(s^{j+1}_\eta)_\xi=\eta$, $t_{j+1}\vtl s^{j+1}_\eta$, and put
$\name{r}^j_\eta= p^{j+1}_{s^{j+1}_\eta}(\xi)$.\\
{\em Then\/} the condition $p^{i+1}_t$ forces in $\bbP_\xi$ the following: 
\begin{quotation}
there is a partial play $\langle s_j,\bar{q}^j,\bar{r}^j:i_0\leq j\leq
i\rangle$ of the game $\Gsg(i_0,\name{r},\name{\bbQ}_\xi)$ in which the
Generic player uses the strategy $\nst_\xi(i_0,\name{r})$ and, for $i_0\leq
j\leq i$, 
\[s_j\cap {}^{j+1}\lambda=\psi_{j+1}(t_{j+1})\quad\mbox{ and }\quad
\bar{r}^j=\langle\name{r}^j_\eta:\eta\in s_j\cap {}^{j+1}\lambda\rangle.\] 
\end{quotation}
\end{enumerate}
Let us describe how the construction of $\langle\cT_i:i<\lambda\rangle$ and
$\langle\bar{p}^{i+1}:i<\lambda\rangle$ is carried out. We start with
letting $T_0=\{\langle\rangle\}$, $p^0_{\langle\rangle}=p$ (as in
$(\delta)$). Now suppose that we have defined $\cT_j,\bar{p}^j$ for
$j<i<\lambda$ so that clauses $(\alpha)$--$(\theta)$ are satisfied. If $i$
is a limit ordinal, then we let $\cT_i=\lim\limits^{\leftarrow}(\langle
\cT_j:j<i\rangle)\in N$ ($\bar{p}^i,\bar{p}^i_*$ are not defined). It is
straightforward to verify conditions $(\alpha)$--$(\gamma)$ (use the
inductive hypothesis), clauses $(\delta)$--$(\theta)$ are not relevant. 

So suppose now that $i$ is a successor ordinal, say $i=i_0+1$. First we let
$\cT^*$ be a the largest standard $(w_i,i)^\gamma$--tree such that
$\proj^{w_i,i}_{w_{i_0},i_0}(\cT^*)=\cT_{i_0}$, and if $t=\langle
(t)_\zeta:\zeta\in w_i\cap \rk^*(t)\rangle\in T^*$, then $(t)_\zeta(i_0)< 
\lambda_{i_0}$, and if $\pi(\zeta)=i_0$, then $(t)_\zeta\rest i_0\equiv
0$. (Plainly $\cT^*\in N$ and $|T^*|<\lambda$.)  Next, for each $t\in T^*$
we define a condition $q_t\in\bbP_{\rk^*(t)}\cap N$ and names
$\name{\alpha}^t(\xi)$ for ordinals (for $\xi\in w_i\cap\rk^*(t)$). For this
let us fix $t\in T^*$ and let $t_j=\proj^{w_i,i}_{w_j,j}(t)\in T_j$ for
$j<i$. Put  
\[\Dom(q_t)=\big(w_i\cup\bigcup\{\Dom(p^j_{t_j}): j<i\mbox{ is not a limit
}\}\big)\cap \rk^*(t),\]
and for $\zeta\in\Dom(q_t)$ let $q_t(\zeta)$ be a $\bbP_\zeta$--name for a
member of $\name{\bbQ}_\zeta$ chosen as follows. If
$\zeta\in\Dom(q_t)\setminus w_i$, then $q_t(\zeta)$ is the $<^*_\chi$--first
$\bbP_\zeta$--name such that
\[\forces_{\bbP_\zeta}\mbox{`` if possible, then $q_t(\zeta)$ is an upper
bound to }\{p^j_{t_j}(\zeta):j<i\mbox{ is non-limit }\}\mbox{ ''}.\]
If $\zeta\in\Dom(q_t)\cap w_i$, then $\name{\alpha}^t(\zeta)\in N$ is a
$\bbP_\zeta$--name for an element of $\lambda_i$ and $q_t(\zeta)$ is the
$<^*_\chi$--first $\bbP_\zeta$--name for a condition in $\name{\bbQ}_\zeta$
with the following property.

Let $\name{r}$ be the $<^*_\chi$--first $\bbP_\zeta$--name for a member of
$\name{\bbQ}_\zeta$ such that 
\[\begin{array}{ll}
\forces_{\bbP_\zeta}&\mbox{`` if possible, then $\name{r}$ is an upper
bound to }\{p^j_{t_j}(\zeta):j\leq\pi(\zeta)\mbox{ is non-limit }\},\\
&\ \ \mbox{ else }\name{r}=p(\zeta)\mbox{ ''.}
  \end{array}\]
Now, suppose that $G_\zeta\subseteq\bbP_\zeta$ is a generic filter over
$\bV$ and $p^{j+1}_{t_{j+1}}\rest\zeta\in G_\zeta$ for all $j<i_0$, and work
in $\bV[G_\zeta]$. Then, by clause $(\theta)$, there is a partial play $\langle
s_j,\bar{q}^j,\bar{r}^j:\pi(\zeta)\leq j< i_0\rangle$ of the game $\Gsg(
\pi(\zeta),\name{r}^G_\zeta,\name{\bbQ}^{G_\zeta}_\zeta)$ in which Generic
uses $\nst_\zeta(\pi(\zeta),\name{r})^{G_\zeta}$, and $s_j\cap{}^{j+1}
\lambda=\psi_{j+1}(t_{j+1}\rest\zeta)$ and $\bar{r}^j=\langle r^j_\eta:\eta
\in s_j\cap{}^{j+1}\lambda\rangle$, where $r^j_\eta=\big(p^{j+1}_{s^{j+
1}_\eta}(\zeta)\big)^{G_\zeta}$ for $s^{j+1}_\eta\in T_{j+1}$ such that
$t_{j+1}\rest \zeta\vtl s^{j+1}_\eta$, $(s^{j+1}_\eta)_\zeta=\eta$ and
$\rk_{j+1}(s^{j+1})=\gamma$. So we may look at the answer $s_{i_0}$,
$\bar{q}^{i_0}=\langle q^{i_0}_\nu:\nu\in s_{i_0}\cap {}^{i_0+1}\lambda
\rangle$ to this play according to the strategy $\nst_\zeta(\pi(\zeta),
\name{r})^{G_\zeta}$. Then, $q_t(\zeta)^{G_\zeta}$ is a condition stronger
than all $r^j_{(t_{j+1})_\zeta}$ for $j<i_0$, and such that 
\begin{quotation}
if $(t_i)_\zeta\in s_{i_0}$, then $q_t(\zeta)^{G_\zeta}=q^{i_0}_{
(t_i)_\zeta}$.
\end{quotation}
Also, $\name{\alpha}^t(\zeta)=\{\alpha<\lambda_i:(t_{i_0})_\zeta\conc
\langle\alpha\rangle\in s_{i_0}\}$.\\
{[If $\pi(\zeta)=i_0$, then we do not have the partial play we started with
- the game just begins and we look at the first move of Generic, requiring
that $q_t(\zeta)^{G_\zeta}$ is stronger than $\name{r}^{G_\zeta}$ and if
$(t_i)_\zeta\in s_{i_0}$ then $q_t(\zeta)^{G_\zeta}=q_{(t_i)_\zeta}^{
i_0}$.]} 

This finishes the definition of $\bar{q}=\langle q_t:t\in T^*\rangle$. One
easily checks that $\bar{q}\in N$ is a tree of conditions (remember the
choice of ``the $<^*_\chi$--first names'') . Also, by induction on $\zeta
\in\Dom(q_t)$, one verifies that $\bar{p}^j\leq^{w_i,i}_{w_j,j}\bar{q}$ for
all non-limit $j\leq i_0$. (Note that if $\pi(\zeta)=i_0$, $t\in T^*$, and
$\rk^*(t)>\zeta$, then in the inductive process we know that by clause
$(\zeta)$  
\[q_t\rest\zeta\forces_{\bbP_\zeta}\mbox{`` there is a common upper bound to
}\{p^j_{t_j}(\zeta):j\leq\pi(\zeta)\mbox{ is non-limit }\}\mbox{ ''}\]
and thus $q_t\rest\zeta$ forces that the respective condition $\name{r}$ is
stronger than all $p^j_{t_j}(\zeta)$ (for non-limit $j\leq\pi(\zeta)$).)

Next, we use \ref{pA.7} to pick a standard tree of conditions $\bar{p}^*=
\langle p^*_t:t\in T^*\rangle\in N$ such that $\bar{q}\leq\bar{p}^*$ and for
each $t\in T^*$ with $\rk^*(t)=\gamma$ the condition $p^*_t$ decides the
values of all names $\name{\alpha}^{t'}(\zeta)$ for $t'\in T^*$, $\zeta\in
w_i\cap\rk(t')$ and the value of $\name{\tau}(i_0)$ (and let $p^*_t\forces
\mbox{`` }\name{\tau}(i_0)=\tau^t_{i_0}\mbox{ ''}$), and such that
$p^*_t\in\cI_\xi$ for all $\xi\leq i_0$. For $t\in T^*$ with $\rk^*(t)=
\gamma$ and for $\zeta\in w_i$ let $\alpha^t(\zeta)$ be the value forced to 
$\name{\alpha}^t(\zeta)$ by $p^*_t$. Since $\name{\alpha}^t(\zeta)$ is a
$\bbP_\zeta$--name, we have that 
\[t_0\vtl t_1\in T^*\ \&\ \rk^*(t_1)=\gamma\ \&\ \zeta\in w_i\cap\rk^*(t_0)\
\qquad\Rightarrow\qquad\ p^*_{t_0}\forces\name{\alpha}^{t_0}(\zeta)=
\alpha^{t_1}(\zeta).\] 
So we may naturally define $\alpha^t(\zeta)$ also for $t\in T^*$ with
$\rk(t)<\gamma$. Now we let 
\[T_i=T_{i_0+1}=\{t\in T^*: (\forall\zeta\in w_i\cap\rk^*(t))((t)_\zeta
(i_0)<\alpha^t(\zeta)\}\]
and $p^i_{*,t}=p^*_t$ for $t\in T_i$ (thus defining $\bar{p}^i_*$). Plainly,
$T_i\in N$ is a standard $(w_i,i)^\gamma$--tree satisfying
$(\alpha)$--$(\gamma)$, $\bar{p}^i_*\in N$. Finally, using the properties of
the strategies $\st^*$ stated in \ref{pA.6}(b) (and the clause $(\zeta)$
from earlier stages) we may pick a standard tree of conditions
$\bar{p}^i=\langle p^i_t:t\in T_i\rangle$ such that $\bar{p}^*\leq
\bar{p}^i$ and  
\begin{enumerate}
\item[if]   $t\in T_i$, $\rk_i(t)=\gamma$, $ t_j=\proj^{w_i,i}_{w_j,j}(t)$
for non-limit $j\leq i$,
\item[then] $\langle p^{j+1}_{t_{j+1}}:j<i\rangle$ is a sequence of answers
of Complete in some partial play of $\Game^\lambda_0(\bbP_\gamma,\emptyset,
p)$ in which she uses the winning strategy $\st^*(\gamma,p)$.
\end{enumerate}
Now one easily verifies that $\cT_i,\bar{p}^i,\bar{p}^i_*$ satisfy
requirements $(\alpha)$--$(\theta)$, thus the construction is complete. 

Let $\cT_\lambda=\lim\limits^{\leftarrow}(\langle\cT_j:j<\lambda\rangle)$. 
We will consider this standard $(N\cap\gamma,\lambda)^\gamma$--tree in
universes $\bV^{\bbP_\xi}$ (for $\xi\leq\gamma$), so let us note that
forcings $\bbP_\xi$ may add new branches in $\cT_\lambda$. But if (in 
$\bV^{\bbP_\xi}$) $t\in\cT_\lambda$ and $i<\lambda$, then 
\[t|i\stackrel{\rm def}{=}\langle (t)_\zeta\rest i: \zeta\in w_i\cap
\rk_\lambda(t)\rangle=\proj^{N\cap\gamma,\lambda}_{w_i,i}(t)\in \bV.\]
Also if $i<\lambda$ is limit, then the equality $\cT_i=
\lim\limits^{\leftarrow}(\langle\cT_j:j<i\rangle)$ holds in $\bV^{\bbP_\xi}$
as well. 

We are going to define a condition $r\in\bbP_\gamma$ such that
$\Dom(r)=N\cap\gamma$ and the names $r(\alpha)$ are defined by induction on
$\alpha\in N\cap\gamma$. For $\alpha\in N\cap\gamma$ we will also choose
$\bbP_{\alpha+1}$--names $\name{t}_\alpha$ for functions in $\bairel$, and
we will put $\name{t}^\alpha=\langle\name{t}_\beta:\beta<\alpha\ \&\
\beta\in N\rangle$. The construction will be carried out so that (for each
$\alpha\in N\cap (\gamma+1)$):
\begin{enumerate}
\item[(i)$_\alpha$]  $r\rest \alpha\forces_{\bbP_\alpha}\mbox{`` }
\name{t}^\alpha\in\cT_\lambda\mbox{ ''}$, 
\item[(ii)$_\alpha$] $r\rest \alpha\forces_{\bbP_\alpha}\mbox{`` } (\forall
i<\lambda)\big((p^{i+1}_{\name{t}^\alpha|i+1})\rest \alpha\in\Gamma_{
\bbP_\alpha}\big)\mbox{ ''}$. 
\end{enumerate}
\smallskip

Arriving at a limit stage $\alpha\in N\cap (\gamma+1)$, we have defined
$r\rest \alpha$ and $\name{t}^\alpha$, and we should only check that
conditions (i)$_\alpha$, (ii)$_\alpha$ hold (assuming (i)$_\beta$,
(ii)$_\beta$ hold for $\beta<\alpha$, $\beta\in N$).

\noindent RE: (i)$_\alpha$:\qquad $\cT_\lambda$ is a standard tree, so every
chain in $T_\lambda$ has a $\vartriangleleft$--bound. Now, the first
condition follows immediately from the inductive hypothesis. 

\noindent RE: (ii)$_\alpha$:\qquad Suppose that $G_\alpha\subseteq
\bbP_\alpha$ is generic over $\bV$ and $r\rest \alpha\in G_\alpha$. Let
$i<\lambda$ and for $\beta\leq\alpha$ let $t^\beta_i=(\name{t}^\beta |
i)^{G_\alpha\cap\bbP_\beta}\in T_i$. Then, by (ii)$_\beta$, we know that
$p^{i+1}_{t^\beta_{i+1}}\rest\beta\in G_\alpha \cap\bbP_\beta$ (for each
$\beta\in\alpha\cap N$). But $p^{i+1}_{t^\alpha_{i+1}}\rest\beta=
p^{i+1}_{t^\beta_{i+1}}\rest\beta$ (as $t^\beta_{i+1}\vtl t^\alpha_{i+1}$),
so remembering that $p^{i+1}_{t^\alpha_{i+1}}\in N$ we conclude
$p^{i+1}_{t^\alpha_{i+1}}\rest \alpha\in G_\alpha$. 
\smallskip

Now suppose that we arrived at stage $\alpha+1\in N\cap (\gamma+1)$ and we
have defined $r\rest\alpha$, $\name{t}^\alpha$ so that (i)$_\alpha$+ 
(ii)$_\alpha$ hold. Let $G_\alpha\subseteq\bbP_\alpha$ be generic over
$\bV$, $r\rest\alpha\in G_\alpha$. For $i<\lambda$ let $t^\alpha_i=(
\name{t}^\alpha | i)^{G_\alpha}\in T_i$ (remember (i)$_\alpha$). Plainly,
$t^\alpha_j=\proj^{w_i,i}_{w_j,j}(t^\alpha_i)$ for $j<i<\lambda$. By
(ii)$_\alpha$ we get $p^{i+1}_{t^\alpha_{i+1}}\rest\alpha\in G_\alpha$ for
all $i<\lambda$. 
\begin{enumerate}
\item[$(\boxplus)_\alpha$] Let $i_0=\pi(\alpha)$ and let $\name{r}$ be the
$<^*_\chi$--first $\bbP_\alpha$--name for an element of $\name{\bbQ}_\alpha$
such that ($\name{r}\in\bV$, of course, and) 
\[\begin{array}{ll}
\forces_{\bbP_\alpha}&\mbox{`` if there is a common upper bound to }
\{p^j_{t^\alpha_j}(\alpha):j\leq i_0\mbox{ is non-limit }\}\\
&\ \mbox{ then $\name{r}$ is such an upper bound, else $\name{r}=p(\alpha)$
''.} 
  \end{array}\]
\end{enumerate}
(Note: for each $j^*<\lambda$ the sequence $\langle t^\alpha_j:j<j^*\rangle$
belongs to the ground model $\bV$, and even to $N$.) 

Fix $j^*<\lambda$, $j^*>i_0$ for a moment. In $\bV$, for each $i\in [i_0,
j^*]$ and $\eta\in \psi_{i+1}(t^\alpha_{i+1})$ let us choose
$s^{i+1}_\eta\in T_{i+1}$ such that $t^\alpha_{i+1}\vtl s^{i+1}_\eta$,
$(s^{i+1}_\eta)_\alpha=\eta$. Now work in $\bV[G_\alpha]$. Since
$p^{j^*+1}_{t^\alpha_{j^*+1}}\in G_\alpha$, we may use clause $(\theta)$ of
the construction to claim that there is a partial play $\bar{\sigma}^{j^*}=
\langle s_i,\bar{q}^i,\bar{r}^i: i_0\leq i\leq j^*\rangle$ of
the game $\Gsg(i_0,\name{r}^{G_\alpha},(\name{\bbQ}_\alpha)^{G_\alpha})$ in 
which Generic uses $\nst_\alpha(i_0,\name{r}^{G_\alpha})$ and
$s_i\cap{}^{i+1}\lambda=\psi_{i+1}(t^\alpha_{i+1})$ and $\bar{r}^i=\langle
(p^{i+1}_{s^{i+1}_\eta}(\alpha))^{G_\alpha}: \eta\in s_i\cap {}^{i+1}\lambda
\rangle$. 

It should be clear that (in $\bV[G_\alpha]$) $\bar{\sigma}^{j^*}\vtl
\bar{\sigma}^{j^{**}}$ for $i_0<j^*<j^{**}<\lambda$, so we have a play
$\bar{\sigma}=\bigcup\limits_{i_0<j^*<\lambda} \bar{\sigma}^{j^*}=\langle
s_i,\bar{q}^i,\bar{r}^i: i_0\leq i<\lambda\rangle$ of the game $\Gsg(i_0, 
\name{r}^{G_\alpha},(\name{\bbQ}_\alpha)^{G_\alpha})$ with the respective
properties. This play is won by Generic, so there are a condition $q\in
(\name{\bbQ}_\alpha)^{G_\alpha}$ and a $(\name{\bbQ}_\alpha)^{
G_\alpha}$--name $\name{\rho}$ for a member of $\bairel$ such that $q\geq
\name{r}^{G_\alpha}$ and 
\[(\otimes)\quad q\forces_{(\name{\bbQ}_\alpha)^{G_\alpha}}\mbox{`` }
(\forall i\in [i_0,\lambda))(\name{\rho}\rest (i+1)\in\psi_{i+1}( 
t^\alpha_{i+1})\ \&\ p^{i+1}_{s^{i+1}_{\name{\rho}\rest (i+1)}}
(\alpha)^{G_\alpha}\in\Gamma_{(\name{\bbQ}_\alpha)^{G_\alpha}})\mbox{ ''.}\]
Let $r(\alpha),\name{t}_\alpha$ be names for the $q,\name{\rho}$ as above
(i.e., $r(\alpha)$ is a $\bbP_\alpha$--name of a member of
$\name{\bbQ}_\alpha$ and $\name{t}_\alpha$ is a $\bbP_{\alpha+1}$--name of a
member of $\bairel$ and $r\rest\alpha$ forces that they have the property
stated in $(\otimes)$). It follows from our choices that (i)$_{\alpha+1}$ +
(ii)$_{\alpha+1}$ hold true, finishing the inductive construction of
$r\in\bbP_\gamma$ and $\name{t}_\alpha$'s. 
\medskip

For $\alpha<\lambda$ let $a_\alpha=\{\tau^t_\alpha:t\in T_{\alpha+1}\ \&\
\rk_{\alpha+1}(t)=\gamma\}$ (remember: $\tau^t_\alpha$ is the value forced
to $\name{\tau}(\alpha)$ by $p^{\alpha+1}_t$). Plainly, $|a_\alpha|<
\lambda_\alpha$ for each $\alpha<\lambda$. 

The proof of the iteration theorem will be complete once we show the
following claim. 

\begin{claim}
\label{cl3}
The condition $r\in\bbP_\gamma$ (defined earlier) is stronger than $p$, it
is $(N,\bbP_\gamma)$--generic and $r\forces_{\bbP_\gamma}\mbox{`` }(\forall
\alpha<\lambda)(\name{\tau}(\alpha)\in a_\alpha)\mbox{ ''}$. 
\end{claim}

\begin{proof}[Proof of the Claim]
First, by induction on $\alpha\in N\cap (\gamma+1)$ we are showing that
$p\rest\alpha\leq r\rest\alpha$. There is nothing to do at limit stages, so
let us deal with non-limit ones. Assume we have shown $p\rest\alpha\leq
r\rest\alpha$. 

Suppose that $G_\alpha\subseteq\bbP_\alpha$ is generic over $\bV$,
$r\rest\alpha\in G_\alpha$. Let $t^\alpha_j=(\name{t}^\alpha | j)^{G_\alpha}
\in T_j$, and let $i_0,\name{r}$ be defined as in $(\boxplus)_\alpha$. 
Since, by (ii)$_\alpha$, $p^j_{t^\alpha_j}\rest\alpha\in G_\alpha$ (for
non-limit $j\leq i_0$) and by the clause $(\zeta)$ of the construction, we
get 
\[\bV[G_\alpha]\models\mbox{`` there is a common upper bound to }\{
p^j_{t^\alpha_j}(\alpha)^{G_\alpha}: j\leq i_0\mbox{ is non-limit }\} 
\mbox{ '',}\]
and thus
\[\bV[G_\alpha]\models\mbox{`` }(\forall j<i_0)(p^{j+1}_{t^\alpha_{j+1}}
(\alpha)^{G_\alpha}\leq \name{r}^{G_\alpha})\mbox{ ''.}\]
By the choice of $r(\alpha)$ we have $r(\alpha)^{G_\alpha}\geq
\name{r}^{G_\alpha}\geq p(\alpha)^{G_\alpha}$. 

Hence $r\rest\alpha\forces p(\alpha)\leq r(\alpha)$, as needed.
\medskip

Now, let $G\subseteq\bbP_\gamma$ be generic over $\bV$, $r\in G$. For
$i<\lambda$ let $t_i=(t^\gamma|i)^G\in T_i$. By (ii)$_\gamma$ we know that
$p^{i+1}_{t_{i+1}}\in G$. By clause $(\eta)$ we have $p^{i+1}_{t_{i+1}}\in
\cI_i$ and (by the definition of $a_i$) $p^{i+1}_{t_{i+1}}\forces
\name{\tau}\in a_i$. The former implies that $G$ intersects $\cI\cap N$ for
each open dense subset of $\bbP_\gamma$ from $N$, the latter gives
$\name{\tau}^G(i)\in a_i$.  
\end{proof}

\noindent (b)\qquad Included in the proof of (a).
\end{proof}

\begin{definition}
\label{dbd}
Let $\bbP$ be a forcing notion.
\begin{enumerate}
\item For a condition $p\in\bbP$ and an ordinal $i_0<\lambda$ we define a
game $\Gbd(i_0,p,\bbP)$ like $\Gsg(i_0,p,\bbP)$, but demand
\ref{da2}($1(\vare)$) is replaced by 
\begin{enumerate}
\item[$(\vare)^-$] \quad $|s_i\cap {}^{i+1}\lambda|<\lambda$.
\end{enumerate}
\item $\bbP$ has {\em the strong $\lambda$--bounding property\/} if Generic
has a winning strategy in the game $\Gbd(i_0,p,\bbP)$ for every
$i_0<\lambda$, $p\in \bbP$.
\item $\bbP$ has {\em the $\lambda$--bounding property\/} if for every
$p\in\bbP$ and a $\bbP$--name $\name{\tau}$ such that $p\forces$``
$\name{\tau}\longrightarrow\bV$ '', there are a condition $q\geq p$ and a
sequence $\langle a_\alpha:\alpha<\lambda\rangle$ such that
$|a_\alpha|<\lambda$ (for $\alpha<\lambda$) and $q\forces$`` $(\forall
\alpha<\lambda)(\name{\tau}(\alpha)\in a_\alpha)$ ''. 
\end{enumerate}
\end{definition}

\begin{remark}
\begin{enumerate}
\item All the remarks stated in \ref{remSacks} have their (obvious)
parallels for the $\lambda$--bounding properties. 
\item Clearly, (strong) $\bar{\lambda}$--Sacks property implies (strong,
respectively) $\lambda$--bounding property. 
\end{enumerate}
\end{remark}

\begin{theorem}
\label{itbound}
Suppose that $\bar{\bbQ}=\langle \bbP_\alpha,\name{\bbQ}_\alpha: \alpha<
\gamma\rangle$ is a $\lambda$--support iteration such that for all
$\alpha<\lambda$: 
\[\forces_{\bbP_\alpha}\mbox{`` $\name{\bbQ}_\alpha$ has the strong
$\lambda$--bounding property ''.}\]
Then:
\begin{enumerate}
\item[(a)] $\bbP_\gamma$ has the $\lambda$--bounding property.
\item[(b)] If $N\prec (\cH(\chi),\in,<^*_\chi)$, $|N|=\lambda$,
${}^{ <\lambda}N\subseteq N$ and $\lambda,p,\bar{\bbQ},\bbP_\gamma,\ldots
\in N$, $p\in \bbP_\gamma$, then there is an $(N,\bbP_\gamma)$--generic
condition $r\in\bbP_\gamma$ stronger than $p$. 
\end{enumerate}
\end{theorem}
 
\begin{proof}
Basically the same as for \ref{itSacks}, just replacing each occurrence of
$\lambda_i$ by $\lambda$. 
\end{proof}

\section{Fuzzy properness over $\lambda$}
A properness-type property preserved in $\lambda$--support iterations, so
called {\em properness over semi-diamonds}, was introduced in Ros{\l}anowski
and Shelah \cite{RoSh:655}. That property worked for any uncountable regular
cardinal $\lambda$ satisfying $\lambda^{<\lambda}=\lambda$ (not necessarily
strongly inaccessible), so because of the known ZFC limitations a number of
natural forcing notions were not covered. For the context considered in this
paper we may do much better: {\em fuzzy properness\/} introduced in this
section captures more examples. Even though we do not prove a real
preservation in $\lambda$--support iterations, our iteration theorem
\ref{1.6} is satisfactory for most applications (see sections B.4 and B.8
later). 

In this section we fix $\lambda^*,A,W$ and $D$ such that 

\begin{context}
\label{Con3.1}
\begin{enumerate}
\item $\lambda^*>\lambda$ is a regular cardinal, $A\subseteq\cH_{<\lambda}( 
\lambda^*)$, $W\subseteq [A]^{\textstyle\lambda}$, and if $a\in W$, $w\in
[a]^{\textstyle {<}\lambda}$, $f:w\longrightarrow a$, then $f\in a$ (hence
also $0\in a$ for $a\in W$),
\item for every $x\in \cH(\chi)$ there is a model $N\prec (\cH(\chi),\in,
<^*_\chi)$ such that $|N|=\lambda$, ${}^{ <\lambda}N\subseteq N$,
$x\in N$ and $N\cap A\in W$,
\item $D$ is a normal filter on $\lambda$ such that there is a $D$--diamond 
(see \ref{Def3.2}).
\end{enumerate}
\end{context}

\begin{definition}
\label{Def3.2}
\begin{enumerate}
\item We say that $\bar{F}=\langle F_\delta:\delta\in S\rangle$
is {\em a $D$--pre-diamond sequence\/} if 
\begin{itemize}
\item $S\in D^+$ contains all successor ordinals below $\lambda$, $\lambda
\setminus S$ is unbounded in $\lambda$, $0\notin S$, and 
\item $F_\delta\in {}^\delta\delta$ for all $\delta\in S$.
\end{itemize}
\item {\em A convenient $D$--diamond\/} is a $D$--pre-diamond
$\bar{F}=\langle F_\delta:\delta\in S\rangle$ such that
\[(\forall f\in\bairel)(\{\delta\in S: F_\delta\subseteq f\}\in D^+).\] 
\end{enumerate}
\end{definition}

\begin{definition}
\label{base}
Let $\bbP$ be a forcing notion. {\em A $\lambda$--base for $\bbP$ over
$W$} is a pair $(\gR,\bar{\gY})$ such that
\begin{enumerate}
\item[(a)] $\gR\subseteq \bbP\times \lambda\times A$ is a relation such that 
\begin{center}
if $(p,\delta,x)\in \gR$ and $p\leq_{\bbP} p'$, then $(p',\delta,x)\in\gR$,
\end{center}
\item[(b)] $\bar{\gY}=\langle \gY_a:a\in W\rangle$ where, for each $a\in W$,
$\gY_a:\lambda\longrightarrow [a]^{\textstyle{<}\lambda}$,
\item[(c)] if $q\in\bbP$, $a\in W$, and $\delta<\lambda$ is a limit
ordinal,\\ 
then there are $p\geq_{\bbP} q$ and $x\in\gY_a(\delta)$ such that
$(p,\delta,x)\in \gR$. 
\end{enumerate}
If $\gR$ is understood and $(p,\delta,x)\in \gR$, then we may say {\em
$p$ obeys $x$ at $\delta$}.
\end{definition}

\begin{definition}
\label{pre1.2}
Let $\bbP$ be a forcing notion and let $(\gR,\bar{\gY})$ be a
$\lambda$--base for $\bbP$ over $W$. Also let a model $N\prec ({\mathcal
H}(\chi),{\in},{<^*_\chi})$ be such that $|N|=\lambda$, ${}^{<\lambda}N
\subseteq N$, $a\stackrel{\rm def}{=}N\cap A\in W$ and $\{\lambda,\bbP,D,
\gR\}\in N$. Furthermore, let $h:\lambda\longrightarrow N$ be such that the 
range $\Rang(h)$ of the function $h$ includes $\bbP\cap N$ and let
$\bar{F}=\langle F_\delta:\delta\in S\rangle$ be a $D$--pre-diamond
sequence.  
\begin{enumerate}
\item Let $\bcI=\langle\cI_\alpha:\alpha<\lambda\rangle\subseteq N$ list all
open dense subsets of $\bbP$ from $N$. A sequence $\bar{p}=\langle p_\alpha:
\alpha<\delta\rangle$ of conditions from $\bbP\cap N$ of length $\delta\leq
\lambda$ is called {\em $\bcI$--exact\/} if 
\[(\forall\xi<\delta)(\exists\alpha<\delta)(p_\alpha\in\cI_\xi).\]
\item We say that $\bar{F}$ is {\em a quasi $D$--diamond sequence for
$(N,h,\bbP)$} if for some (equivalently: all) list $\bcI=\langle\cI_\alpha:
\alpha<\lambda\rangle$ of all open dense subsets of $\bbP$ from $N$, for
every $\leq_\bbP$--increasing sequence $\bar{p}=\langle p_\alpha:\alpha<
\lambda\rangle\subseteq \bbP\cap N$ such that   
\[E\stackrel{\rm def}{=}\{\delta<\lambda:\langle p_\alpha:\alpha<\delta
\rangle\mbox{ is $\bcI$--exact }\}\in D\]
(equivalently: $\bar{p}$ is $\bcI$--exact) we have 
\[\{\delta\in E: (\forall\alpha<\delta)(h\circ F_\delta(\alpha)=p_\alpha)\}
\in D^+.\]

\item For a limit ordinal $\delta\in S$ we define $\cY(\delta)=\cY(N,\bbP,h,
\bar{F},\gR,\bar{\gY},\delta)$ as the set 
\[\begin{array}{ll}
\big\{x\in\gY_a(\delta):&
\mbox{if }\langle h\circ F_\delta(\alpha):\alpha<\delta\rangle\mbox{ is a
$\le_{\bbP}$--increasing sequence}\\
&\mbox{of conditions from $\bbP$,}\\
&\mbox{then there is a condition $p\in\bbP$ such that}\\
&(\forall\alpha<\delta)(h\circ F_\delta(\alpha)\leq_{\bbP} p)\mbox{ and }
(p,\delta,x)\in \gR\big\}
  \end{array}\]
\item Let $\bcI=\langle\cI_\alpha:\alpha<\lambda\rangle\subseteq N$ list all
open dense subsets of $\bbP$ from $N$. A sequence $\bar{q}=\langle
q_{\delta,x}:\delta\in S\mbox{ limit }\ \&\ x\in \cX_\delta\rangle\subseteq
N\cap\bbP$ is called {\em a weak fuzzy candidate over $\bar{F}$ for
$(N,h,\bbP,\gR,\bar{\gY},\bcI)$} whenever $\emptyset\neq
X_\delta\subseteq\cY(\delta)$ (for limit $\delta\in S$) and
\begin{enumerate}
\item[$(\alpha)$] $\{\delta\in S: (\forall x\in\cX_\delta)(q_{\delta,x}\in
\cI_\alpha)\}=S\mod D$ for each $\alpha<\lambda$, and
\item[$(\beta)$]  if $\delta\in S$ is a limit ordinal, $x\in \cX_\delta$,
and $\langle h\circ F_\delta(\alpha):\alpha<\delta\rangle$ is a
$\le_{\bbP}$--increasing $\bcI$--exact sequence of members of $\bbP\cap
N$,\\  
then $(\forall\alpha<\delta)(h\circ F_\delta(\alpha)\leq_{\bbP}
q_{\delta,x})$ and $(q_{\delta,x},\delta,x)\in \gR$.
\end{enumerate}
If above $X_\delta=\cY(\delta)$ for each limit $\delta\in S$, then $\bar{q}$ 
is called {\em a fuzzy candidate over $\bar{F}$ for $(N,h,\bbP,\gR,
\bar{\gY},\bcI)$}. 

Omitting $\bcI$ means ``for some $\bcI$''.
\item Let $\bar{q}=\langle q_{\delta,x}:\delta\in S\mbox{ limit }\ \&\ x\in
\cX_\delta\rangle$ be a weak fuzzy candidate over $\bar{F}$ for $(N,h,
\bbP,\gR, \bar{\gY},\bcI)$, and $r\in\bbP$. We define a game $\Gfl(r,N,
\bcI,h,\bbP,\bar{F},\bar{q})$ of two players, the {\em Generic player} and
the {\em Antigeneric player}, as follows. A play lasts $\lambda$ moves, in
the $i^{\rm th}$ move a condition $r_i\in\bbP$ and a set $C_i\in D$ are
chosen such that $(\forall j<i)(r\leq r_j\leq r_i)$, and 
Generic chooses $r_i,C_i$ if $i\in S=\dom(\bar{F})$, and Antigeneric chooses
$r_i,C_i$ if $i\notin S$. In the end Generic wins the play if
\begin{enumerate}
\item[$(\alpha)$] $(\forall\alpha<\lambda)(\exists i<\lambda)(\exists p\in
\bbP\cap N)(p\in\cI_\alpha\ \&\ p\leq r_i)$, and 
\item[$(\beta)$]  if $\delta\in S\cap\bigcap\limits_{i<\delta} C_i$ is a
limit ordinal, $\langle h\circ F_\delta(\alpha):\alpha<\delta\rangle$ is
a $\leq_{\bbP}$--increasing $\bcI$--exact sequence and $(\forall\alpha<
\delta)(\exists i<\delta)(h\circ F_\delta(\alpha)\leq r_i)$,\\ 
then for some $x\in\cX_\delta$ we have $q_{\delta,x}\leq r_\delta$.  
\end{enumerate}
\item Let $\bar{q}$ be a weak fuzzy candidate over $\bar{F}$ for
$(N,h,\bbP,\gR,\bar{\gY},\bcI)$. We say that a condition $r\in{\mathbb
P}$ is {\em $(\gR,\bar{\gY})$--fuzzy generic for $\bar{q}$ (over
$(N,\bcI,h,\bbP,\bar{F})$)} if Generic has a winning strategy in the game
$\Gfl(r,N,\bcI,h,\bbP,\bar{F},\bar{q})$.  
\end{enumerate}
\end{definition}

\begin{remark}
\begin{enumerate}
\item For any two lists $\bcI^1,\bcI^2$ of open dense subsets of
$\bbP$ from $N$, on a club $E$ of $\lambda$ we have  
\[\{\cI^1_\xi:\xi<\delta\}=\{\cI^2_\xi:\xi<\delta\}\]
for $\delta\in E$. Thus the corresponding notions of exactness agree for
$\delta\in E$. As the generic player can choose $C_i\subseteq E$, in
\ref{pre1.2}(4,5,6) we allow not mention $\bcI$ as a parameter. 
\item Plainly, every fuzzy candidate is a weak fuzzy candidate.
\end{enumerate}
\end{remark}

\begin{definition}
\label{1.2} 
Let $\bbP$ be a $\lambda$--complete forcing notion.
\begin{enumerate}
\item We say that $\bbP$ is {\em fuzzy proper over quasi $D$--diamonds for
$W$} whenever for some $\lambda$--base $(\gR,\bar{\gY})$ for $\bbP$ over $W$
and for some $c\in\cH(\chi)$,   
\begin{enumerate}
\item[$(\circledast)$\quad {\em if}]
\begin{itemize}
\item $N\prec(\cH(\chi),\in,<^*_\chi)$, $|N|=\lambda$, ${}^{<\lambda}N
\subseteq N$, $\lambda,\bbP,c,\gR\in N$, and $a\stackrel{\rm def}{=}N\cap
A\in  W$, $p\in\bbP\cap N$,  
\item $h:\lambda\longrightarrow N$ satisfies $\bbP\cap N\subseteq\Rang(h)$,
and 
\item $\bar{F}$ is a quasi $D$--diamond for $(N,h,\bbP)$ and $\bar{q}$ is a
fuzzy candidate over $\bar{F}$,  
\end{itemize}
\item[{\em then}] there is $r\in\bbP$ stronger than $p$ and such that $r$
is $(\gR,\bar{\gY})$--fuzzy generic for $\bar{q}$.
\end{enumerate}
(We may call $(\gR,\bar{\gY})$ and $c$ {\em witnesses for fuzzy
properness}.) 
\item $\bbP$ is {\em strongly fuzzy proper over quasi $D$--diamonds}
whenever for some $\lambda$--base $(\gR,\bar{\gY})$ for $\bbP$ over $W$
and for some $c\in\cH(\chi)$,   
\begin{enumerate}
\item[$(\circledast)^+$\quad {\em if}] 
\begin{itemize}
\item $N\prec(\cH(\chi),\in,<^*_\chi)$, $|N|=\lambda$, ${}^{<\lambda}N
\subseteq N$, $\lambda,\bbP,c,\gR\in N$, and $a\stackrel{\rm def}{=}N\cap
A\in  W$, $p\in\bbP\cap N$,  
\item $h:\lambda\longrightarrow N$ satisfies $\bbP\cap N\subseteq\Rang(h)$,
\item $\bar{F}$ is a quasi $D$--diamond for $(N,h,\bbP)$ and $\bar{q}$ is a
weak fuzzy candidate over $\bar{F}$,  
\end{itemize}
\item[{\em then}] there is a condition $r\in\bbP$ stronger than $p$ such
that $r$ is $(\gR,\bar{\gY})$--fuzzy generic for $\bar{q}$. 
\end{enumerate}
\item $\bbP$ is {\em weakly fuzzy proper over quasi $D$--diamonds} whenever  
for some $\lambda$--base $(\gR,\bar{\gY})$ for $\bbP$ over $W$
and for some $c\in\cH(\chi)$,   
\begin{enumerate}
\item[$(\circledast)^-$\quad {\em if}]
\begin{itemize}
\item $N\prec(\cH(\chi),\in,<^*_\chi)$, $|N|=\lambda$, ${}^{<\lambda}N
\subseteq N$, $\lambda,\bbP,c,\gR\in N$, and $a\stackrel{\rm def}{=}N\cap
A\in  W$, $p\in\bbP\cap N$,  
\item $h:\lambda\longrightarrow N$ satisfies $\bbP\cap N\subseteq\Rang(h)$, 
\end{itemize}
\item[{\em then}] for some quasi $D$--diamond $\bar{F}$ for $(N,h,\bbP)$ and
a weak fuzzy candidate $\bar{q}$ over $\bar{F}$, there is a condition
$r\in\bbP$ stronger than $p$ such that $r$ is $(\gR,\bar{\gY})$--fuzzy
generic for $\bar{q}$. 
\end{enumerate}
\item $\bbP$ is {\em fuzzy proper for $W$} if it is fuzzy proper over
quasi $D'$--diamonds for every normal filter $D'$ on $\lambda$ (which has 
diamonds). Similarly for {\em strongly fuzzy} and {\em weakly fuzzy
proper}. 
\end{enumerate}
\end{definition}

\begin{remark}
Strong fuzzy properness is very close to {\em properness over
semi-diamonds\/} of Ros{\l}anowski and Shelah \cite{RoSh:655} and even
closer to {\em properness over diamonds\/} introduced by Eisworth
\cite{Ei0x}. (Note that considering the condition \ref{1.2}$(\circledast)^+$
we may assume that the weak fuzzy candidate $\bar{q}=\langle
q_{\delta,x}:\delta\in S\mbox{ is limit }\&\ x\in \cX_\delta\rangle$ is such
that $|\cX_\delta|=1$ for each relevant $\delta$, so one may treat it as  
$\bar{q}=\langle q_\delta:\delta\in S\mbox{ is limit }\rangle$.) Thus fuzzy
properness has a flavour of a weaker property. However, the differences in
technical details of the conditions introduced in this section and those in
\cite{RoSh:655} and/or \cite{Ei0x} make it unclear if there are any
implications between the ``properness conditions'' in this section and those
in the other two papers.
\end{remark}

\begin{proposition}
\label{1.2A} 
Let $N,\bbP,h,\bcI,\gR,\bar{\gY}$ be as in \ref{pre1.2}, $\bar{F}=\langle
F_\delta:\delta\in S\rangle$ be a $D$--pre-diamond. 
\begin{enumerate}
\item If the forcing notion $\bbP$ is $\lambda$-complete, then there is a
fuzzy candidate $\bar{q}$ over $\bar{F}$ for $(N,h,\bbP,\gR,\bar{\gY},
\bcI)$. In fact we can even demand:  
\begin{enumerate}
\item[$(+)$] for every $\alpha<\lambda$, for every large enough $\delta\in
S$, $q_{\delta,x}\in\cI_\alpha$ for all $x\in\cY(\delta)$.
\end{enumerate}
\item If $r$ is $(\gR,\bar{\gY})$--fuzzy generic for some weak fuzzy
candidate $\bar{q}$, then $r$ is $(N,\bbP)$--generic (in the standard
sense).   
\item  Assume that a condition $r$ is $(N,\bbP)$--generic (in the
standard sense), $\bar{F}$ is a quasi $D$--diamond and $\bar{q}$ is a
weak fuzzy candidate over $(N,\bcI,h,\bbP,\bar{F})$.  Suppose that
Generic has a strategy in the game $\Gfl(r,N,\bcI,h,\bbP,\bar{F},
\bar{q})$ which guarantees that the result $\langle r_i,C_i:i<\lambda
\rangle$ of the play satisfies \ref{pre1.2}(5)$(\beta)$. Then she has
a winning strategy in $\Gfl(r,N,\bcI,h,\bbP,\bar{F},\bar{q})$ (i.e.,
one ensuring $(\alpha)+(\beta)$ of \ref{pre1.2}(5)). 
\item If $\bbP$ is fuzzy proper over quasi $D$--diamonds, then it is weakly
fuzzy proper over quasi $D$--diamonds.  If $\bbP$ is strongly fuzzy proper
over quasi $D$--diamonds, then it is fuzzy proper over quasi $D$--diamonds. 
\item Assume that $\bbP$ is weakly fuzzy proper over quasi $D$--diamonds,
$\mu\geq\lambda$, $Y\subseteq [\mu]^{\le\lambda}$, $A^*\subseteq\cH(\chi)$,
$W^*\subseteq [A^*]^{\textstyle\lambda}$ ($Y,A^*,W^*\in\bV$). Then: 
\begin{enumerate}
\item[(a)] forcing with $\bbP$ does not collapse $\lambda^+$,
\item[(b)] forcing with $\bbP$ preserves the following two properties: 
\begin{enumerate}
\item[(i)]  $Y$ is a cofinal subset of $[\mu]^{\leq\lambda}$ (under
inclusion), 
\item[(ii)] for every $x\in\cH(\chi)$ there is $N\prec (\cH(\chi),\in,
<^*_\chi)$ such that $|N|=\lambda$, ${}^{<\lambda}N\subseteq N$, $N\cap
A^*\in W^*$ (i.e., the stationarity of $W^*$ under the relevant filter).
\end{enumerate}
\end{enumerate}
\end{enumerate}
\end{proposition}

\begin{proof}
(1)\quad Immediate (by the $\lambda$--completeness of $\bbP$; remember
\ref{base}(c) and that $\gR\in N$; note that $\gY_a(\delta)\in N$). 

\smallskip

\noindent (2)\quad Remember that $0 \notin S$, so in the game $\Gfl(r,N,\bcI,
h,\bbP,\bar{F},\bar{q})$ the condition $r_0$ is chosen by Antigeneric. So if
the conclusion fails, then for some $\bbP$--name $\name{\alpha}\in N$ for an
ordinal we have $r\not\forces\mbox{`` }\name{\alpha}\in N\mbox{ ''}$. Thus
Antigeneric can choose $r_0\geq r$  so that $r_0\forces\mbox{`` }
\name{\alpha}=\alpha_0\mbox{ ''}$ for some ordinal $\alpha_0\notin N$, what
guarantees him to win the play (remember clause $(\alpha)$ of
\ref{pre1.2}(5)).   
\smallskip

\noindent (3)\quad The Generic player modifies her original strategy as
follows. During the play she builds aside a $\leq_\bbP$--increasing sequence
of conditions $\langle p_i:i\in \lambda\setminus S\rangle 
\subseteq\bbP\cap N$. Arriving to stage $i+1$, $i\in \lambda\setminus S$,
she has two sequences: $\langle r_j,C_j:j\leq i\rangle$ (of the play) and
$\langle p_j: j\in i\setminus S\rangle$ such that $p_j\leq r_j$. Now Generic
picks $p_i\in \bbP\cap N$ such that 
\[(\forall j\in i\setminus S)(p_j\leq p_i)\quad\mbox{ and }\quad (\forall\xi
<i)(p_i\in\cI_\xi),\]
and $p_i,r_i$ are compatible (remember: we assumed that $r$ is
$(N,\bbP)$--generic). Next she replaces $r_i$ by a common upper bound of
$p_i$ and $r_i$, pretending that that was the condition played by her
opponent, and then she plays according to her original strategy.  One easily
verifies that this is a winning strategy for the Generic player.
\smallskip

\noindent (4)\quad Straightforward (remember that, by \ref{Con3.1}(3), there 
is a quasi $D$--diamond and by \ref{1.2A}(1) there is a fuzzy candidate over
it). 
\smallskip

\noindent (5)\quad Follows from (2).
\end{proof}

\begin{proposition}
\label{lambdaplus}
$\lambda^+$--complete forcing notions are strongly fuzzy proper for $W$.  
\end{proposition}

\begin{proof}
This is essentially a variant of \cite[2.5]{RoSh:655}, but since we did not 
give the proof there, we will present it fully here.

So suppose that a forcing notion $\bbP$ is $\lambda^+$--complete. Let
$\gR^{\rm tr}=\gR^{\rm tr}(\bbP)$ be the trivial relation consisting of all
triples $(p,\delta,0)$ such that $p\in\bbP$ and $\delta<\lambda$ and let
$\bar{\gY}^{\rm tr}$ be such that $\gY^{\rm tr}_a(\delta)=\{0\}$ (for each
$\delta<\lambda$, $a\in W$).  Assume now that
\begin{itemize}
\item $N\prec(\cH(\chi),\in,<^*_\chi)$, $|N|=\lambda$, ${}^{<\lambda}N
\subseteq N$, $\lambda,\bbP\in N$, and $a\stackrel{\rm def}{=}N\cap
A\in  W$,
\item  $p\in\bbP\cap N$, and $h:\lambda\longrightarrow N$ satisfies
$\bbP\cap N\subseteq\Rang(h)$, 
\item $\bar{F}=\langle F_\delta:\delta\in S\rangle$ is a quasi
$D$--diamond for $(N,h,\bbP)$ and $\bar{q}$ is a weak fuzzy candidate over
$\bar{F}$. Since $\gY^{\rm tr}_a(\delta)$ has a one member only we may think 
of $\bar{q}$ as a sequence $\langle q_\delta:\delta\in S\mbox{ is limit}\;
\rangle$. 
\end{itemize}
Let $\bcI=\langle\cI_\xi:\xi<\lambda\rangle$ list all open dense subsets of
$\bbP$ from $N$.

We are going to build a condition $r\in\bbP$ stronger than $p$ which is
$(\gR^{\rm tr},\bar{\gY}^{\rm tr})$--fuzzy generic for $\bar{q}$. For this 
we inductively build a $\leq_{\bbP}$--increasing sequence $\langle r_i':i< 
\lambda\rangle \subseteq\bbP \cap N$ such that 
\begin{itemize}
\item $r_0'=p$, $r_{i+1}'\in \bigcap\limits_{\xi\leq i}\cI_\xi$,
\item if there is an upper bound to $\{r_j':j<i\}\cup\{q_i\}$, then $r_i'$
is such an upper bound.
\end{itemize}
Then we pick any upper bound $r$ to the sequence $\langle r_i':i<\lambda
\rangle$ (remember: $\bbP$ is $\lambda^+$--complete). Now we want to argue
that Generic has a winning strategy in the game $\Gfl(r,N,\bcI,h,\bbP,
\bar{F},\bar{q})$. Since $r$ is $(N,\bbP)$--generic it is enough to give a
strategy for the Generic player which ensures that the result of the play
satisfies \ref{pre1.2}(5)$(\beta)$ (by \ref{1.2A}(3)). To this end note that
there is a club $E_0$ of $\lambda$ such that
\begin{itemize}
\item every member of $E_0$ is a limit of ordinals from $\lambda\setminus
S$, 
\item for every $\delta\in E_0$ and $i<\delta$, 
\[\{q\in\bbP:q\geq r_i'\ \mbox{ or } q,r_i'\mbox{ are incompatible }\}\in
\{\cI_\xi:\xi<\delta\},\]
\end{itemize}
Let Generic play so that arriving to a stage $\delta\in S$ of the play she
puts the $<^*_\chi$--first upper bound to the conditions played so far and
$E_0$. Why does this strategy work? Let $\langle r_i,C_i:i<\lambda\rangle$
be the result of the play in which Generic plays as described above and let 
$\delta\in S\cap\bigcap\limits_{i<\delta}C_i$ be a limit ordinal such that 
$\langle h\circ F_\delta(\alpha):\alpha<\delta\rangle$ is a
$\leq_{\bbP}$--increasing $\bcI$--exact sequence and 
\[(\forall\alpha<\delta)(\exists i<\delta)(h\circ F_\delta(\alpha)\leq
r_i).\] 
Then no $r_i',h\circ F_\delta(\alpha)$ (for $i,\alpha<\delta$) can be
incompatible, so (since $\delta\in E_0$ and $\langle h\circ
F_\delta(\alpha):\alpha<\delta\rangle$ is $\bcI$--exact) we have also
\[(\forall i<\delta)(\exists\alpha<\delta)(r_i'\leq h\circ
F_\delta(\alpha)),\] 
and hence $q_\delta$ is stronger than all $r_i'$ (for $i<\delta$). Therefore
$q_\delta\leq r_\delta'\leq r_\delta$.
\end{proof}

\begin{theorem}
\label{1.6} 
Let $A,W,D$ be as in \ref{Con3.1} and let $\bar{\bbQ}=\langle\bbP_\alpha,
\name{\bbQ}_\alpha:\alpha<\zeta^*\rangle$ be a $\lambda$--support
iteration of $\lambda$--complete forcing notions, and assume that
$\zeta^*\subseteq A$. Suppose also that for each
$\zeta<\zeta^*$ we have $\bar{\gY}^\zeta$ and $\bbP_\zeta$--names
$\name{\gR}_\zeta,\name{c}_\zeta$ such that   
\[\begin{array}{ll}
\forces_{\bbP_\zeta}&\mbox{`` $\name{\bbQ}_\zeta$ is fuzzy proper over
quasi $D$--diamonds for $W$}\\
&\qquad\mbox{ with witnesses $(\name{\gR}_\zeta,\bar{\gY}^\zeta)$ and 
$\name{c}_\zeta$''.}
  \end{array}\] 
Then $\bbP_{\zeta^*}=\lim(\bar{\bbQ})$ is weakly fuzzy proper over quasi
$D$--diamonds.   
\end{theorem}

\begin{proof}  
By \ref{first}, the forcing notion $\bbP_{\zeta^*}$ is $\lambda$--complete,
so we have to concentrate on showing clause \ref{1.2}(3)($(\circledast)^-$)
for it. The proof, though not presented as such, is by induction on
$\zeta^*$. However, the inductive hypothesis is used only to be able to
claim that $A,W,D$ are as in \ref{Con3.1} when considered in the
intermediate universes $\bV^{\bbP_\zeta}$ (for $\zeta<\zeta^*$) ---
remember \ref{1.2A}(5). Thus our assumptions on $\name{\bbQ}_\zeta$'s are
meaningful. 

Let us fix a convenient $D$--diamond sequence $\bar{F}'=\langle F_\delta':
\delta\in S\rangle$ (so in particular, $S\in D^+$ contains all successors,
$\lambda\setminus S$ is unbounded in $\lambda$ and $0\notin S$). Put 
\[E_0\stackrel{\rm def}{=}\{\delta<\lambda:\delta\mbox{ is a limit of points
from $\lambda\setminus S$ }\},\qquad E_1\stackrel{\rm
def}{=}(\lambda\setminus S)\cup E_0.\]  
Plainly, $E_0,E_1$ are clubs of $\lambda$. Let $\langle i_\alpha:\alpha<
\lambda\rangle$ be the increasing enumeration of $E_1$ and $E_2=E_0\cap\{
\alpha<\lambda:i_\alpha=\alpha\}$ (So $E_2$ is a club of $\lambda$ too).

For each $a\in W$ fix a one-to-one mapping $\pi_a:a\cap\zeta^*
\longrightarrow\lambda$ such that $\pi_a(0)=0$ (say, $\pi_a$ is the
$<^*_\chi$--first such function), and for $\alpha<\lambda$ let $w^a_\alpha
=(\pi_a)^{-1}[i_\alpha]$ (so $a\cap\zeta^*=\bigcup\limits_{\alpha<\lambda}
w^a_\alpha$).  

For $\zeta\leq\zeta^*$ let $\gR^{[\zeta]}$ consist of all triples
$(p,\delta,\bar{x})\in \bbP_\zeta\times\lambda\times A$ such that for some 
non-empty $w\in [\zeta\cap A]^{\textstyle {<}\lambda}$ we have
\[\bar{x}=\langle x_\vare:\vare\in w\rangle\quad\mbox{ and }\quad (\forall 
\vare\in w)(p\rest\vare\forces_{\bbP_\vare}\mbox{`` }(p(\vare),\delta,
x_\vare)\in\name{\gR}_\vare\mbox{ ''}).\]
Next, for $\zeta\leq\zeta^*$, $a\in W$ and $\delta<\lambda$ we put
\[\gY^{[\zeta]}_a(\delta)=\prod\{\gY^\vare_a(\delta):\vare\in w^a_\delta
\cap\zeta\},\]
thus defining $\gY^{[\zeta]}_a$ and $\bar{\gY}^{[\zeta]}=\langle \gY^{[
\zeta]}_a:a\in W\rangle$. If $\zeta=\zeta^*$ we will omit it (so then we
write $\gR$ and $\bar{\gY}$).  

\begin{claim}
\label{cl7}
For each $\zeta\leq \zeta^*$, $(\gR^{[\zeta]},\bar{\gY}^{[\zeta]})$ is a
$\lambda$--base for $\bbP_\zeta$ over $W$. 
\end{claim}

\begin{proof}[Proof of the Claim]
Immediately by the definition of $\gR^{[\zeta]},\bar{\gY}^{[\zeta]}$ we see
that clauses (a), (b) of \ref{base} hold (note: $\gY^{[\zeta]}_a(\delta)
\subseteq a$ by \ref{Con3.1}(1)). Now, to verify \ref{base}(c), suppose
$q\in\bbP_\zeta$, $a\in W$ and $\delta<\lambda$ is limit. For each
$\vare\in w^a_\delta\cap\zeta$ let $p'(\vare),\name{x}_\vare'$ be
$\bbP_\vare$--names such that 
\[q\rest\vare\forces_{\bbP_\vare}\mbox{`` }p'(\vare)\geq q(\vare)\ \&\
\name{x}_\vare'\in \gY^\vare_a(\delta)\ \&\ (p'(\vare),\delta,
\name{x}_\vare')\in\name{\gR}_\vare\mbox{ '',}\]
and for $\vare\in\Dom(q)\setminus w^a_\delta$ let
$p'(\vare)=q(\vare)$. This defines a condition $p'\in\bbP_\zeta$ stronger
than $q$ (and names $\name{x}_\vare'$). Since $\bbP_\zeta$ is
$\lambda$--complete we may find a condition $p\geq p'$ and $x_\vare\in
\gY^\vare_a(\delta)$ (for $\vare\in w^a_\delta\cap\zeta$) such that
$p\rest\vare\forces_{\bbP_\vare}$`` $\name{x}_\vare'=x_\vare$ '' (for
$\vare\in w^a_\delta\cap\zeta$). Then, by \ref{base}(a), we have $p\rest
\vare\forces_{\bbP_\vare}$`` $(p(\vare),\delta,x_\vare)\in\name{\gR}_\vare$
'' (for each $\vare\in w^a_\delta\cap\zeta$), and hence $(p,\delta,\langle
x_\vare:\vare\in w^a_\delta\cap\zeta\rangle)\in \gR^{[\zeta]}$.  
\end{proof}

Our aim now is to show that $\bbP_{\zeta^*}$ is weakly fuzzy proper with
witnesses $(\gR,\bar{\gY})$ and $c=(\langle \name{c}_\vare:\vare<\zeta^* 
\rangle, \langle\name{\gR}_\vare:\vare<\zeta^*\rangle,\gR,\bar{\gY},S,D,
\bar{F}',\bar{\bbQ})$. So suppose that a model $N\prec ({\mathcal H}
(\chi),\in,<^*_\chi)$ satisfies 
\[|N|=\lambda,\quad {}^{<\lambda}N\subseteq N,\quad \lambda,\bbP_{\zeta^*},
c\in N,\quad a\stackrel{\rm def}{=}N\cap A\in W,\]
and $p\in\bbP_{\zeta^*}\cap N$, and $h:\lambda\longrightarrow N$ is such
that $\bbP_{\zeta^*}\cap N\subseteq \Rang(h)$. To simplify the notation
later, let $\pi=\pi_a$, $w_\alpha=w^a_\alpha$ (for $\alpha<\lambda$).

Let us fix a list $\bcI=\langle\cI_\alpha:\alpha<\lambda\rangle$ of all open
dense subsets of $\bbP_{\zeta^*}$ from $N$. For $\zeta\in (\zeta^*+1)\cap
N$, let $\bcI^{[\zeta]}=\langle \cI^{[\zeta]}_\alpha:\alpha<\lambda\rangle$,
where $\cI^{[\zeta]}_\alpha=\{p\restriction\zeta:p\in\cI_\alpha\}$. (Note
that $\bcI^{[\zeta]}$ lists all open dense subsets of $\bbP_\zeta$ from
$N$.) Also for $\zeta\in\zeta^*\cap N$ let $\cJ_\zeta=\{p\in\bbP_{\zeta^*}:
p\rest\zeta\forces p(\zeta)\neq\name{\emptyset}_{\name{\bbQ}_\zeta}\}$ (so
$\cJ_\zeta$ is an open dense subset of $\bbP_{\zeta^*}$ from $N$) and let 
$E_3=\{\alpha\in E_2:(\forall\zeta\in w_\alpha)(\exists\beta<\alpha)(
\cJ_\zeta=\cI_\beta)\}$. Clearly, $E_3$ is a club of $\lambda$.

Now, using the diamond $\bar{F}'$  fixed earlier, we are going to define the
quasi $D$--diamond sequence $\bar{F}$ (and then a weak fuzzy candidate
$\bar{q}$ over it) that will be as required by $(\circledast)^-$ of
\ref{1.2}(3).  So, for each $\delta\in S$ we let  
\[\begin{array}{ll}
Z(\delta)=\big\{\zeta\in(\zeta^*+1)\setminus\{0\}:&\langle (h\circ F_\delta'
(\alpha)\restriction\zeta:\alpha<\delta\rangle\mbox{ is a
$\leq_{\bbP_\zeta}$--increasing}\\  
&\bcI^{[\zeta]}\mbox{--exact sequence of members of }N\cap\bbP_\zeta\ \big\} 
  \end{array}\]
and if $Z(\delta)\neq\emptyset$ then we put $\gamma(\delta)=\sup(Z(
\delta))$. Note that $Z(\delta)\in N$ and thus $\gamma(\delta)\in N$ (when
defined). Now, the pre--diamond $\bar{F}=\langle F_\delta:\delta\in
S\rangle$ is picked so that for a limit $\delta\in S$: 
\begin{enumerate}
\item[$(\odot)_1$] if $Z(\delta)\neq\emptyset$, then $h\circ F_\delta(
\alpha)=\big(h\circ F_\delta'(\alpha)\big)\restriction\gamma(\delta)$ for
all $\alpha<\delta$; 
\item[$(\odot)_2$] if $Z(\delta)=\emptyset$, then $h\circ F_\delta(\alpha
)=\emptyset_{\bbP_{\zeta^*}}$ for all $\alpha<\delta$.
\end{enumerate}
Then easily $\bar{F}$ is a quasi $D$--diamond for $(N,h,\bbP_{\zeta^*})$ and
for each limit $\delta\in S$, $\langle h\circ F_\delta(\alpha):\alpha<\delta
\rangle$ is a $\leq_{\bbP_{\zeta^*}}$--increasing sequence of conditions
from $\bbP_{\zeta^*}\cap N$.

Just for notational simplicity, we will identify a sequence $\bar{\sigma}=
\langle\sigma_0\rangle$ with its (only) term $\sigma_0$. Thus below, when we
talk about a standard $(w,1)^{\zeta^*}$--tree $\cT$, we think of $T$ as a
set of sequences $t=\langle (t)_\zeta:\zeta\in w\cap\rk(t)\rangle$ where 
$(t)_\zeta$'s do not have to be sequences.

Now we are going to define sequences $\bar{p}=\langle p_i:i<\lambda\rangle
\subseteq\bbP_{\zeta^*}\cap N$, $\langle \cT_\delta:\delta\in S\mbox{ is
limit }\rangle$, and $\langle q_{\delta,t}: \delta\in S\mbox{ is limit }\
\&\ t\in T_\delta\rangle\subseteq\bbP_{\zeta^*}\cap N$ such that for a limit
ordinal $\delta\in S$: 
\begin{enumerate}
\item[(i)] $\cT_\delta=(T_\delta,\rk_\delta)$ is a standard $(w_\delta,
1)^{\zeta^*}$--tree, and (under the identification mentioned earlier)
$\{t\in T_\delta:\rk_\delta(t)=\zeta\}\subseteq \gY^{[\zeta]}_a(\delta)$ for 
$\zeta\in w_\delta\cup\{\zeta^*\}$,  
\item[(ii)] $\langle q_{\delta,t}:t\in T_\delta\rangle$ is a standard tree
of conditions in $\bar{\bbQ}$, 
\item[(iii)] $p\leq p_i\leq p_j$ for $i<j<\lambda$,
\item[(iv)]  if $j<\lambda$, then $w_j\subseteq\Dom(p_j)$ and $(\forall
\vare\in w_j)(\forall j'>j)(p_j(\vare)=p_{j'}(\vare))$,
\item[(v)] if $t\in T_\delta$, $\rk_\delta(t)=\zeta$, then $q_{\delta,t}\in
\bbP_\zeta\cap N$ is such that 
\begin{enumerate}
\item[(a)] $\Big(\bigcup\limits_{\alpha<\delta}\Dom(h\circ F_\delta(\alpha))
\cup\bigcup\limits_{i<\delta}\Dom(p_i)\Big)\cap\zeta\subseteq\Dom(
q_{\delta,t})$, 
\item[(b)] $(\forall\alpha<\delta)\big( (h\circ F_\delta)(\alpha)\rest\zeta
\leq q_{\delta,t}\big)$, and 
\item[(c)] if $\vare\in\Dom(q_{\delta,t})\setminus w_\delta$, then  
\[\begin{array}{ll}
q_{\delta,t}\restriction\vare\forces&\mbox{`` if the set }\{p_i(
\vare):i<\delta\}\cup \{\big(h\circ F_\delta(\alpha)\big)(
\vare):\alpha<\delta\}\\
&\mbox{ \ has an upper bound in }\name{\bbQ}_\vare,\\
&\mbox{ \ then $q_{\delta,t}(\vare)$ is such an upper bound '',}
  \end{array}\]
\item[(d)] if $\vare\in\Dom(q_{\delta,t})\cap w_\delta$, then  
\[\begin{array}{ll}
q_{\delta,t}\restriction\vare\forces&\mbox{`` if the set }\{p_i(
\vare):i<\delta\}\cup \{\big(h\circ F_\delta(\alpha)\big)(
\vare):\alpha<\delta\}\\
&\mbox{ \ has an upper bound which obeys }(t)_\vare\mbox{ at }\delta,\\ 
&\mbox{ \ then $q_{\delta,t}(\vare)$ is such an upper bound,}\\
&\mbox{ \ else $q_{\delta,t}(\vare)$ is an upper bound of }\{\big(h\circ
F_\delta(\alpha)\big)(\vare):\alpha<\delta\}\\
&\mbox{ \ which obeys $(t)_\vare$ at $\delta$ '',}
  \end{array}\]
\item[(e)] $q_{\delta,t}\in \bigcap\limits_{\xi<\delta}\cI_\xi^{[\zeta]}$,  
\end{enumerate}
\item[(vi)] if $t\in T_\delta$, $\zeta=\rk_\delta(t)<\zeta^*$, $\zeta'\in
w_\delta\cup\{\zeta^*\}$ is the successor of $\zeta$ in $w_\delta\cup\{
\zeta^*\}$ and $t',t''\in T_\delta$ are such that $\rk_\delta(t)=\rk_\delta
(t'')=\zeta'$, $t\vtl t'$, $t\vtl t''$ and $t'\neq t''$, then 
\[q_{\delta,t}\forces_{\bbP_\zeta}\mbox{`` the conditions $q_{\delta,t'}(
\zeta)$ and $q_{\delta,t''}(\zeta)$ are incompatible '',}\]
\item[(vii)]  if $t\in T_\delta$ and $\vare\in\Dom(q_{\delta,t})\setminus
w_\delta$, then $\vare\in\Dom(p_\delta)$ and
\[\begin{array}{ll}
p_\delta\rest\vare\forces&\mbox{`` if }q_{\delta,t}\rest\vare\in
\Gamma_{\bbP_\vare}\mbox{ and }\{p_i(\vare):i<\delta\}\cup \{q_{\delta,
t}(\vare)\}\mbox{ has an upper bound in }\name{\bbQ}_\vare,\\ 
&\mbox{ \ then $p_\delta(\vare)$ is such an upper bound '',}
  \end{array}\]
\item[(viii)] if $t\in T_\delta$, $\rk_\delta(t)=\zeta<\zeta^*$, $x\in
\gY^\zeta_a(\delta)$ and 
\[\begin{array}{ll}
q_{\delta,t}\not\forces_{\bbP_\zeta}&\mbox{`` there is no condition stronger
than all }\\
&\ \ (h\circ F_\delta(\alpha))(\zeta)\mbox{ for }\alpha<\delta\mbox{ which
obeys $x$ at }\delta\mbox{ '',} 
  \end{array}\]
then there is $t'\in T_\delta$ such that $t\vtl t'$ and $(t')_\zeta=x$.
\end{enumerate}
Assume $\delta<\lambda$ and we have defined $p_i,\cT_i,q_{i,t}$ for relevant
$i<\delta$ and $t$. If $\delta$ is not a limit ordinal from $S$, then only
$p_\delta\in \bbP_{\zeta^*}\cap N$ needs to be defined, and clauses (iii),
(iv) can be easily taken care of. So suppose that $\delta\in S$ is limit. 

First we let $\cT_\delta'$ be a standard $(w_\delta,1)^{\zeta^*}$--tree
such that $\{t\in T_\delta':\rk_\delta(t)=\zeta\}=\gY^{[\zeta]}_a(\delta)$
(for $\zeta\in w_\delta\cup\{\zeta^*\}$). For $t\in T_\delta'$ we define a
condition $r_t\in\bbP_{\rk_\delta'(t)}$ so that 
\[\Dom(r_t)=\Big(\bigcup\limits_{\alpha<\delta}\Dom(h\circ F_\delta(\alpha))
\cup\bigcup\limits_{i<\delta}\Dom(p_i)\cup w_\delta\Big)\cap \rk_\delta'
(t),\] 
and for each $\zeta\in\Dom(r_t)$, $r_t(\zeta)$ is the $<^*_\chi$--first
$\bbP_\zeta$--name for a condition in $\name{\bbQ}_\zeta$ such that:

if $\zeta\in w_\delta$, then 
\[\begin{array}{ll}
r_t\rest\zeta\forces_{\bbP_\zeta}&\mbox{`` if the family }\{p_i(\zeta):i<
\delta\}\cup \{\big(h\circ F_\delta(\alpha)\big)(\zeta):\alpha<\delta\}\\
&\mbox{ \ has an upper bound which obeys }(t)_\zeta\mbox{ at }\delta,\\
&\mbox{ \ then $r_t(\zeta)$ is such an upper bound,}\\
&\mbox{ \ if the previous is impossible, but there is an upper bound of}\\
&\ \ \{\big(h\circ F_\delta(\alpha)\big)(\zeta):\alpha<\delta\}\mbox{ which
obeys $(t)_\zeta$ at $\delta$}\\
&\mbox{ \ then $r_t(\zeta)$ is such an upper bound,}\\
&\mbox{ \ if neither from the previous two possibilities holds,}\\
&\mbox{ \ then $r_t(\zeta)$ is an upper bound of }\{\big(h\circ F_\delta(
\alpha)\big)(\zeta):\alpha<\delta\}\mbox{ '',} 
  \end{array}\]
and if $\zeta\notin w_\delta$, then 
\[\begin{array}{ll}
r_t\rest\zeta\forces_{\bbP_\zeta}&\mbox{`` if the family }\{p_i(\zeta):i<
\delta\}\cup \{\big(h\circ F_\delta(\alpha)\big)(\zeta):\alpha<\delta\}\\
&\mbox{ \ has an upper bound, then $r_t(\zeta)$ is such an upper bound,}\\
&\mbox{ \ if this is not possible, then $r_t(\zeta)$ is just an upper bound
of }\\
&\ \ \{\big(h\circ F_\delta(\alpha)\big)(\zeta):\alpha<\delta\}\mbox{ ''.} 
  \end{array}\]
Plainly, $|T_\delta'|<\lambda$ and $\bar{r}=\langle r_t:t\in T_\delta'
\rangle$ is a standard tree of conditions, and it belongs to $N$ (remember:
${}^{<\lambda}N\subseteq N$). So using \ref{pA.7} in $N$ we may pick a
standard tree of conditions $\bar{r}^*=\langle r^*_t:t\in T_\delta'\rangle 
\in N$ such that $\bar{r}\leq\bar{r}^*$ and for each $t\in T_\delta'$ and
$\zeta\in w_\delta\cap\rk_\delta'(t)$ the condition $r^*_t\rest \zeta$
decides the truth value of the sentence 
\[\mbox{`` }r^*_t(\zeta)\mbox{ obeys }(t)_\zeta\mbox{ at }\delta\mbox{ (with
respect to $\name{\gR}_\zeta$) ''}\] 
(remember the choice of $r_t(\zeta)$ for $\zeta\in w_\delta$ and
\ref{base}(a)). Put   
\[T_\delta=\big\{t\in T_\delta':\mbox{ for each }\zeta\in w_\delta\cap
\rk_\delta'(t),\ r^*_t\rest \zeta\forces_{\bbP_\zeta}\mbox{``
}r^*_t(\zeta)\mbox{ obeys }(t)_\zeta\mbox{ at }\delta\mbox{ '' }\big\},\] 
and notice that $T_\delta\in N$ is a standard $(w_\delta,1)^{\zeta^*
}$--tree. Note also that for each $t\in T_\delta$ there is $t'\in T_\delta$
such that $t\trianglelefteq t'$ and $\rk_\delta(t')=\zeta^*$. Finally, using
\ref{obsA.4} and next \ref{pA.7} in $N$ we may choose a standard tree of
conditions $\langle q_{\delta,t}:t\in T_\delta\rangle\in N$ which satisfies
clauses (vi) and (v)(e) and such that $r^*_t\leq q_{\delta,t}$ (for $t\in
T_\delta$). It should be clear that then $\cT_\delta$ and  $\langle
q_{\delta,t}:t\in T_\delta\rangle$ satisfy all the relevant demands from our
list ((i)--(viii)). Now finding a condition $p_\delta\in \bbP_{\zeta^*}\cap
N$ which satisfies (iii)+(iv)+(vii) is straightforward. (Note that, by (vi),
if $t',t''\in T_\delta$, $\vare\in\Dom(q_{\delta,t'})\cap\Dom(q_{\delta,
t''})$, $\vare\notin w_\delta$ and the conditions $q_{\delta,t'}\rest\vare,
q_{\delta,t''}\rest\vare$ are compatible, then $q_{\delta,t'}\rest(\vare+1)
=q_{\delta,t''}\rest(\vare+1)$.)
\medskip

For a limit ordinal $\delta\in S$ and $\zeta\in N\cap(\zeta^*+1)$ we let 
\begin{itemize}
\item $\cX_\delta^{[\zeta]}=\{t\rest\zeta:t\in T_\delta\ \&\ \rk_\delta(t)=
\zeta^*\}$; 
\item if $s\in \cX^{[\zeta]}_\delta$, then $q^{[\zeta]}_{\delta,s}=q_{
\delta,t}\rest\zeta$ for some (equivalently: all) $t\in T_\delta$ such that
$s\trianglelefteq t$ and $\rk_\delta(t)=\zeta^*$; 
\item $\bar{q}^{[\zeta]}=\langle q_{\delta,s}^{[\zeta]}:\delta\in S\mbox{ is
limit }\ \&\ s\in \cX^{[\zeta]}_\delta\rangle$;
\item $h^{[\zeta]}:\lambda\longrightarrow N$ is such that $h^{[\zeta]}(
\gamma)= (h(\gamma))\restriction\zeta$ provided $h(\gamma)$ is a function,
and $h^{[\zeta]}(\gamma)=\emptyset_{\bbP_\zeta}$ otherwise.  
\end{itemize}

Plainly, $\emptyset\neq\cX_\delta^{[\zeta]}\subseteq\gY_a^{[\zeta]}
(\delta)$ (remember (i)) and $h^{[\zeta]}:\lambda\longrightarrow N$ is such
that $\bbP_\zeta\cap N\subseteq\Rang(h^{[\zeta]})$. Moreover, one easily
verifies the following claim. 

\begin{claim}
\label{cl1}
Let $\zeta\in N\cap (\zeta^*+1)$. Then $\bar{F}$ is a quasi $D$--diamond
sequence for $(N,h^{[\zeta]},\bbP_\zeta)$ and $\bar{q}^{[\zeta]}$ is a weak
fuzzy candidate over $\bar{F}$ for $(N,h^{[\zeta]},\bbP_\zeta,\gR^{[\zeta]},
\bar{\gY}^{[\zeta]},\bcI^{[\zeta]})$.
\end{claim}

We may write $\bar{q}, q_{\delta,t},\cX_\delta$  for $\bar{q}^{[\zeta^*]},
q_{\delta,t}^{[\zeta^*]},\cX_\delta^{[\zeta^*]}$, respectively. Also note
that, in the context of our definitions, the functions $h$ and
$h^{[\zeta^*]}$ behave the same, so we may identify them. Of course, we are
going to define an $(\gR,\bar{\gY})$--fuzzy generic condition $r\in
\bbP_{\zeta^*}$ for $\bar{q}$ over $\bar{F}$, but before that we have to
introduce more notation used later and prove some important facts.  

For $\zeta\in \zeta^*\cap N$ we define a function $h^{\langle\zeta\rangle}$
and $\bbP_\zeta$--names $\name{S}^{\langle\zeta\rangle},\name{\cX}^{\langle
\zeta\rangle}_\delta,\name{\cI}_\alpha^{\langle\zeta\rangle},\name{\bar{F}
}^{\langle\zeta\rangle},\name{\bcI}^{\langle\zeta\rangle}$ and
$\name{\bar{q}}^{\langle\zeta\rangle}$ so that: 
\begin{itemize}
\item $h^{\langle\zeta\rangle}:\lambda\longrightarrow N$ is such that if
$h(\gamma)$ is a function, $\zeta\in\Dom(h(\gamma))$ and
$(h(\gamma))(\zeta)$ is a $\bbP_\zeta$--name then $h^{\langle\zeta\rangle}
(\gamma)=((h(\gamma))(\zeta)$, otherwise $h^{\langle\zeta\rangle}(\gamma)=
\name{\emptyset}_{\name{\bbQ}_\zeta}$;
\item $\forces_{\bbP_\zeta}$``\/$\name{S}^{\langle\zeta\rangle}=\{\delta 
\in S:$ if $\delta$ is limit then $\delta>\pi(\zeta)\ \&\ (\exists s\in
\cX_\delta^{[\zeta]})(q_{\delta,s}^{[\zeta]}\in\Gamma_{\bbP_\zeta})\}$\/'';  
\item $\forces_{\bbP_\zeta}$`` if $\delta\in\name{S}^{\langle\zeta\rangle}$
is limit, then $\name{\cX}^{\langle\zeta\rangle}_\delta=\{x\in\gY_a^\zeta(
\delta): (\exists t\in\cX_\delta)(q_{\delta,t}\rest\zeta\in
\Gamma_{\bbP_\zeta}\ \&\ (t)_\zeta=x)\}$ '';  
\item $\forces_{\bbP_\zeta}\mbox{`` }\name{\bar{q}}^{\langle\zeta\rangle}=
\langle\name{q}_{\delta,x}^{\langle\zeta\rangle}:\delta\in\name{S}^{\langle
\zeta\rangle}\mbox{ is limit }\ \&\ x\in\name{\cX}_\delta^{\langle\zeta
\rangle}\rangle$, where:
\item $\forces_{\bbP_\zeta}$`` if $\delta\in\name{S}^{\langle\zeta
\rangle}$ is limit and $x\in\name{\cX}_\delta^{\langle\zeta\rangle}$, then
$\name{q}_{\delta,x}^{\langle\zeta\rangle}=q_{\delta,t}(\zeta)$ for some
(equivalently: all) $t\in\cX_\delta$ such that $q_{\delta,t}\rest\zeta\in
\Gamma_{\bbP_\zeta}$ and $(t)_\zeta=x$ '';   
\item $\forces_{\bbP_\zeta}$`` $\name{\bar{F}}^{\langle\zeta\rangle}=\langle
F_\delta:\delta\in \name{S}^{\langle\zeta\rangle}\rangle$;
\item $\forces_{\bbP_\zeta}$`` $\name{\bcI}^{\langle\zeta\rangle}=\langle
\name{\cI}^{\langle\zeta\rangle}_\alpha:\alpha<\lambda\rangle$ '', where: 
\item $\forces_{\bbP_\zeta}$`` $\name{\cI}^{\langle\zeta\rangle}_\alpha=\{ 
p(\zeta): p\in\cI_\alpha\ \&\ p\restriction\zeta\in\Gamma_{\bbP_\zeta}\}$
''. 
\end{itemize}
Naturally, we treat $h^{\langle\zeta\rangle}$ as a $\bbP_\zeta$--name
for a function from $\lambda$ to $N[\Gamma_{\bbP_\zeta}]$. Observe that 
$\forces_{\bbP_\zeta}$`` $N[\Gamma_{\bbP_\zeta}]\cap\name{\bbQ}_\zeta
\subseteq\Rang(h^{\langle\zeta\rangle})$ '', and 

$\forces_{\bbP_\zeta}$`` $\name{\bcI}^{\langle\zeta\rangle}$ lists all open
dense subsets of $\name{\bbQ}_\zeta$ from $N[\Gamma_{\bbP_\zeta}]$ ''. 

\begin{claim}
\label{cl2}
Assume that $\zeta\in N\cap\zeta^*$ and $r\in\bbP_\zeta$ is a $(\gR^{[
\zeta]},\bar{\gY}^{[\zeta]})$--fuzzy generic condition for $\bar{q}^{[
\zeta]}$ over $(N,\bcI^{[\zeta]},h^{[\zeta]},\bbP_\zeta,\bar{F})$. Then 
\begin{enumerate}
\item $r\forces_{\bbP_\zeta}$`` $\name{S}^{\langle\zeta\rangle}\in D^+$ '', 
\item $r\forces_{\bbP_\zeta}$`` $\name{\bar{F}}^{\langle\zeta\rangle}$ is a 
quasi $D$--diamond for $(N[\Gamma_{\bbP_\zeta}],h^{\langle\zeta\rangle},
\name{\bbQ}_\zeta)$ '', and  
\item $r\forces_{\bbP_\zeta}$`` $\name{\bar{q}}^{\langle\zeta\rangle}$ is a
fuzzy candidate for $(N[\Gamma_{\bbP_\zeta}],h^{\langle\zeta\rangle},
\name{\bbQ}_\zeta,\name{\gR}_\zeta,\bar{\gY}^\zeta)$ over
$\name{\bar{F}}^{\langle\zeta\rangle}$ ''.    
\end{enumerate}
\end{claim}

\begin{proof}
(1)\quad Will follow from (2).
\medskip

\noindent (2)\quad Assume that this fails. Then we can find a condition
$r^*\in\bbP_\zeta$, a $\bbP_\zeta$--name $\name{\bar{q}}'=\langle
\name{q}_\alpha':\alpha<\lambda\rangle$ for an increasing sequence of
conditions from $\name{\bbQ}_\zeta\cap N[\Gamma_{\bbP_\zeta}]$, and
$\bbP_\zeta$--names $\name{A}_\xi,\name{B}_\xi$ for 
members of $D\cap\bV$ such that $r\leq_{\bbP_\zeta} r^*$ and 
\[\begin{array}{ll}
r^*\forces_{\bbP_\zeta}&\mbox{`` }(\forall\delta\in
\mathop{\triangle}\limits_{\xi<\lambda}\name{A}_\xi)(\langle
\name{q}_\alpha':\alpha<\delta\rangle\mbox{ is $\name{\bcI}^{\langle\zeta
\rangle}$--exact })\quad\mbox{ and}\\
&\quad (\forall\delta\in\name{S}^{\langle\zeta\rangle}\cap
\mathop{\triangle}\limits_{\xi<\lambda}\name{B}_\xi)(\langle h^{\langle
\zeta\rangle}\circ F_\delta(\alpha):\alpha<\delta\rangle\neq\name{\bar{q}}'
\restriction \delta)\mbox{ ''.}
  \end{array}\]  
Consider a play $\langle r_j,C_j:j<\lambda\rangle$ of the game $\Gfl(r,N,
\bcI^{[\zeta]},h^{[\zeta]},\bbP_\zeta,\bar{F},\bar{q}^{[\zeta]})$ in which
Generic uses her winning strategy and Antigeneric plays as follows.

Together with choosing $r_j$ (for $j\in\lambda\setminus S$), the Antigeneric 
player chooses also side conditions $p_j\in N\cap\bbP_\zeta$, sets $A_\xi, 
B_\xi\in D$ and $\bbP_\zeta$--names $\name{q}^*_\xi\in N$ for elements of
$\name{\bbQ}_\zeta$ (for $\xi\leq j$) such that 
\begin{itemize}
\item $r_j\geq r^*$ (so $r_0\geq r^*$; remember Antigeneric plays at 0),
$r_j\geq r_i$ (for $i<j$), and $r_j\geq p_j$ and 
\item $r_j\forces_{\bbP_\zeta}$`` $(\forall\xi\leq j)(\name{A}_\xi=A_\xi\ \&\
\name{B}_\xi=B_\xi\ \&\ \name{q}'_\xi=\name{q}^*_\xi)$ '', and 
\item if $j'<j$ are from $\lambda\setminus S$, then $p_j\geq p_{j'}$, and  
$p_j\in\bigcap\limits_{\xi<j}\cI^{[\zeta]}_\xi$, and 
\item $p_j\forces_{\bbP_\zeta}$`` $(\forall\xi_0<\xi_1\leq j)(q^*_{\xi_0}
  \leq q^*_{\xi_1})$ '', and 
\item if $\delta<j$, $\delta\in\bigcap\limits_{\xi<\delta}A_\xi$, then
${p_j}^{\frown}\langle\name{q}_j^*\rangle\in\cI_\xi^{[\zeta+1]}$ for all
$\xi<\delta$. 
\end{itemize}
[The $C_j$'s are not that important for our argument, so we do not specify
any requirements on them. Regarding the choice of the $p_j$'s, remember
\ref{1.2A}(2); for the last two demands remember that $\name{q}_j'$'s are
(forced to be) increasing.] After the play, the Antigeneric player completes 
$\langle p_j:j\in\lambda\setminus S\rangle$ to a
$\leq_{\bbP_\zeta}$--increasing sequence $\langle p_j:j<\lambda\rangle
\subseteq N\cap\bbP_\zeta$ letting $p_j=p_{\min(\lambda\setminus S\setminus
(j+1))}$ for $j\in S$.  

Note that if $\delta\in E_0$ is a limit of elements of
$\mathop{\triangle}\limits_{\xi<\lambda} A_\xi$, then the sequence $\langle
{p_j}^\frown\langle\name{q}_j^*\rangle:j<\delta\rangle$ is $\bcI^{[\zeta+
1]}$--exact and increasing (and $\langle p_j:j<\delta\rangle$ is
$\bcI^{[\zeta]}$--exact).  So, as $D$ is normal and $A_\xi,B_\xi,C_j\in D$
and $\bar{F}$ is a quasi $D$--diamond for $(N,h^{[\zeta+1]},\bbP_{\zeta+1})$
(by \ref{cl1}), we may find an ordinal $\delta\in S\cap E_0\cap
\mathop{\triangle}\limits_{\xi<\lambda}A_\xi\cap
\mathop{\triangle}\limits_{\xi<\lambda} B_\xi\cap\mathop{\triangle}\limits_{
j<\lambda}C_j\setminus(\pi(\zeta)+1)$ which is a limit of members of
$\mathop{\triangle}\limits_{\xi<\lambda}A_\xi$ and such that $\langle
h^{[\zeta+1]}\circ F_\delta(j):j<\delta\rangle=\langle {p_j}^\frown
\langle\name{q}_j^* \rangle:j<\delta\rangle$. By the choice of $\bar{F}$ we
know that $h(F_\delta(j))$ is a condition in $\bbP_{\zeta^*}$ (so a
function) and hence $h^{\langle\zeta\rangle}(F_\delta(j))=\name{q}^*_j$ for
all $j<\delta$. Also 
\[(\forall i<\delta)(\exists j\in\delta\setminus S)(h^{[\zeta]}\circ
F_\delta(i)\leq_{\bbP_\zeta} h^{[\zeta]}\circ F_\delta(j)=p_j
\leq_{\bbP_\zeta} r_j).\]
Since the play is won by Generic, for some $s\in\cX_\delta^{[\zeta]}$ we
have $q_{\delta,s}^{[\zeta]}\leq r_\delta$. But then   
\[r_\delta\forces\mbox{`` }\delta\in S^{\langle\zeta\rangle}\cap
\mathop{\triangle}_{\xi<\lambda}\name{B}_\xi\quad \&\quad \langle h^{\langle
\zeta\rangle}\circ F_\delta(\alpha):\alpha<\delta\rangle=\bar{q}'\rest
\delta\mbox{ '',}\] 
a contradiction.
\medskip

\noindent (3)\quad Straightforward (remember the choice of $q_{\delta,t}$'s,
specifically clauses (v)(b,d,e) and (viii)). 
\end{proof}

Now we are going to define an $(\gR,\bar{\gY})$--fuzzy generic condition 
$r\in\bbP$ for $\bar{q}$ over $(N,\bcI,h,\bbP_{\zeta^*},\bar{F})$ in the
most natural way. Its domain is $\Dom(r)=\zeta^*\cap N$ and for each $\zeta\in
\zeta^*\cap N$ 
\[\begin{array}{ll}
r\restriction\zeta\forces&\mbox{`` }r(\zeta)\geq p_{\pi(\zeta)+1}(\zeta)
\mbox{ is an $(\name{\gR}_\zeta,\bar{\gY}^\zeta)$--fuzzy generic condition for
}\bar{q}^{\langle\zeta\rangle}\\
&\ \ \mbox{ over }(N[\Gamma_{\bbP_\zeta}],\name{\bcI}^{\langle\zeta\rangle},
h^{\langle\zeta\rangle},\name{\bbQ}_\zeta,\name{\bar{F}}^{\langle\zeta
\rangle})\mbox{ ''.}
  \end{array}\]
[So $r\geq p_i$ for all $i<\lambda$.]

\begin{claim}
\label{cl6}
For every $\zeta\in (\zeta^*+1)\cap N$, the Generic player has a winning
strategy in the game $\Gfl(r\restriction\zeta,N,\bcI^{[\zeta]},h^{[\zeta]},
\bbP_\zeta,\bar{F},\bar{q}^{[\zeta]})$.
\end{claim}

\begin{proof}[Proof of the Claim]
We will prove the claim by induction on $\zeta\in (\zeta^*+1)\cap N$. After
we are done with stage $\zeta\in(\zeta^*+1)\cap N$, we know that $r\rest
\zeta$ is $(\gR^{[\zeta]},\bar{\gY}^{[\zeta]})$--fuzzy generic for
$\bar{q}^{[\zeta]}$ over $(N,\bcI^{[\zeta]},h^{[\zeta]},\bbP_\zeta,
\bar{F})$. For $\zeta\in\zeta^*\cap N$ this implies that $r(\zeta)$ is
well-defined (remember \ref{cl2}). Of course for $\zeta=\zeta^*$ we finish
the proof of the theorem.    
\medskip

Suppose that $\zeta\in (\zeta^*+1)\cap N$ and we know that $r\restriction
\zeta'$ is $(\gR^{[\zeta']},\bar{\gY}^{[\zeta']})$--fuzzy generic for
$\bar{q}^{[\zeta']}$ over $(N,\bcI^{[\zeta']},h^{[\zeta']},\bbP_{\zeta'},
\bar{F})$ for all $\zeta'\in N\cap\zeta$. We are going to define a winning
strategy $\st$ for Generic in the game $\Gfl(r\rest\zeta,N,\bcI^{[\zeta]},
h^{[\zeta]},\bbP_\zeta,\bar{F},\bar{q}^{[\zeta]})$. First, for $\vare\in
\zeta\cap N$ fix a $\bbP_\vare$--name $\nst_\vare$ such  
\[\begin{array}{ll}
r\rest\vare\forces&\mbox{``  $\nst_\vare$ is a winning strategy of the 
Generic player}\\
&\ \mbox{ in the game }\Gfl(r(\vare),N[\Gamma_{\bbP_\vare}],
\name{\bcI}^{\langle\vare\rangle},h^{\langle\vare\rangle}, 
\name{\bbQ}_\vare,\bar{F}^{\langle\vare\rangle},\bar{q}^{\langle\vare
\rangle})\mbox{ ''.}
  \end{array}\]
We will think of $\nst_\vare$ as a name for a function from
${<}\lambda$--sequences of members of $\name{\bbQ}_\vare\times D$ (thought
of as pairs of sequences of the same length $<\lambda$) to 
$\name{\bbQ}_\vare\times D$ such that if $(\name{\bar{\sigma}},
\name{\bar{C}})\in\Dom(\nst_\vare)$ and $\name{\bar{\sigma}}$ has an upper
bound, then the first coordinate of $\st_\vare(\name{\bar{\sigma}},
\name{\bar{C}})$ is such an upper bound (and, of course, any play according
to $\nst_\vare$ is won by Generic). (In a play of $\Gfl(r(\vare),N[
\Gamma_{\bbP_\vare}],\name{\bcI}^{\langle\vare\rangle},h^{\langle\vare
\rangle},\name{\bbQ}_\vare,\bar{F}^{\langle\vare\rangle},\bar{q}^{\langle
\vare\rangle})$ only the values of $\nst_\vare$ at ``legal partial plays
according to $\nst_\vare$'' matter, but it is notationally convenient to
have $\nst_\vare$ giving values for all sequences of elements of
$\name{\bbQ}_\vare\times D$, even if first coordinates are not increasing,
as well as for sequences after which Antigeneric should play.) 

We will describe the strategy $\st$ by giving the answers of Generic on
intervals $S\cap [i_\alpha,i_{\alpha+1})$ (for $\alpha<\lambda$), where,
remember, $\langle i_\alpha:\alpha<\lambda\rangle$ is the increasing
enumeration of $E_1$. Aside the Generic player will construct sequences
$\langle\name{r}_{j'}'(\vare):j'<\lambda,\ \vare\in\zeta\cap N\rangle$ and
$\langle\name{C}_{j'}^\xi(\vare):j',\xi<\lambda,\ \vare\in\zeta\cap
N\rangle$ so that   
\begin{enumerate}
\item[$(*)_1$] $\name{r}_{j'}'(\vare)$ is a $\bbP_\vare$--name for a
member of $\name{\bbQ}_\vare$, $\name{C}_{j'}^\xi(\vare)$ is a 
$\bbP_\vare$--name for a member of $D\cap\bV$, and
\item[$(*)_2$] if $\alpha<\lambda$, $\delta=\min\big(S\cap [i_\alpha,
i_{\alpha+1})\big)$, and $\vare\in w_{\alpha+1}\cap\zeta$, then after the
$\delta$-th move (which is a move of the Generic player) the terms
$\name{r}_{j'}'(\vare),\langle\name{C}_{j'}^\xi(\vare):\xi<\lambda\rangle$
are defined for all $j'<i_{\alpha+1}$, and 
\item[$(*)_3$] if $\alpha<\lambda$, $\vare\in w_{\alpha+1}\cap\zeta$ and
$p^*\in \bbP_\vare$ is stronger than all $r_j\rest\vare$ for $j\in (i_\alpha
+1)\setminus S$ ($r_j$ are the conditions played in the game), then 
\[p^*\forces_{\bbP_\vare}\mbox{`` }(\forall j\in (i_\alpha+1)\setminus S)(
r_j(\vare)\leq \name{r}_{i_\alpha}'(\vare))\mbox{ '',}\]
\item[$(*)_4$] if $\alpha<\lambda$,  $\vare\in w_{\alpha+1}\cap\zeta$ and
$r_{i_\alpha+1}$ is the condition played by Generic at stage $i_\alpha+1\in
S$, then 
\[r_{i_\alpha+1}\rest\vare\forces_{\bbP_\vare}\mbox{`` }(\forall j'<i_{
\alpha+1})(\name{r}_{j'}'(\vare)\leq r_{i_\alpha+1}(\vare))\mbox{ '',}\]
\item[$(*)_5$] for each $\vare\in N\cap\zeta$, 
\[\begin{array}{ll}
r\rest\vare\forces_{\bbP_\vare}&\mbox{`` }\langle \name{r}'_j(\vare),
\mathop{\triangle}\limits_{\xi<\lambda}\name{C}_j^\xi(\vare):j<\lambda\rangle 
\mbox{ is a play of the game}\\
&\ \Gfl(r(\vare),N[\Gamma_{\bbP_\vare}],\bcI^{\langle\vare
\rangle},h^{\langle\vare\rangle},\name{\bbQ}_\vare,\bar{F}^{\langle\vare
\rangle},\name{\bar{q}}^{\langle\vare\rangle})\mbox{ in which Generic uses
$\nst_\vare$ ''.}
  \end{array}\]
\end{enumerate}
So suppose that $\alpha<\lambda$, $\delta=\min\big(S\cap [i_\alpha,i_{\alpha
+1})\big)$ and $\langle r_j,C_j:j<\delta\rangle$ is the result of the play
so far. Now Generic looks at ordinals $\vare\in w_{\alpha+1}\cap\zeta$. She 
lets the $\bbP_\vare$--names $\name{r}_{j'}'(\vare),\name{C}_{j'}^\xi(
\vare)$ be such that $\langle \name{r}_{j'}'(\vare),
\mathop{\triangle}\limits_{\xi<\lambda}\name{C}_{j'}^\xi(\vare):j'<i_{\alpha
+1}\rangle$ is forced by $r\rest\vare$ to be a play of $\Gfl(r(\vare),N[
\Gamma_{\bbP_\vare}],\bcI^{\langle\vare\rangle},h^{\langle\vare\rangle}, 
\name{\bbQ}_\vare,\bar{F}^{\langle\vare\rangle},\name{\bar{q}}^{\langle\vare
\rangle})$ in which the moves are determined as follows. If $\vare\in
w_\alpha$, then we have already the play below $i_\alpha$ and the names
$\name{r}_{i_\alpha}'(\vare),\name{C}_{i_\alpha}^\xi(\vare)$ are such that  
\begin{itemize}
\item if $i_\alpha=\delta$ (i.e., $i_\alpha\in S$ and $r_i,C_i$ have been
chosen for $i<i_\alpha$ only), then   
\[\begin{array}{ll}
r\rest\vare\forces_{\bbP_\vare}&\mbox{`` }(\name{r}_{i_\alpha}'(\vare), 
\mathop{\triangle}\limits_{\xi<\lambda}\name{C}_{i_\alpha}^\xi(\vare))
\mbox{ is the value of $\nst_\vare$}\\
&\mbox{ applied to }\langle r'_j(\vare),\mathop{\triangle}\limits_{\xi<
\lambda}\name{C}_j^\xi(\vare):j<i_\alpha\rangle\mbox{, ''}
  \end{array}\]
\item if $i_\alpha<\delta$ (i.e., $i_\alpha\notin S$ so $r_i,C_i$ are already
chosen for $i\leq i_\alpha$), then   
\[\begin{array}{ll}
r\rest\vare\forces_{\bbP_\vare}&\mbox{`` if }(\forall j<i_\alpha)
(\name{r}_j'(\vare)\leq r_{i_\alpha}(\vare))\mbox{ then } 
\name{r}_{i_\alpha}'(\vare)=r_{i_\alpha}(\vare)\\
&\mbox{ otherwise $\name{r}_{i_\alpha}'(\vare)$ is the first coordinate of
$\nst_\vare$ applied to}\\
&\mbox{ the play so far, and }\name{C}_{i_\alpha}^\xi(\vare)=
\bigcap\limits_{j\leq i_\alpha}C_j\mbox{ for all $\xi<\lambda$ ''.}  
  \end{array}\]
\end{itemize}
Then for $j\in (i_\alpha,i_{\alpha+1})$ (and $\vare\in w_\alpha\cap\zeta$)
the names $\name{r}_j'(\vare),\name{C}_j^\xi(\vare)$ are determined by
applying successively $\nst_\vare$, that is 
\[\begin{array}{ll}
r\rest\vare\forces_{\bbP_\vare}&\mbox{`` }(\name{r}_j'(\vare),
\mathop{\triangle}\limits_{\xi<\lambda}\name{C}_j^\xi(\vare))\mbox{ is the
value of $\nst_\vare$}\\ 
&\mbox{ applied to }\langle r'_{j'}(\vare),\mathop{\triangle}\limits_{\xi<
\lambda}\name{C}_{j'}^\xi(\vare):j'<j\rangle\mbox{. ''}  
  \end{array}\]
If $\vare\in (w_{\alpha+1}\setminus w_\alpha)\cap\zeta$, then the Generic
player defines the names $\name{r}_j'(\vare),\name{C}_j^\xi(\vare)$ somewhat
like above, but starting with subscript $j=0$. Thus 
\begin{itemize}
\item if $i_\alpha=\delta$, then   
\[\begin{array}{ll}
r\rest\vare\forces_{\bbP_\vare}&\mbox{`` }\name{r}_0'(\vare)\mbox{ is the
first coordinate of the value of }\nst_\vare\\ 
&\mbox{ at }\langle r_{j'}(\vare),\bigcap\limits_{i<i_\alpha}C_i:j'<i_\alpha
\rangle\\
&\mbox{ and }\name{C}_0^\xi(\vare)=\bigcap\limits_{i<i_\alpha} C_i\mbox{ for
all $\xi<\lambda$ '',} 
  \end{array}\]
\item if $i_\alpha<\delta$, then   
\[\begin{array}{ll}
r\rest\vare\forces_{\bbP_\vare}&\mbox{`` if }(\forall j<i_\alpha)(r_j(\vare)
\leq r_{i_\alpha}(\vare))\mbox{ then }\name{r}_0'(\vare)=r_{i_\alpha}(
\vare)\\
&\mbox{ otherwise $\name{r}_0'(\vare)$ is the first coordinate of the value
of $\nst_\vare$}\\ 
&\mbox{ at }\langle r_j(\vare),\bigcap\limits_{i<i_\alpha}C_i:j<i_\alpha
\rangle\\ 
&\mbox{ and }\name{C}_0^\xi(\vare)=\bigcap\limits_{j\leq i_\alpha}C_j\mbox{
for all $\xi<\lambda$ ''.}  
  \end{array}\]
\end{itemize}
Lastly, for $0<j<i_{\alpha+1}$ (and $\vare\in (w_{\alpha+1}\setminus
w_\alpha)\cap\zeta$) the names $\name{r}_j'(\vare),\name{C}_j^\xi(\vare)$
are determined by applying successively $\nst_\vare$ (like earlier).  

This defines the names $\name{r}_j'(\vare),\name{C}_j^\xi(\vare)$ for
$j<i_{\alpha+1}, \xi<\lambda$ and $\vare\in w_{\alpha+1}\cap\zeta$. It is
straightforward to check that the requirements $(*)_1$--$(*)_3$ and $(*)_5$
restricted to $\vare\in w_{\alpha+1}\cap\zeta$ (and with ``$j<\lambda$''
replaced by ``$j<i_{\alpha+1}$'') are satisfied. 

Next, using the fact that $\bbP_\zeta$ is $\lambda$--closed and $(*)_3$
of the choice above, Generic picks a condition $r^*\in\bbP_\zeta$ such that  
\begin{enumerate}
\item[$(*)_6$] $r^*\geq r_j$ for every $j<\delta$,
\item[$(*)_7$] $r^*\rest\vare\forces\mbox{`` }\name{r}_{j'}'(\vare)\leq
r^*(\vare)\mbox{ ''}$ for every $j'<i_{\alpha+1}$ and $\vare\in w_{\alpha+1}
\cap\zeta$, 
\item[$(*)_8$] $r^*\in\bigcap\limits_{\xi<i_{\alpha+1}}\cI_\xi$, and 
\item[$(*)_9$] for every $j',\xi<i_{\alpha+1}$ and $\vare\in w_{\alpha+1}
\cap\zeta$, the condition $r^*\rest\vare$ decides the value of
$\name{C}_{j'}^\xi(\vare)$, say $r^*\rest\vare\forces\mbox{`` }
\name{C}_{j'}^\xi(\vare)=C_{j'}^\xi(\vare)\mbox{ ''}$, where $C_{j'}^\xi(
\vare)\in D\cap\bV$. 
\end{enumerate}
If $i_\alpha\in S$ (so $\delta=i_\alpha$ is a limit ordinal), then Generic
picks a condition $r^+\in\bbP_\zeta$ stronger than $r^*$ and such that for
every $t\in\cX_\delta^{[\zeta]}$ and $\vare\in (w_\delta\cap\zeta)\cup
\{\zeta\}$ we have:  
\begin{enumerate}
\item[$(*)_{10}$] either the conditions $r^+\rest\vare$ and
$q_{\delta,t}^{[\zeta]}\rest\vare$ are incompatible, or  $q_{\delta,t}^{
[\zeta]}\rest\vare\leq_{\bbP_\vare}r^+\rest\vare$, 
\item[$(*)_{11}$] if $\vare\in w_\delta\cap\zeta$ and $q_{\delta,t}^{[
\zeta]}\rest\vare\leq_{\bbP_\vare}r^+\rest\vare$, and $s\in\cX_\delta^{[
\zeta]}$ is such that $t\rest\vare=s\rest\vare$, then either $r^+\rest\vare
\forces$`` $q_{\delta,s}^{[\zeta]}(\vare),r^+(\vare)$ are incompatible '' or
$r^+\rest\vare\forces$``\/$q_{\delta,s}^{[\zeta]}(\vare)\leq r^+(\vare)$\/''.
\end{enumerate}
If $\delta>i_\alpha$ (i.e., $i_\alpha\notin S$) then Generic lets $r^+=r^*$. 
Finally, for every $j\in [i_\alpha,i_{\alpha+1})\cap S$ she plays 
\[r_j=r^+\quad\mbox{ and }\quad C_j=E_3\cap\bigcap\{C_{j'}^\xi(\vare):j',
\xi<i_{\alpha+1},\ \vare\in w_{\alpha+1}\cap\zeta\}.\]   
Plainly, $r_{i_\alpha+1}=r^+$ satisfies clause $(*)_4$.
\medskip

Why does the strategy described above work?\\
Suppose that $\langle r_j,C_j:j<\lambda\rangle$ is a play of the game
$\Game(r\restriction\zeta,N,\bcI^{[\zeta]}h^{[\zeta]},\bbP_{\zeta},\bar{F},
\bar{q}^{[\zeta]})$ in which the Generic player used this strategy and let
$\langle \name{r}_{j'}'(\vare):j'<\lambda,\ \vare\in\zeta\cap N\rangle$ and
$\langle\name{C}_{j'}^\xi(\vare):j',\xi<\lambda,\ \vare \in\zeta\cap
N\rangle$ be the sequences she constructed aside.
\medskip

First let us argue that condition \ref{pre1.2}(5)$(\beta)$ holds. We will
show slightly more than actually needed to help later with clause
$(\alpha)$. Below remember that ordinals $\gamma(\delta)$ were defined when
we picked our quasi $D$--diamond $\bar{F}$, and if $\vare<\gamma(\delta)$
then the sequence $\langle h^{[\vare+1]}\circ F_\delta(\alpha):\alpha<\delta 
\rangle$ is $\bcI^{[\vare+1]}$--exact. Now, suppose that a limit ordinal
$\delta\in S\cap\bigcap\limits_{j<\delta}C_j$ (so in particular $\delta\in
E_3$) is such that  
\begin{enumerate}
\item[$(\boxplus)_\delta$]\quad  $w_\delta\cap\zeta\subseteq\gamma(\delta)$
and $(\forall\alpha<\delta)(\exists j<\delta)(h^{[\zeta]}\circ F_\delta(
\alpha)\leq r_j)$.
\end{enumerate}
(So then $(\forall \alpha<\delta)(h^{[\zeta]}\circ F_\delta(\alpha)\leq
r_\delta)$. Note also that by the choice of $E_3$, if $\langle
h^{[\vare]}\circ F_\delta(\alpha):\alpha<\delta\rangle$ is
$\bcI^{[\vare]}$--exact, then $w_\delta\cap\zeta\subseteq\gamma(\delta)$.)

We are going to choose $t\in\cX_\delta^{[\zeta]}$ and show that $q_{\delta,
t}^{[\zeta]}\leq r_\delta$. We do this by induction on $\vare\in(\zeta+1)
\cap N$, defining $t\rest\vare\in T_\delta$ and showing that $q_{\delta,t}
\rest\vare=q_{\delta,t\rest\vare}\rest\vare\leq r_\delta\rest\vare$ (and for
$\vare=\zeta$ we get the desired conclusion). Limit stages and the initial
stage $\vare=0$ are trivial, so assume that we have defined $t\rest\vare$
and have shown $q_{\delta,t\rest \vare}\rest\vare= q_{\delta,t}\rest\vare
\leq r_\delta\rest\vare$ (where $\vare\in\zeta\cap N$), and let us consider
the restrictions to $\vare+1$. 

If $\vare\notin w_\delta$ then $t\rest (\vare+1)=t\rest\vare$ (so it has
been already defined). Suppose also that $\vare\in\Dom(q_{\delta,t\rest
\vare})$ (otherwise there is nothing to do). Look at the clause (v)(c) of the
choice of $q_{\delta,t\rest\vare}$ at the beginning: $r_\delta\geq r\geq
p_\delta$ (and $(\boxplus)_\delta$) implies that 
\[r_\delta\restriction\vare\forces\mbox{`` }q_{\delta,t\rest\vare}(\vare)
\mbox{ is an upper bound to }\{p_i(\vare):i<\delta\}\mbox{ ''.}\]
But then also by the clause (vii) there, $r_\delta\rest\vare\forces
\mbox{`` }q_{\delta,t\rest\vare}(\vare)\leq p_\delta(\vare)\leq
r_\delta(\vare)\mbox{ ''}$, so we are done. 

Suppose now that $\vare\in w_\delta\cap\zeta$ (and thus $\vare<\gamma(
\delta)$). Since $\delta\in E_3$, we know that arriving to stage $\delta$ of
the game, Generic has already defined $\name{r}_j'(\vare),\name{C}_j^\xi(
\vare)$ for $j<\delta$ and $\xi<\lambda$ (remember $(*)_2$). Moreover, the
condition $r_\delta\rest\vare$ forces that (remember: $\langle
h^{[\vare+1]}\circ F_\delta(\alpha):\alpha<\delta\rangle$ is
$\bcI^{[\vare+1]}$--exact):    
\begin{itemize}
\item the sequence $\langle h^{\langle\vare\rangle}\circ F_\delta(\alpha):
\alpha<\delta\rangle$ is $\leq_{\name{\bbQ}_\vare}$--increasing
$\name{\bcI}^{\langle\vare\rangle}$--exact, and 
\item $\langle\name{r}_j'(\vare),\mathop{\triangle}\limits_{\xi<\lambda}
\name{C}_j^\xi(\vare):j<\delta\rangle$ is a play according to $\nst_\vare$ 
(by $(*)_5$), and
\item  $\delta\in\name{S}^{\langle\vare\rangle}\cap\bigcap\limits_{j,\xi
<\delta}\name{C}_j^\xi(\vare)$ (remember $(*)_9$ and the choice of
$C_{i_\alpha+1}$ for $\alpha<\delta$), and hence also $\delta\in
\name{S}^{\langle\vare\rangle}\cap \bigcap\limits_{j<\delta}
\mathop{\triangle}\limits_{\xi<\lambda}\name{C}_j^\xi(\vare)$, 
\item $(\forall j<\delta)(\exists j'<\delta)(r_j(\vare)\leq\name{r}_{j'}'(
\vare))$ and $(\forall j<\delta)(\exists j'<\delta)(\name{r}_j'(\vare)\leq
r_{j'}(\vare))$ (by $(*)_3+(*)_4$), so also 
\[(\forall\alpha<\delta)(\exists j<\delta)(h^{\langle\zeta\rangle}\circ
F_\delta(\alpha)\leq_{\name{\bbQ}_\vare} \name{r}_j'(\vare)).\]
\end{itemize}
Since $\nst_\vare$ is a name for a winning strategy, we may conclude that  
(by $(*)_7$)
\[r_\delta\rest\vare\forces_{\bbP_\vare}\mbox{`` }(\exists x\in\name{\cX}^{
\langle\vare\rangle}_\delta)(q^{\langle\vare\rangle}_{\delta,x}\leq
\name{r}'_\delta(\vare)\leq r_\delta(\vare))\mbox{ ''.}\]
Now look at $(*)_{11}$ remembering clause (vi) of the choice of $\bar{q}$:
by them there is a unique $x\in\gY^\vare_a(\delta)$ such that letting
$(t)_\vare=x$ we get $t\rest (\vare+1)\in T_\delta$ satisfying 
$q_{\delta, t\rest (\vare+1)}\rest (\vare+1)\leq r_\delta\rest (\vare+1)$. 

This completes the inductive proof of \ref{pre1.2}(5)$(\beta)$.
\medskip

Why does \ref{pre1.2}(5)$(\alpha)$ hold? To show this condition, it is
enough to prove that $(\boxplus)_\delta$ holds for unboundedly many
$\delta\in S\cap\mathop{\triangle}\limits_{\xi<\lambda}C_j$ (remember
clause (v)(e) of the choice of $q_{\delta,t}$'s and what we have already
shown).  We do this considering various characters of $\zeta$.   
\medskip

\noindent {\bf $\zeta$ is a limit ordinal of cofinality
$\cf(\zeta)<\lambda$}.\\
Pick a closed set $u\subseteq\zeta$ such that $u\in N$, $0\in u$, $\otp(u)=
\cf(\zeta)$ and $\sup(u)=\zeta$. For $\alpha<\lambda$ let $\vare_\alpha\in
u$ be such that  $\alpha=\otp\big(u\cap \vare_\alpha\big)\mod\cf(\zeta)$.
Now, by induction on $\alpha<\lambda$ we choose conditions $s_\alpha\in N
\cap\bbP_\zeta$ such that 
\begin{enumerate}
\item[(a)$_\alpha$] $(\exists j<\lambda)(s_\alpha\leq r_j)$,
\item[(b)$_\alpha$] $s_\alpha\in\bbP_{\vare_\alpha}\cap N$, 
\item[(c)$_\alpha$] if $\beta<\alpha<\lambda$, then $s_\beta\rest (
\vare_\alpha\cap\vare_\beta)\leq s_\alpha\rest(\vare_\alpha\cap
\vare_\beta)$, 
\item[(d)$_\alpha$] $s_\alpha\in\bigcap\limits_{\gamma<\alpha}\cI_\gamma^{
[\vare_\alpha]}$. 
\end{enumerate}
So suppose that we have defined $s_\beta$'s for $\beta<\alpha$. For $\beta<
\alpha$ let 
\[\begin{array}{ll}
\cI_{\alpha,\beta}=\{s\in\bbP_{\vare_\alpha}:&\mbox{either }s_\beta\rest
\vare_\alpha\leq s,\\  
&\mbox{or } s_\beta\rest\vare_\alpha, s\mbox{ are incompatible }\}. 
  \end{array}\]
Clearly $\cI_{\alpha,\beta}\in N$ is an open dense subset of $\bbP_{
\vare_\alpha}$. Since the condition $r\rest\vare_\alpha$ is $(N,
\bbP_{\vare_\alpha})$--generic and the increasing sequence $\langle r_j
\rest\vare_\alpha:j<\lambda\rangle$ enters all open dense subsets of
$\bbP_{\vare_\alpha}$ from $N$ (by $(*)_8$), we may find $s_\alpha\in
\bigcap\limits_{\beta<\alpha}\cI_{\alpha,\beta}\cap\bigcap\limits_{j<\alpha}
\cI_j^{[\vare_\alpha]}$ such that $s_\alpha\leq r_j\rest\vare_\alpha$ for
all large enough $j<\lambda$. By (a)$_\beta$ (for $\beta<\alpha$) we
conclude that $s_\alpha$ and $s_\beta\rest\vare_\alpha$ cannot be
incompatible, and hence clauses (a)$_\alpha$--(d)$_\alpha$ are satisfied.  

Now, let conditions $s_\alpha'\in\bbP_\zeta\cap N$ (for $\alpha<\lambda$) be
such that $\Dom(s_\alpha')=\bigcup\limits_{\beta\leq\alpha}\Dom(s_\beta)$
and $s_\alpha'(\vare)$ (for $\vare\in\Dom(s_\alpha')$) is the
$<^*_\chi$--first $\bbP_\vare$--name for a condition in $\name{\bbQ}_\vare$
satisfying 
\[\begin{array}{ll}
s_\alpha'\rest\vare\forces_{\bbP_\vare}&\mbox{`` }(\forall\beta<\alpha)
(s_\beta'(\vare)\leq s_\alpha'(\vare))\ \mbox{ and}\\
&\mbox{ if there is a }\gamma\in [\alpha,\lambda)\mbox{ such that }
(\forall\beta<\alpha)(s_\beta'(\vare)\leq s_\gamma(\vare))\\
&\mbox{ then }s_\alpha'(\vare)=s_\gamma(\vare)\mbox{ for the first such
$\gamma$ ''.}
  \end{array}\]
Then the sequence $\langle s_\alpha':\alpha<\lambda\rangle$ is
$\leq_{\bbP_\zeta}$--increasing and 
\[(\forall\alpha<\lambda)(\forall\vare<\vare_\alpha)(s_\alpha\rest\vare
\forces_{\bbP_\vare}\mbox{`` }s_\alpha(\vare)=s_\alpha'(\vare)\mbox{
 ''}).\]
So it follows from (d)$_\alpha$ that for each $\vare\in u$ there is a club
$C_\vare'\subseteq\lambda$ such that $\langle s_\alpha'\rest\vare:\alpha<
\delta\rangle$ is $\bcI^{[\vare]}$--exact for all $\delta\in
C_\vare'$. Also, as clearly $(\forall\alpha<\lambda)(\exists j<\lambda)(
 s_\alpha'\leq r_j)$, we may pick a club $C^*$ of $\lambda$ such that 
\[(\forall \delta\in C^*)(\forall\alpha<\delta)(\exists j<\delta)(s_\alpha'
 \leq r_j).\]
Now, as $\bar{F}'$ is a $D$--diamond, for unboundedly many $\delta\in S\cap
\mathop{\triangle}\limits_{j<\lambda} C_j\cap\bigcap\limits_{\vare\in
u}C'_\vare\cap C^*$ we have $\langle s_\alpha':\alpha<\delta\rangle=\langle
h\circ F_\delta'(\alpha):\alpha<\delta\rangle$. Plainly, defining $F_\delta$
for those $\delta$ we had clause $(\odot)_1$ with $\gamma(\delta)\geq\zeta$
and hence $\langle s_\alpha':\alpha<\delta\rangle=\langle h^{[\zeta]}\circ
F_\delta(\alpha):\alpha<\delta\rangle$. Therefore $(\boxplus)_\delta$ holds
for those $\delta$ (remember the choice of $C^*$).  
\medskip

\noindent {\bf $\zeta$ is a limit ordinal of cofinality $\geq\lambda$}.\\
Let $\langle\vare_\alpha:\alpha<\lambda\rangle\subseteq\zeta\cap N$
be an increasing continuous sequence cofinal with $\zeta\cap N$,
$\vare_0=0$. By induction on $\alpha<\lambda$ choose conditions
$s_\alpha$ such that    
\begin{enumerate}
\item[(a)$_\alpha$] $(\exists j<\lambda)(s_\alpha\leq r_j)$,
\item[(b)$_\alpha$] $s_\alpha\in\bbP_{\vare_\alpha}\cap N$, 
\item[(c)$_\alpha$] if $\beta<\alpha<\lambda$, then $s_\beta\leq s_\alpha$, 
\item[(d)$_\alpha$] $s_\alpha\in\bigcap\limits_{\gamma<\alpha}\cI_\gamma^{
[\vare_\alpha]}$. 
\end{enumerate}
[Possible as $r\rest\vare_\alpha$ is $(N,\bbP_{\vare_\alpha})$--generic and
by $(*)_8$.] For each $\alpha<\lambda$, for some club $C_\alpha'$ of
$\lambda$ we have  
\[(\forall\delta\in C_\alpha')(\langle s_\gamma\rest\vare_\alpha:\gamma<
\delta\rangle\mbox{ is  $\bcI^{[\vare_\alpha]}$--exact\/}).\]
Also for a club $C^*$ of $\lambda$ we have $(\forall\delta\in C^*)(w_\delta
\cap\zeta\subseteq\vare_\delta)$. Like before, as $\bar{F}'$ is a
$D$--diamond, for unboundedly many $\delta\in S\cap\mathop{\triangle}
\limits_{j<\lambda} C_j\cap\mathop{\triangle}\limits_{\alpha<\lambda}
C_\alpha'\cap C^*$ we have $\langle s_\alpha:\alpha<\delta\rangle=\langle
h\circ F_\delta'(\alpha):\alpha<\delta\rangle$. Plainly, for those $\delta$
we have $\gamma(\delta)=\vare_\delta$ and also $\langle s_\alpha:\alpha<
\delta\rangle=\langle h^{[\zeta]}\circ F_\delta(\alpha):\alpha<\delta
\rangle$, and thus $(\boxplus)_\delta$ holds (remember clause
(a)$_\alpha$). 
\medskip

\noindent {\bf $\zeta$ is a successor ordinal}.\\
Like before (remember that, letting $\zeta=\zeta'+1$, the condition $r
\rest\zeta'$ is $(N,\bbP_{\zeta'})$--generic and it forces that $r(\zeta')$
is $(N[\name{G}_{\bbP_{\zeta'}}],\name{\bbQ}_{\zeta'})$--generic). 
\end{proof}
\end{proof}

\begin{remark}
In \ref{Con3.1} we may have $\bar{S}=\langle S_a:a\in W\rangle$ and
$\bar{D}=\langle D_a:a\in W\rangle$ be such that each $D_a$ is a normal
filter on $\lambda$, $S_a\in D_a^+$ satisfies the relevant demands of
\ref{Def3.2}(1), and require that there is a $D_a$ diamond $\langle
F^a_\delta: \delta\in S_a\rangle$. Then in all definitions and results
we may replace $D,S$ by $D_a,S_a$, where $a=N\cap A$. In particular, this
way we get the notions of {\em fuzzy properness over quasi
$\bar{D}$--diamonds\/} which behave nicely in iterations.
\end{remark}

\PART{B}{Building suitably proper forcing notions}

\section{A creature--free example}
In this section we show that a natural forcing notion uniformizing
colourings on ladder systems is fuzzy proper. (This forcing is a relative of
$\bbQ^*$ from \cite[4.6--4.8]{RoSh:655}.)   

Here we assume that: 

\begin{context}
\label{ConEx655}
\begin{enumerate}
\item $\lambda^*>\lambda$ is a regular cardinal, $A=\cH_{<\lambda}(
\lambda^*)$ and $W\subseteq [A]^{\textstyle\lambda}$ is as in \ref{Con3.1},
and $\lambda\subseteq a$ for each $a\in W$,   
\item $\xi^*<\lambda$, $S^*\subseteq S^{\lambda^+}_\lambda
\stackrel{\rm def}{=}\{\delta<\lambda^+:\cf(\delta)=\lambda\}$ and for
$\beta\in S^*$: 
\begin{enumerate}
\item[$(\alpha)$] $B_\beta\subseteq\beta$ is a club of $\beta$ of order
type $\otp(B_\beta)=\lambda$, 
\item[$(\beta)$] $h_\beta:B_\beta\longrightarrow\xi^*$.
\end{enumerate}
Let $\bar{B}=\langle B_\beta:\beta\in S^*\rangle$, $\bar{h}=\langle
h_\beta:\beta\in S^*\rangle$.
\end{enumerate}
\end{context}

The forcing notion $\bbQ^*=\bbQ^*(S^*,\bar{B},\bar{h})$ is defined as
follows: 
\smallskip

\noindent{\bf a condition in $\bbQ^*$}\quad is a tuple $p=(u^p,v^p,
\bar{e}^p,h^p)$ such that 
\begin{enumerate}
\item[(a)] $u^p\in [\lambda^+]^{\textstyle {<}\lambda}$, $v^p\in
[S^*]^{\textstyle {<}\lambda}$, 
\item[(b)] $\bar{e}^p=\langle e^p_\beta:\beta\in v^p\rangle$, where each
$e^p_\beta$ is a closed bounded subset of $B_\beta$, and $e^p_\beta
\subseteq u^p$, and
\item[(c)] $\sup(e^p_\beta)=\sup(u^p\cap\beta)$ (for $\beta\in v^p$), and
if $\beta_1<\beta_2$ are from $v^p$, then  
\[\sup(e^p_{\beta_2})>\beta_1\quad\mbox{ and }\quad \sup(e^p_{\beta_1})>
\sup(B_{\beta_2}\cap\beta_1),\] 
\item[(d)] $h^p:u^p\longrightarrow\xi^*$ is such that 
\[(\forall\beta\in v^p)(\forall\alpha\in e^p_\beta)(h^p(\alpha)=h_\beta(
\alpha));\]
\end{enumerate}
\noindent{\bf the order $\leq$ of $\bbQ^*$}\quad is such that $p\leq q$ if
and only if $u^p\subseteq u^q$, $h^p\subseteq h^q$, $v^p\subseteq v^q$, and
for each $\beta\in v^p$ the set $e^q_\beta$ is an end-extension of
$e^p_\beta$.
\smallskip

A tuple $p=(u^p,v^p,\bar{e}^p,h^p)$ satisfying clauses (a), (b) and (d)
above will be called {\em a pre-condition}. Note that every pre-condition
can be extended to a condition in $\bbQ^*$.

\begin{proposition}
\label{easyprop}
\begin{enumerate}
\item The forcing notion $\bbQ^*$ is $\lambda$--complete, it satisfies the 
$\lambda^+$--chain condition and $|\bbQ^*|=\lambda^+$.  
\item If $p\in\bbQ^*$, $\alpha<\lambda^+$, $\beta\in S^*$ and $\delta<
\lambda$, then there is a condition $q\geq p$ such that  
\[\alpha\in u^q,\quad \beta\in v^q\quad\mbox{ and }\quad (\forall\beta' 
\in v^q)(\otp(e^q_{\beta'})>\delta).\]
\end{enumerate}
\end{proposition}

\begin{proof}
(1)\quad Verification of the chain condition is a straightforward
application of the $\Delta$--lemma. To check that $\bbQ^*$ is
$\lambda$--complete suppose that $\langle p_i:i<j\rangle$ is a
$\leq_{\bbQ^*}$--increasing sequence of conditions from $\bbQ^*$,
$j<\lambda$. Let $r=(u^r,v^r,\bar{e}^r,h^r)$ be such that 
\[\begin{array}{ll}
v^r&=\bigcup\limits_{i<j} v^{p_i},\ \mbox{ and for }\beta\in v^r\\ 
e^r_\beta&=\bigcup\{e^{p_i}_\beta:\beta\in v^{p_i}\ \&\ i<j\}\cup
\{\sup\big(\bigcup\{e^{p_i}_\beta:\beta\in v^{p_i}\ \&\ i< j\}
\big)\}\\
u^r&=\bigcup\limits_{i<j} u^{p_i}\cup\bigcup\{e^r_\beta:\beta\in v^r\}\\ 
h^r&\supseteq\bigcup\limits_{i<j}h^{p_i},\\
  \end{array}\]
and if $\alpha\in e^r_\beta\setminus\bigcup\{e^{p_i}_\beta:\beta\in
v^{p_i},\ i<j\}$, then $h^r(\alpha)=h_\beta(\alpha)$. Using clause (c) for
$p_i$'s one easily sees that $r$ is a pre-condition. Extend it to a
condition $q\in\bbQ^*$. 
\medskip

\noindent (2)\quad Should be clear.
\end{proof}

\begin{proposition}
\label{fuzzyexample}
$\bbQ^*$ is fuzzy proper for $W$.
\end{proposition}

\begin{proof}
Suppose that $D$ is a normal filter on $\lambda$ such that there is a
$D$--diamond. We will show that $\bbQ^*$ is fuzzy proper over quasi
$D$--diamonds. First we define a $\lambda$--base $(\gR^*,\bar{\gY}^*)$ for
$\bbQ^*$ over $W$. We let $\gR^*$ be the set of all triples $(p,\delta,x)$
such that $p\in\bbQ^*$, $\delta\in\lambda$ and $x$ is a function with 
$\Dom(x)\subseteq u^p$ and $(\forall\alpha\in \Dom(x))(h^p(\alpha)=
x(\alpha))$.   

Now suppose that $a\in W$ and let $\pi_a$ be the $<^*_\chi$--first
one-to-one mapping from $a\cap\lambda^+$ to $\lambda$.  For a limit ordinal 
$\delta<\lambda$ we put
\[x_0^\delta=(\pi_a)^{-1}[\delta]\cup\{\alpha<\lambda^+:\alpha=\sup\big(
\alpha\cap (\pi_a)^{-1}[\delta]\big)\},\]
and then 
\[\gY_a^*(\delta)=\{x:\ x\mbox{ is a function from }x_0^\delta\cap a \mbox{
to }\xi^*\}.\]    
This defines $\gY_a^*$ and $\bar{\gY}^*=\langle\gY_a^*:a\in W\rangle$. It 
is easy to check that $(\gR^*,\bar{\gY}^*)$ is a $\lambda$--base for
$\bbQ^*$ (for \ref{base}(c) use repeatedly \ref{easyprop}). Assume now that  
\begin{itemize}
\item $N\prec(\cH(\chi),\in,<^*_\chi)$, $|N|=\lambda$, ${}^{<\lambda}N
\subseteq N$, $\lambda,\bbQ^*,\bar{B},\bar{h},S^*,\gR^*\in N$, and $a
\stackrel{\rm def}{=}N\cap A\in W$, $p\in\bbQ^*\cap N$,  
\item $\bcI=\langle\cI_\xi:\xi<\lambda\rangle$ lists all open dense subsets
of $\bbQ^*$ from $N$,  
\item $h:\lambda\longrightarrow N$ satisfies $\bbQ^*\cap N\subseteq\Rang(
h)$, and 
\item $\bar{F}=\langle F_\delta:\delta\in S\rangle$ is a quasi $D$--diamond
for $(N,h,\bbQ^*)$ and $\bar{q}$ is a fuzzy candidate over $\bar{F}$.  
\end{itemize}
For limit $\delta\in S$ let $\cY(\delta)=\cY(N,\bbQ^*,h,\bar{F},\gR^*,
\bar{\gY}^*,\delta)$ be as defined in \ref{pre1.2}(3) (and thus $\bar{q}=
\langle q_{\delta,x}:\delta\in S\mbox{ is limit }\&\ x\in\cY(\delta)
\rangle$). Also let $E_0$ be the set of all $\delta<\lambda$ which are
limits of members of $\lambda\setminus S$ (so it is a club of $\lambda$). 

We are going to show that the condition $p$ is $(\gR^*,\bar{\gY}^*)$--fuzzy
generic for $\bar{q}$. Note that, as $\bbQ^*$ satisfies the $\lambda^+$--cc,
the condition $p$ is $(N,\bbQ^*)$--generic (in the standard sense). So, by
\ref{1.2A}(3), it is enough to give a strategy of the Generic player in the
game $\Gfl(p,N,\bcI,h,\bbQ^*,\bar{F},\bar{q})$ which guarantees that the
result $\langle r_i,C_i:i<\lambda\rangle$ of the play satisfies
\ref{pre1.2}(5)($\beta$). 

Suppose that we arrive to a stage $\delta\in S$ and $\langle r_i,C_i:i<
\delta\rangle$ is the sequence played so far. First, Generic picks the
$<^*_\chi$--first condition $r_\delta'$ stronger than all $r_i$'s played so
far and such that 
\smallskip

{\em if}\qquad $\delta$ is limit and $\big(\exists x\in \cY(\delta)\big)
\big(\exists r\in\bbQ^*\big)\big(q_{\delta,x}\leq r\ \&\ (\forall i<\delta) 
(r_i\leq r)\big)$,

{\em then\/} $q_{\delta,x}\leq r_\delta'$ for some $x\in \cY(\delta)$.
\smallskip

\noindent Then she plays the $<^*_\chi$--first condition $r_\delta$ above
$r_\delta'$ such that  
\begin{enumerate}
\item[$(*)_1$] if $\beta\in v^{r_\delta}$, then $\otp(e^{r_\delta}_\beta)>
\delta$, and  
\item[$(*)_2$] $(\pi_a)^{-1}[\delta]\subseteq u^{r_\delta}$ and $(\pi_a)^{
-1}[\delta]\cap S^*\subseteq v^{r_\delta}$. 
\end{enumerate}
The set $C_\delta$ played at this stage is $(\alpha,\lambda)\cap E_0$, where 
$\alpha$ is the first ordinal such that 
\begin{enumerate}
\item[$(*)_3$] $\pi_a[u^{r_\delta}\cap N]\subseteq\alpha$, and the set 
\[\{q\in\bbQ^*:(\pi_a)^{-1}[\delta]\subseteq u^q\ \&\
(\pi_a)^{-1}[\delta]\cap S^*\subseteq v^q\ \&\ (\forall\beta\in v^q)(\otp(
e^q_\beta)>\delta)\}\]
(which is an open dense subset of $\bbQ^*$ from $N$; remember
\ref{easyprop}) is in $\{\cI_\xi:\xi<\alpha\}$,
\item[$(*)_4$] $\otp(B_\beta\cap (\sup(e^{r_\delta}_\beta)+1))<\alpha$ for
all $\beta\in v^{r_\delta}$,  
\item[$(*)_5$] if $\beta\in v^{r_\delta}$ and $a\cap\beta\setminus (\sup(
e^{r_\delta}_\beta)+1)\neq\emptyset$, then there is $\gamma\in a\cap\beta
\setminus(\sup(e^{r_\delta}_\beta)+1)$ with $\pi_a(\gamma)<\alpha$.
\end{enumerate}

Why does this strategy work (i.e., why does it ensure
\ref{pre1.2}(5)($\beta$))?\\ 
Let $\langle r_i,C_i:i<\lambda\rangle$ be a play according to this strategy,
and suppose that $\delta\in S\cap \mathop{\triangle}\limits_{i<\lambda} C_i$
is a limit ordinal such that $\langle h\circ F_\delta(\alpha):\alpha<\delta
\rangle$ is a $\leq_{\bbQ^*}$--increasing $\bcI$--exact sequence of
conditions from $\bbQ^*\cap N$ such that $(\forall\alpha<\delta)(\exists i<
\delta)(h\circ F_\delta(\alpha)\leq r_i)$. Note that then  
\begin{enumerate}
\item[$(*)_6$] if $\beta\in\bigcup\limits_{i<\delta} v^{r_i}$, then
$\otp(\bigcup\limits_{i<\delta} e^{r_i}_\beta)=\delta$ and 
$\bigcup\limits_{i<\delta} e^{r_i}_\beta$ is an unbounded subset of
$\{\varepsilon\in B_\beta:\otp(\varepsilon\cap B_\beta)<\delta\}$, and 
\item[$(*)_7$] $\bigcup\limits_{i<\delta} u^{r_i}\cap N=(\pi_a)^{-1}[
\delta]=\bigcup\limits_{\alpha<\delta} u^{h\circ F_\delta(\alpha)}$ and 
$\bigcup\limits_{i<\delta} v^{r_i}\cap N=(\pi_a)^{-1}[\delta]\cap S^*=
\bigcup\limits_{\alpha<\delta} v^{h\circ F_\delta(\alpha)}$, 
\item[$(*)_8$] if $\beta\in(\pi_a)^{-1}[\delta]\cap S^*$, then 
\[\bigcup\{e^{h\circ F_\delta(\alpha)}_\beta:\alpha<\delta\ \&\ \beta\in
v^{h\circ F_\delta(\alpha)}\}=\bigcup\{e^{r_i}_\beta:i<\delta\ \&\ \beta\in 
v^{r_i}\}.\]
\end{enumerate}
We want to show that 
\begin{enumerate}
\item[$(\boxdot)$] for some $x\in\cY(\delta)$, there is a common upper
bound to $\{r_i:i<\delta\}\cup\{q_{\delta,x}\}$ 
\end{enumerate}
(which, by the definition of our strategy, will finish the proof). For
$\beta\in S^*$ let $\gamma_\beta \in B_\beta$ be such that $\otp(B_\beta
\cap\gamma_\beta)=\delta$. Now, let a pre-condition $r'=(u^{r'},v^{r'},
\bar{e}^{r'},h^{r'})$ be such that 
\begin{itemize}
\item $v^{r'}=\bigcup\limits_{i<\delta} v^{r_i}$, $u^{r'}=\bigcup\limits_{i<
\delta}u^{r_i}\cup\{\gamma_\beta:\beta\in v^{r'}\}$, 
\item $e^{r'}_\beta=\bigcup\{e^{r_i}_\beta:i<\delta\ \&\ \beta\in v^{r_i}\}
\cup\{\gamma_\beta\}$  (for $\beta\in v^{r'}$), and 
\item $h^{r'}:u^{r'}\longrightarrow\xi^*$ is such that $\bigcup\limits_{i<
\delta} h^{r_i}\subseteq h^{r'}$ and $h^{r'}(\gamma_\beta)=h_\beta(
\gamma_\beta)$. 
\end{itemize}
One easily verifies that the above conditions indeed define a pre-condition
(remember $(*)_6$). Also, note that if $\beta\in v^{r'}$, then $\gamma_\beta 
\in e^{r'}_\beta\setminus\bigcup\{e^{r_i}_\beta:i<\delta\ \&\ \beta\in
v^{r_i}\}$ and each $e^{r_i}_\beta$ is a proper subset of $\bigcup\{
e^{r_i}_\beta:i<\delta\ \&\ \beta\in v^{r_i}\}$. Moreover, if $\beta\in
v^{r'}$ and $\gamma_\beta\in N$, then $\gamma_\beta=\sup(u^{r'}\cap N\cap
\gamma_\beta)=\sup((\pi_a)^{-1}[\delta]\cap\gamma_\beta)$ (by $(*)_5+(*)_2$; 
remember also $(*)_7$). Now, extend $r'$ to a pre-condition $r''$ such that
$u^{r''}=u^{r'}\cup x^\delta_0$, $v^{r''}=v^{r'}$ and $e^{r''}_\beta=
e^{r'}_\beta$ for $\beta\in v^{r''}$ (clearly possible). Let $x=h^{r''}\rest
x^\delta_0$ (note that $x^\delta_0\subseteq a$).  Since $r''$ is stronger
than all $h\circ F_\delta(\alpha)$ (for $\alpha<\delta$), any condition
stronger than $r''$ witnesses that $x\in\cY(\delta)$. Now we put
\begin{itemize}
\item $u^*=u^{q_{\delta,x}}\cup u^{r''}$, $v^*=v^{q_{\delta,x}}\cup v^{r''}$,
$h^*= h^{q_{\delta,x}}\cup h^{r''}$, 
\item if $\beta\in v^{q_{\delta,x}}$, then $e^*_\beta=e^{q_{\delta,x}
}_\beta$, and if $\beta\in v^{r''}\setminus N$, then
$e^*_\beta=e^{r''}_\beta$.
\end{itemize}
Note that $h^*$ is a function from $u^*$ to $\xi^*$ by $(*)_7$ (remember the
choice of $x$ and that $q_{\delta,x}\in N$ is stronger than all $h\circ
F_\delta(\alpha)$'s). Also, if $\beta\in v^{r''}\cap N$ then ($\beta\in
v^{q_{\delta,x}}$ and) $e^{q_{\delta,x}}_\beta$ is an end-extension of
$e^{r''}_\beta$ (remember $(*)_7+(*)_8$). Hence $(u^*,v^*,\bar{e}^*,h^*)$ is
a pre-condition stronger than both $q_{\delta,x}$ and $r''$. Extending it to
a condition in $\bbQ^*$ we conclude $(\boxdot)$, thus completing the proof
of \ref{fuzzyexample}. 
\end{proof}

\begin{corollary}
Assume that $\lambda$ is a strongly inaccessible cardinal, $2^\lambda=
\lambda^+$, $2^{\lambda^+}=\lambda^{++}$ and $D$ is a normal filter on
$\lambda$ such that there is a $D$--diamond. Then there is a forcing 
notion $\bbP$ such that:
\begin{itemize}
\item $\bbP$ is $\lambda$--complete weakly fuzzy proper over quasi
$D$--diamonds for $W$ and it satisfies the $\lambda^{++}$--cc,
\item in $\bV^{\bbP}$, $2^\lambda=2^{\lambda^+}=\lambda^{++}$ and for every
$\xi^*,S^*,\bar{B},\bar{h}$ as in \ref{ConEx655}(2) there is $h:\lambda^{++}
\longrightarrow\xi^*$ such that for every $\beta\in S^*$ the set  
\[\{\alpha\in B_\beta:h_\beta(\alpha)=h(\alpha)\}\]
contains a club.  
\end{itemize}
\end{corollary}

\begin{proof}
The forcing notion $\bbP$ is the limit $\bbP_{\lambda^{++}}$ of a
$\lambda$--support iteration $\bar{\bbQ}=\langle\bbP_\alpha,
\name{\bbQ}_\alpha:\alpha<\lambda^{++}\rangle$, where each
$\name{\bbQ}_\alpha$ is forced to be $\bbQ^*(S^{
\lambda^+}_\lambda,\name{\bar{B}}_\alpha,\name{\bar{h}}_\alpha)$ for some  
$\name{\bar{B}}_\alpha,\name{\bar{h}}_\alpha$. Then, by \ref{lppcc},
\ref{1.6} and \ref{fuzzyexample} we are sure that $\bbP$ satisfies the 
$\lambda^{++}$--cc, it is weakly fuzzy proper over quasi $D$--diamonds
for $W$ and it has a dense subset of size $\lambda^{++}$. Consequently we
may arrange suitable bookkeeping to take care of all
$\bbP_{\zeta^++}$--names $\name{\bar{B}},\name{\bar{h}}$ for objects as in
\ref{ConEx655}(2) -- the details and the rest should be clear. 
\end{proof}

\section{Trees and creatures} 
Let us introduce the notation used in the forcing notions we want to
build. The terminology here is somewhat parallel to that of \cite[\S 1.2,
1.3]{RoSh:470}, but there are some differences as the context is
different. We start with the tree case. 

\begin{definition}
\label{treecreature}
Let $\bH:\lambda\longrightarrow\cH(\lambda^+)$.
\begin{enumerate}
\item 
A {\em $\lambda$--tree creature for $\bH$} is a tuple
\[t=(\eta,\dis,\pos,\nor)=(\eta[t],\dis[t],\pos[t],\nor[t])\]
such that $\dis\in\cH(\lambda^+)$, $\nor\in\lambda+1$,
\[\eta\in\bigcup\limits_{\alpha<\lambda}\prod\limits_{\beta<\alpha}\bH(
\beta),\quad \mbox{ and }\quad \emptyset\neq\pos\subseteq\{\nu\in
\bigcup\limits_{\alpha<\lambda}\prod\limits_{\beta<\alpha}\bH(\beta):\eta
\vtl\nu\}.\] 
$\TCR[\bH]$ is the family of all $\lambda$--tree creatures for $\bH$. 

For $\eta\in\bigcup\limits_{\alpha<\lambda}\prod\limits_{\beta<\alpha} 
\bH(\beta)$ we let $\TCR_\eta[\bH]=\{t\in\TCR[\bH]: \eta[t]=\eta\}$.

\item Let $K\subseteq\TCR[\bH]$. A {\em tree--composition operation on $K$}
is a mapping $\Sigma$ with values in $\cP(K)$ and the domain consisting of
systems $\langle t_\nu:\nu\in \hat{T}\rangle$ such that
\begin{itemize}
\item $T$ is a complete $\lambda$--quasi tree of height $\rht(T)<\lambda$,
\item for each $\nu\in \hat{T}$, $t_\nu\in K$ satisfies $\suc_T(\nu)=
\pos[t_\nu]$,
\end{itemize}
and 
\begin{itemize}
\item if $t\in\Sigma(t_\nu:\nu\in\hat{T})$, then $\eta[t]=\mrot(T)$ and
$\pos[t]\subseteq\max(T)$, 
\item if $t\in\Sigma(t_\nu:\nu\in\hat{T})$ and $t_\nu\in\Sigma(s^\nu_\rho:
\rho\in\hat{T}_\nu)$ (for $\nu\in\hat{T}$), then $t\in\Sigma(s^\nu_\rho:
\rho\in\bigcup\limits_{\nu\in \hat{T}} \hat{T}_\nu)$, and 
\item for each $t\in K$ we have $\langle t\rangle\in\dom(\Sigma)$ and
$t\in\Sigma(t)$. 
\end{itemize}
Then $(K,\Sigma)$ is called {\em a $\lambda$--tree creating pair} (for
$\bH$).  
\item A $\lambda$--tree creating pair $(K,\Sigma)$ is {\em local} if
\begin{itemize}
\item $(t_\nu:\nu\in T)\in\dom(\Sigma)$ implies $\rht(T)=\lh(\mrot(T))+1$
  (and so $T=\{\mrot(T)\}\cup\pos[t_{\mrot(T)}]$), and 
\item $t'\in\Sigma(t)$ implies $\nor[t']\leq\nor[t]$.
\end{itemize}
We say that $(K,\Sigma)$ is {\em very local} if, additionally, for every
$\nu\in\bigcup\limits_{\alpha<\lambda}\prod\limits_{\beta<\alpha}\bH(\beta)$
such that $K\cap\TCR_\nu[\bH]\neq\emptyset$ there is $t^*_\nu\in K\cap
\TCR_\nu[\bH]$ satisfying $(\forall t\in K\cap\TCR_\nu[\bH])(t\in
\Sigma(t^*_\nu))$. The tree creature $t^*_\nu$ may be called {\em the
minimal creature at $\nu$}. 
\item If $(K,\Sigma)$ is a very local $\lambda$--tree creating pair, then
{\em the minimal tree $T^*$ for $(K,\Sigma)$} and {\em the minimal condition
$p^*$ for $(K,\Sigma)$} are defined by
\[\begin{array}{lll}
T^*&=T^*(K,\Sigma)&=\{\eta\in\bigcup\limits_{\alpha<\lambda}\prod\limits_{
\beta<\lambda}\bH(\beta): (\forall\alpha<\lh(\eta))(\eta\rest(\alpha+1)\in
\pos[t^*_{\eta\rest\alpha}])\}\\ 
p^*&=p^*(K,\Sigma)&=\langle t^*_\nu:\nu\in T^*\rangle.
  \end{array}\]
[Note that, in the general case, $T^*$ could be empty, but in real
applications this can be easily avoided.]
\end{enumerate}
\end{definition} 

\begin{definition}
\label{treeforcing} 
Let $(K,\Sigma)$ be a $\lambda$--tree creating pair for $\bH$.
\begin{enumerate}
\item We define the forcing notion $\bbQ^{\tree}_1(K,\Sigma)$ by:

\noindent{\bf conditions } are systems $p=\langle t_\eta:\eta\in T\rangle$  
such that 
\begin{enumerate}
\item[(a)] $\emptyset\neq T\subseteq\bigcup\limits_{\alpha<\lambda}
\prod\limits_{\beta<\alpha}\bH(\beta)$ is a complete $\lambda$--quasi tree
with $\max(T)=\emptyset$,  
\item[(b)] $t_\eta\in\TCR_\eta[\bH]\cap K$ and $\pos[t_\eta]=\suc_T(\eta)$,
\item[(c)${}_1$] for every $\eta\in\llim(T)$, $\lim(\nor[t_{\eta\rest
\alpha}]:\alpha<\lambda,\, \eta\rest\alpha\in T)=\lambda$;  
\end{enumerate}
\noindent{\bf the order} be given by:

\noindent $\langle t^1_\eta: \eta\in T^1\rangle\leq\langle t^2_\eta: \eta\in
T^2\rangle$ if and only if 

\noindent $T^2\subseteq T^1$ and for each $\eta\in T^2$ there is a complete 
$\lambda$--quasi tree $T_{0,\eta}\subseteq (T^1)^{[\eta]}$ of height
$\rht(T_{0,\eta})<\lambda$ such that $t^2_\eta\in\Sigma(t^1_\nu:\nu\in
\hat{T}_{0,\eta})$. 

If $p=\langle t_\eta:\eta\in T\rangle$ then we write $\mrot(p)=\mrot(T)$,
$T^p= T$, $t^p_\eta = t_\eta$ etc.
  
\item Let $D^*$ be a filter on $\lambda$. The forcing notion
$\bbQ^\tree_{D^*}(K,\Sigma)$ is defined similarly, replacing the condition 
(c)$_1$ by 
\begin{enumerate}
\item[(c)${}_{D^*}$] for some set $Y=Y^p\in D^*$ we have 
\[(\forall \delta\in Y)(\forall\eta\in (T)_\delta)(\nor[t_\eta]\geq
|\delta|).\]  
[The set $Y^p$ above may be called {\em a witness for $p\in\bbQ^\tree_{
D^*}(K,\Sigma)$}.] 
\end{enumerate}
\item The forcing notion $\bbQ^\tree_{\cl}(K,\Sigma)$ is defined by
replacing the condition (c)$_1$ by 
\begin{enumerate}
\item[(c)$_{\cl}$] 
\begin{enumerate}
\item[$(\alpha)$] $(\forall \eta\in T)(\exists\nu\in T)(\eta\vtl\nu\ \&\
\nor[t_\nu]\geq |\lh(\nu)|)$, and 
\item[$(\beta)$] $(\forall \eta\in T)(\nor[t_\eta]=0\ \mbox{ or }\
\nor[t_\eta]\geq |\lh(\eta)|)$, and 
\item[$(\gamma)$] $\nor[t_{\mrot(p)}]\geq|\lh(\mrot(p))|$, and 
\item[$(\delta)$] if $\delta<\lambda$ is a limit ordinal and $\langle\eta_i:
i<\delta\rangle\subseteq T$ is a $\vtl$--increasing sequence such that
$\nor[t_{\eta_i}]\geq |\lh(\eta_i)|$ for each $i<\delta$ and $\eta=
\bigcup\limits_{i<\delta}\eta_i$, then ($\eta\in T$ and) $\nor[t_\eta]\geq
|\lh(\eta)|$. 
\end{enumerate}
\end{enumerate}

\item If $p\in\bbQ^{\tree}_e(K,\Sigma)$ and $\eta\in T^p$, then we let
$p^{[\eta]}=\langle t^p_\nu:\nu\in (T^p)^{[\eta]}\rangle$.
\item For the sake of notational convenience we define partial order
$\bbQ^{\tree}_{\emptyset}(K,\Sigma)$ in the same manner as
$\bbQ^{\tree}_e(K,\Sigma)$ above but we omit the requirement (c)$_e$.
\end{enumerate}
\end{definition}

\begin{definition}
\label{order}
Let $(K,\Sigma)$ be a $\lambda$--tree creating pair for $\bH$, $t\in K$. We
define a relation $\preceq^t_\Sigma$ on $\Sigma(t)$ by

$t'\preceq^t_\Sigma t''$\quad if and only if\quad ($t',t''\in\Sigma(t)$ and)
$t''\in\Sigma(t')$. 

\noindent If $(K,\Sigma)$ is very local, $t^*_\nu$ is the minimal creature
at $\nu$, then $\preceq^{t^*_\nu}_\Sigma$ is also denoted by
$\preceq^\nu_\Sigma$.  
\end{definition}

\begin{remark}
\begin{enumerate}
\item Note that the relation $\preceq^t_\Sigma$ is transitive and
reflexive. 
\item If $(K,\Sigma)$ is local and $p\in\bbQ^\tree_\emptyset(K,\Sigma)$,
then $T^p$ is a complete $\lambda$--tree.
\end{enumerate}
\end{remark}

Now we are going to describe the non-tree case of forcing with
creatures. For sake of simplicity we restrict ourselves to what corresponds
to forgetful creatures of \cite[1.2.5]{RoSh:470}. 

\begin{definition}
\label{creature}
Let $\bH:\lambda\longrightarrow\cH(\lambda^+)$.
\begin{enumerate}
\item A {\em forgetful $\lambda$--creature for $\bH$} is a tuple
\[t=(\ad,\au,\dis,\val,\nor)=(\ad[t],\au[t],\dis[t],\val[t],\nor[t])\]
such that $\dis\in\cH(\lambda^+)$, $\nor\in\lambda+1$, $\ad<\au<\lambda$ and 
$\emptyset\neq\val\subseteq\prod\limits_{\ad\leq\beta<\au}\bH(\beta)$.

$\CR[\bH]$ is the family of all forgetful $\lambda$--creatures for $\bH$. 

Since we will consider only forgetful $\lambda$--creatures, from now on we
will omit the adjective ``forgetful''. 
\item Let $K\subseteq\CR[\bH]$. A {\em composition operation on $K$} is a
mapping $\Sigma$ with values in $\cP(K)$ and the domain consisting of
systems $\langle t_i:i<j\rangle\subseteq K$ such that $j<\lambda$ and 
\[\begin{array}{l}
\au[t_i]=\ad[t_{i+1}]\quad\mbox{ for $i<i+1<j$, and}\\
\sup\{\au[t_{i'}]:i'<i\}=\ad[t_i]\quad\mbox{ for limit }i<j,
  \end{array}\]
and if $t\in\Sigma(t_i:i<j)$, then 
\begin{itemize}
\item $\alpha^-=\ad[t]=\ad[t_0]$, $\alpha^+=\au[t]=\sup\{\au[t_i]:i<j\}$,
and  
\item $\val[t]\subseteq\{\nu\in\prod\limits_{\alpha^-\leq\beta<\alpha^+}\bH(
\beta):(\forall i<j)(\nu\rest [\ad[t_i],\au[t_i])\in\val[t_i])\}$,
\end{itemize}
and 
\begin{itemize}
\item if $t_i\in\Sigma(s^i_\zeta:\zeta<\zeta_i)$ (for $i<j$) and $t\in\Sigma 
(t_i:i<j)$,\\
then $t\in\Sigma(s^i_\zeta:\zeta<\zeta_i,\ i<j)$, and 
\item for each $t\in K$ we have $\langle t\rangle\in\dom(\Sigma)$ and
$t\in\Sigma(t)$. 
\end{itemize}
Then $(K,\Sigma)$ is called {\em a $\lambda$--creating pair} (for $\bH$).  
\item We say that $(K,\Sigma)$ is {\em local\/} if for each $t\in K$
\begin{itemize}
\item $\au[t]=\ad[t]+1$ , and 
\item $t'\in\Sigma(t)$ implies $\nor[t']\leq\nor[t]$.
\end{itemize}
It is {\em very local\/} if, additionally, for each $\alpha<\lambda$ there
is $t^*_\alpha\in K$ such that $\ad[t^*_\alpha]=\alpha$ and for every $t\in
K$ with $\ad[t]=\alpha$ we have $t\in\Sigma(t^*_\alpha)$. The creature
$t^*_\alpha$ will be called {\em the minimal creature $t^*_\alpha$ at
  $\alpha$}.  
\item For $j<\lambda$, {\em a $j$--approximation for $(K,\Sigma)$\/} is a
pair $(w,\langle t_i:i<j\rangle)$ such that $t_i\in K$, 
\[\begin{array}{l}
\au[t_i]=\ad[t_{i+1}]\quad\mbox{ for $i<i+1<j$, and}\\
\sup\{\au[t_{i'}]:i'<i\}=\ad[t_i]\quad\mbox{ for limit }i<j,
  \end{array}\]
and $w\in\prod\limits_{\alpha<\ad[t_0]}\bH(\alpha)$.
\item For a $j$--approximation $(w,\langle t_i:i<j\rangle)$ for $(K,\Sigma)$
we let 
\[\begin{array}{ll}
\pos(w,\langle t_i:i<j\rangle)=\{v\in\prod\limits_{\alpha<\alpha^+}\bH(
\alpha):&
w\vtl v\ \mbox{ and for all }i<j\\
&v\rest [\ad[t_i],\au[t_i])\in\val[t_i]\ \},
  \end{array}\] 
where $\alpha^+=\sup\{\au[t_i]:i<j\}$.
\end{enumerate}
\end{definition}

\begin{definition}
\label{forcing} 
Let $(K,\Sigma)$ be a $\lambda$--creating pair for $\bH$.
\begin{enumerate}
\item We define the forcing notion $\bbQ^*_1(K,\Sigma)$:

\noindent{\bf conditions } are pairs $p=(w,\bar{t})$ such that 
\begin{enumerate}
\item[(a)] $\bar{t}=\langle t_i:i<\lambda\rangle$ is a sequence of
$\lambda$--creatures from $K$ satisfying
\[\begin{array}{ll}
\au[t_i]=\ad[t_{i+1}]&\quad\mbox{ for $i<i+1<\lambda$, and}\\
\sup\{\au[t_{i'}]:i'<i\}=\ad[t_i]&\quad\mbox{ for limit }i<\lambda,
  \end{array}\]
\item[(b)] $w\in \prod\limits_{\alpha<\ad[t_0]}\bH(\alpha)$
\item[(c)${}_1$] $\lim(\nor[t_i]:i<\lambda)=\lambda$ 
\end{enumerate}
\noindent{\bf the order} is given by:

\noindent $(w^1,\langle t^1_i:i<\lambda\rangle)\leq (w^2,\langle t^2_i: i<
\lambda\rangle$ if and only if  

\noindent for some continuous strictly increasing sequence $\langle i_\zeta:
\zeta<\lambda\rangle$ we have 
\[w^2\in\pos(w^1,\langle t^1_i:i<i_0\rangle)\quad\mbox{and}\quad (\forall
\zeta<\lambda)\big(t^2_\zeta\in\Sigma(t^1_i:i_\zeta\leq i<i_\zeta)\big).\]  

If $p=(w,\langle t_i:i<\lambda\rangle)$, then we write $w^p=w$, $t^p_i=t_i$
(for $i<\lambda$).
  
\item Let $D^*$ be a filter on $\lambda$. The forcing notion $\bbQ^*_{D^*}(
K,\Sigma)$ is defined similarly, replacing the condition (c)$_1$ by  
\begin{enumerate}
\item[(c)${}_{D^*}$] for some set $Y=Y^p\in D^*$ we have 
\[(\forall i\in Y)(\nor[t_i]\geq |\ad[t_i]|).\]
[The set $Y^p$ above may be called {\em a witness for $p\in\bbQ^*_{D^*}(K,
\Sigma)$}.] 
\end{enumerate}
\item For the sake of notational convenience we define partial order
$\bbQ^*_{\emptyset}(K,\Sigma)$ in the same manner as $\bbQ^*_e(K,\Sigma)$
above but we omit the requirement (c)$_e$. If $(K,\Sigma)$ is very local,
then {\em the minimal condition $p^*$ for $(K,\Sigma)$} is 
\[p^*=p^*(K,\Sigma)=(\langle\rangle,\langle t^*_\alpha:\alpha<\lambda
\rangle)\in\bbQ^*_\emptyset(K,\Sigma),\]
where $t^*_\alpha$ is the minimal creature at $\alpha$.
\item The relations $\preceq^t_\Sigma$ and $\preceq^{t^*_\alpha}_\Sigma=
\preceq^\alpha_\Sigma$ are defined in a way parallel to \ref{order}. 
\end{enumerate}
\end{definition}

\section{Getting completeness and bounding properties}
In this section we introduce properties of $\lambda$--tree creating pairs
ensuring that the resulting forcing notions are complete or strategically
complete. Next we show that adding bounds on the size of $\bH(\alpha)$
guarantees strong bounding properties from section A.2. Finally we will
introduce parallel completeness conditions for the case of
$\lambda$--creating pairs. 

\begin{definition}
\label{treecomp}
Let $(K,\Sigma)$ be a $\lambda$--tree creating pair for $\bH$, $\kappa$ be a
cardinal (and $\lambda,\bar{\lambda}$ be as in \ref{incon}).
\begin{enumerate}
\item We say that a $\lambda$--tree creature $t\in K$ is {\em
$\kappa$--complete} (for $(K,\Sigma)$) if 
\begin{enumerate}
\item[$(\alpha)$] for every $\preceq^t_\Sigma$--increasing chain $\langle
t_\alpha:\alpha<\delta\rangle\subseteq \Sigma(t)$ with $\delta<\kappa$ and
$\nor[t_\alpha]>0$, there is $t_\delta\in \Sigma(t)$ such that $(\forall
\alpha<\delta)(t_\alpha\preceq^t_\Sigma t_\delta)$ and $\nor[t_\delta]\geq
\min\{\nor[t_\alpha]:\alpha<\delta\}$,  
\item[$(\beta)$] if $t'\in\Sigma(t)$, $\nor[t']=0$, then $|\pos[t']|=1$ and
$\Sigma(t')=\{t'\}$,  
\item[$(\gamma)$] if $\nu\in \pos[t]$, then there is $t'\in\Sigma(t)$ such
that $\pos[t']=\{\nu\}$ and $\nor[t']=0$.
\end{enumerate}
\item $t\in K$ is said to be {\em exactly $\kappa$--complete} if it is
$\kappa$--complete and 
\begin{enumerate}
\item[$(\otimes)$] if $\bar{t}=\langle t_\alpha:\alpha<\kappa\rangle
\subseteq\Sigma(t)$ is a strictly $\preceq^t_\Sigma$--increasing chain, then
$\bar{t}$ has no $\preceq^t_\Sigma$--upper bound in $\Sigma(t)$, but
$\bigcap\limits_{\alpha<\kappa}\pos[t_\alpha]\neq\emptyset$.
\end{enumerate}
\item We say that $(K,\Sigma)$ is {\em $\bar{\lambda}$--complete\/} ({\em
exactly $\bar{\lambda}$--complete}, respectively) if 
\begin{enumerate}
\item[(a)] $(K,\Sigma)$ is very local, and
\item[(b)$_{\bar{\lambda}}$] each minimal creature $t^*_\nu$ is
$\lambda^+_{\lh(\nu)}$--complete (exactly $\lambda_{\lh(\nu)}$--complete,
resp.). 
\end{enumerate}
We say that $(K,\Sigma)$ is {\em just $\lambda$--complete\/} if it
satisfies (a) above and 
\begin{enumerate}
\item[(b)$^\lambda$] each minimal creature $t^*_\nu$ is $\lambda$--complete.
\end{enumerate} 
\end{enumerate}
\end{definition}

\begin{proposition}
\label{getcomp}
Assume that $(K,\Sigma)$ is a very local $\lambda$--tree creating pair for
$\bH$, $D^*$ is a ${<}\lambda$--complete filter on $\lambda$. Let $\bbP$ be
one of the forcing notions $\bbQ^\tree_1(K,\Sigma)$, $\bbQ^\tree_{D^*}(K,
\Sigma)$, or $\bbQ^\tree_{\cl}(K,\Sigma)$. 
\begin{enumerate}
\item If $(K,\Sigma)$ is $\bar{\lambda}$--complete, then $\bbP$ is
  strategically $({<}\lambda)$--complete.   
\item If $(K,\Sigma)$ is just $\lambda$--complete, then $\bbP$ is
  $\lambda$--complete. 
\item If $(K,\Sigma)$ is exactly $\bar{\lambda}$--complete, then $\bbP$ is 
$\lambda$--complete. 
\end{enumerate}
\end{proposition}

\begin{proof}
(1)\qquad Let $\bbP\in\{\bbQ^\tree_1(K,\Sigma),\bbQ^\tree_{D^*}(K,\Sigma),
\bbQ^\tree_{\cl}(K,\Sigma)\}$ and let $r\in\bbP$. Consider the following 
strategy $\st$ of Complete in the game $\Game_0^\lambda(\bbP,\emptyset,r)$. 
At stage $j<\lambda$ of the game, after a sequence $\langle (p_i,q_i):i<j
\rangle\conc\langle p_j\rangle\subseteq\bbP$ has been constructed (so $p_j$ is
the move of Incomplete), she plays the $<^*_\chi$--first condition $q_j\in
\bbP$ stronger than $p_j$ and such that $\lh(\mrot(q_j))>j+\omega$. 

Why is this a winning strategy? Suppose that the players have arrived at a
limit stage $\delta<\lambda$ of the game, Complete has used $\st$ and
$\langle (p_i,q_i):i<\delta\rangle$ is the result of the game so far. Our
aim is to show that the (increasing) sequence $\langle q_i:i<\delta\rangle$
has an upper bound in $\bbP$. To this end we are going to define a condition
$q=\langle t_\eta:\eta\in T\rangle\in\bbP$ inductively defining $(T)_\alpha$
and $t_\eta$ for $\alpha<\lambda$, $\eta\in (T)_\alpha$. First we declare   
$\mrot(T)=\bigcup\limits_{i<\delta}\mrot(q_i)$ and we note that 
\[\mrot(T)\in\bigcap\limits_{i<\delta}T^{q_i}\quad\mbox{and}\quad \delta\leq
\lh(\mrot(T))<\lambda_{\lh(\mrot(T))}^+.\] 
Now we may choose $t_{\mrot(T)}\in\TCR_{\mrot(T)}[\bH]$ so that 
$t_{\mrot(T)}\in\Sigma(t_{\mrot(T)}^{q_i})$ (for all $i<\delta$) and 
$\nor[t_{\mrot(T)}]\geq\min\{\nor[t_{\mrot(T)}^{q_i}]:i<\delta\}$, and we
declare $\pos[t_{\mrot(T)}]\subseteq T$ (thus defining
$(T)_{\lh(\mrot(T))+1}$). Next we proceed inductively in a similar manner:
suppose that $(T)_\alpha$ has been already defined and it is included in
$\bigcap\limits_{i<\delta}T^{q_i}$. For each $\eta\in (T)_\alpha$ choose 
$t_\eta$ such that  
\[(\forall i<\delta)(t_\eta\in\Sigma(t_\eta^{q_i}))\quad\mbox{and}\quad
\nor[t_\eta]\geq\min\{\nor[t_\eta^{q_i}]:i<\delta\},\]
and declare $\pos[t_\eta]\subseteq T$. (So after this step $(T)_{\alpha+1}$
is defined.) If $\alpha<\lambda$ is limit and $(T)_\beta$ has been defined
for $\beta<\alpha$, then we let $(T)_\alpha$ consist of all $\eta$ such that
$\eta\rest\beta\in (T)_\beta$ whenever $\lh(\mrot(T))\leq\beta<\alpha$, and
then we choose $t_\eta$ (for $\eta\in (T)_\alpha$ like above).  

This way we build a condition $\langle t_\eta:\eta\in T\rangle\in
\bbQ^\tree_\emptyset(K,\Sigma)$, and it is very straightforward to verify
that this condition is actually in $\bbP$, and it is stronger than all $q_i$
(for $i<\delta$).
\medskip

\noindent (2), (3)\qquad Similar.
\end{proof}

The exact $\bar{\lambda}$--completeness may seem to be very strange and/or
strong. But as a matter of fact it is easy to modify any
$\bar{\lambda}$--complete $\lambda$--tree creating pair to one that is
exactly complete (and the respective forcing notions are very close).

\begin{definition}
\label{exactivity}
Let $(K,\Sigma)$ be a very local $\bar{\lambda}$--complete $\lambda$--tree
creating pair for $\bH$. We define {\em the $\bar{\lambda}$--exactivity
$(K^\exl,\Sigma^\exl)$ of $(K,\Sigma)$} as follows.

Let $\eta\in \bigcup\limits_{\alpha<\lambda}\prod\limits_{\beta<\alpha}
\bH(\beta)$. We let $K^\exl\cap \TCR_\eta[\bH]$ consist of all
$\lambda$--tree creatures $t$ such that
\begin{itemize}
\item $\eta[t]=\eta$,
\item $\dis[t]=\langle t_\xi:\xi\leq\delta\rangle$, where $t_0=t^*_\eta$ is
the minimal creature at $\eta$ for $(K,\Sigma)$, $\delta<\lambda_{\lh(
\eta)}$, and\quad $\xi<\zeta\leq\delta\quad \Rightarrow\quad t_\xi
\preceq^\eta_\Sigma t_\zeta\ \&\ t_\xi\neq t_\zeta$, 
\item $\pos[t]=\pos[t_\delta]$,
\item $\nor[t]=\min\{\nor[t_\xi]:\xi\leq\delta\}$.
\end{itemize}
Then, for $t',t\in K^\exl\cap\TCR_\eta[\bH]$ we let $t'\in\Sigma^\exl(t)$
if and only if $\dis[t]\trianglelefteq\dis[t']$.
\end{definition}

\begin{proposition}
Assume$(K,\Sigma)$ is a very local $\bar{\lambda}$--complete $\lambda$--tree
creating pair. Then $(K^\exl,\Sigma^\exl)$ is a very local exactly
$\bar{\lambda}$--complete $\lambda$--tree creating pair. The minimal
creature for it at $\eta$ is $t^{**}_\eta$ such that $\dis[t^{**}_\eta]=
\langle t^*_\eta\rangle$.
\end{proposition}

\begin{proof}
Easy.
\end{proof}

\begin{theorem}
\label{getSacks}
Suppose that
\begin{enumerate}
\item[(a)] $(\forall\alpha<\lambda)(|\bH(\alpha)|<\lambda_\alpha)$, and 
\item[(b)] $(K,\Sigma)$ is a $\bar{\mu}$--complete very local
$\lambda$--tree creating pair for $\bH$ for some strictly increasing
sequence $\bar{\mu}=\langle\mu_\alpha:\alpha<\lambda\rangle$ of regular
cardinals such that $\mu_\alpha<\lambda$ (for $\alpha<\lambda$), and   
\item[(c)] $D^*$ is a normal filter on $\lambda$. 
\end{enumerate}
Then the forcing notion $\bbQ^\tree_{D^*}(K,\Sigma)$ has the strong
$\bar{\lambda}$--Sacks property
\end{theorem}

\begin{proof}
Let $i_0<\lambda$  and $p\in\bbQ^\tree_{D^*}(K,\Sigma)$. Just for notational
simplicity we assume that $\bH(\alpha)\in\lambda_\alpha$ for all $\alpha<
\lambda$ and $\lh(\mrot(p))\leq i_0$. We are going to describe a strategy
for the Generic player in the game $\Gsg(i_0,p,\bbQ^\tree_{D^*}(K,
\Sigma))$. In the course of the play she will also choose sets $Y_{i+1}\in
D^*$ and $\lambda$--tree creatures $t_\nu$.    

First Generic picks $\eta\in T^p$ such that $\lh(\eta)>i_0$ and she starts
the game with playing $s_{i_0}=\{\eta\rest (i_0+1)\}$ and $q^{i_0}_{\eta
\rest (i_0+1)}=p^{[\eta]}$. She also picks $t_{\eta\rest i_0}\in
\Sigma(t^p_{\eta\rest (i_0+1)})$ such that $\pos[t_{\eta\rest (i_0+1)}]=\{
\eta\rest (i_0+1)\}$ (remember \ref{treecomp}(1)($\gamma$)). 

Arriving at a successor stage $j=i+1$ of the play the players have
determined $s_i,\bar{q}^i,\bar{p}^i$ and $Y_i$ so that, in addition to the
demands of the game, for each $\eta\in s_i\cap {}^{i+1}\lambda$ we have  
$\eta\trianglelefteq\mrot(q^i_\eta)$. Now for each $\eta\in s_i\cap
{}^{i+1}\lambda$ Generic picks $\nu_\eta\in T^{p^i_\eta}$ strictly extending
$\eta$ and she plays 
\[s_{i+1}=s_i\cup\{\nu_\eta\rest (i+2):\eta\in s_i\},\qquad
q^{i+1}_{\nu_\eta\rest (i+2)}= (p^i_\eta)^{[\nu_\eta]}\quad\mbox{ for }\eta\in
s_i\cap {}^{i+1}\lambda.\]
She also fixes a set $Y_{i+1}\in D^*$ of limit ordinals included in
$\bigcap\limits_{\eta\in s_i} Y^{q^{i+1}_{\nu_\eta\rest (i+2)}}$ (recall
\ref{treeforcing}(2)) and for $\eta\in s_i\cap {}^{i+1}\lambda$ she lets 
$t_\eta\in\Sigma(t^p_\eta)$ be such that $\pos[t_\eta]=\{\nu_\eta\}$.  

Now suppose that the players have arrived to a limit stage $\delta<\lambda$
of the game, and assume that $\delta\notin\bigcap\limits_{i<\delta}
Y_{i+1}$. Generic lets $s_\delta^*$ consist all sequences $\eta$ of length
$\delta$ such that $\eta\rest (i+1)\in s_i$ whenever $i_0\leq i<\delta$. For
each $\eta\in s_\delta^*$ she first picks a condition $r_\eta$ stronger than
all $p^i_{\eta\rest (i+1)}$ for $i_0\leq i<\delta$ (there is one by arguments
as in the proof of \ref{getcomp}(1)) and then she chooses $\nu_\eta\in
T^{r_\eta}$ strictly extending $\eta$. Then she plays
\[s_\delta=s_\delta^*\cup\{\nu_\eta\rest (\delta+1):\eta\in s_\delta^*\},
\qquad q^\delta_{\nu_\eta\rest (\delta+1)}= (r_\eta)^{[\nu_\eta]}\quad\mbox{
  for }\eta\in s_\delta^*.\]  
The $\lambda$--tree creatures $t_\eta$ (for $\eta\in s_\delta^*$) are chosen
as above: $t_\eta\in\Sigma(t^p_\eta)$, $\pos[t_\eta]=\{\nu_\eta\}$.

Finally suppose that we are at a limit stage $\delta<\lambda$ of the game
and $\delta\in\bigcap\limits_{i<\delta} Y_i$. Let $s_\delta^*$ be defined as
above and for each $\eta\in s_\delta^*$ let $r_\eta$ be a condition stronger
than all $p^i_{\eta\rest (i+1)}$ for $i_0\leq i<\delta$ and such that
$\mrot(r_\eta)=\eta$ and $\nor[t^{r_\eta}_\eta]\geq |\delta|$ (there is one
by arguments as in the proof of \ref{getcomp}(1) and the choice of the
$Y_i$'s). Then she plays
\[s_\delta=s_\delta^*\cup \bigcup\{\pos[t^{r_\eta}_\eta]:\eta\in s_\delta^*\},
\qquad q^\delta_\nu= (r_\eta)^{[\nu]}\quad\mbox{ for }\nu\in
\pos[t^{r_\eta}_\eta],\ \eta\in s_\delta^*.\]  
She also lets $t_\eta=t^{r_\eta}_\eta\in\Sigma(t^p_\eta)$ (for $\eta\in
s^*_\delta$). 

It should be clear that the strategy described above always tells Generic to
play legal moves (remember \ref{incon}(c)). It should also be clear that if
$\langle (s_i,\bar{q}^i,\bar{p}^i):i_0\leq i<\lambda\rangle$ is the result
of a play of $\Gsg(i_0,p,\bbQ^\tree_{D^*}(K,\Sigma))$ in which Generic uses
that strategy, then letting $T=\bigcup\{s_i:i_0\leq i<\lambda\}$ and $q=\langle
t_\eta:\eta\in T\rangle$ we get a condition in $\bbQ^\tree_{D^*}(K,\Sigma)$
(as witnessed by $\mathop{\triangle}\limits_{i<\lambda} Y_{i+1}$) stronger
than $p$ and forcing that 
\[\mbox{`` }(\exists\rho\in\bairel)(\forall i\in [i_0,\lambda))(\rho\rest (i+1)
\in s_i\ \&\ q^i_{\rho\rest
  (i+1)}\in\Gamma_{\bbQ^\tree_{D^*}(K,\Sigma)}\big)
\mbox{ ''.}\]
\end{proof}

\begin{theorem}
\label{getbound}
Assume that $(\forall\alpha<\lambda)(|\bH(\alpha)|<\lambda)$, and
$(K,\Sigma), D^*$ satisfy (b), (c) of \ref{getSacks}. Then the forcing
notion $\bbQ^\tree_{D^*}(K,\Sigma)$ has the strong $\lambda$--bounding
property.  
\end{theorem}

\begin{proof}
Similar to \ref{getSacks}.
\end{proof}

The above two theorems are applicable to forcing notions of the type
$\bbQ^\tree_{\cl}(K,\Sigma)$ as they may be treated as a special case (under
the assumptions as there):

\begin{proposition}
\label{same}
Assume that $(\forall\alpha<\lambda)(|\bH(\alpha)|<\lambda)$ and
$(K,\Sigma)$ is a $\bar{\mu}$--complete very local $\lambda$--tree creating
pair for $\bH$ (for some $\bar{\mu}$). Then  the forcing notions 
$\bbQ^\tree_{\cl} (K,\Sigma)$ and $\bbQ^\tree_{\cD_\lambda}(K,\Sigma)$ are
equivalent.  
\end{proposition}

Turning to the case of $\lambda$--creating pairs (and forcing notions of the
form $\bbQ^*_e(K,\Sigma)$), we do easy ways to ensure they are suitably
complete (parallel to \ref{treecomp}, \ref{getcomp}). However at present we
do not know how to get the strong bounding properties of section A.2 (which
were obviously tailored for trees).

\begin{definition}
\begin{enumerate}
\item For a $\lambda$--creating pair $(K,\Sigma)$ and $t\in K$ we define
when $t$ is {\em $\kappa$--complete} and {\em exactly $\kappa$--complete} 
like in \ref{treecomp}(1,2) (but with $\val$ replacing $\pos$). 
\item If $(K,\Sigma)$ is very local, then we say that it is
{\em $\bar{\lambda}$--complete\/} ({\em exactly $\bar{\lambda}$--complete},
resp.) if each minimal creature $t^*_\alpha$ is $\lambda^+_\alpha$--complete
(exactly $\lambda_\alpha$--complete, resp.). 
\end{enumerate}
\end{definition}

\begin{proposition}
Assume that $(K,\Sigma)$ is a very local $\lambda$--creating pair for
$\bH$, $D^*$ is a normal filter on $\lambda$. Let $\bbP$ be either the
forcing notion $\bbQ^*_1(K,\Sigma)$ or $\bbQ^*_{D^*}(K,\Sigma)$.
\begin{enumerate}
\item If $(K,\Sigma)$ is $\bar{\lambda}$--complete, then $\bbP$ is
strategically $({<}\lambda)$--complete.  
\item If $(K,\Sigma)$ is exactly $\bar{\lambda}$--complete, then 
$\bbP$ is $\lambda$--complete. 
\end{enumerate}
\end{proposition}

\section{Getting fuzzy properness}
In this section we show that the forcing notions with trees and creatures
may fit the fuzzy proper framework. Note that even though the forcing
notions covered by Theorems \ref{firstproper} and \ref{secondproper} below
are also covered by Theorem \ref{getbound}, the results here still have
value if we want to iterate that forcing notions with ones which do not have
the strong $\lambda$--bounding property. 

\begin{theorem}
\label{firstproper}
Let $A,W,D$ be as in \ref{Con3.1} and let $D^*$ be a normal filter on
$\lambda$ such that for some $S_0\in D^*$ we have $\lambda\setminus S_0\in
D$. Let $\bar{\mu}$ be an increasing sequence of regular cardinals cofinal
in $\lambda$. Assume that $(K,\Sigma)$ is an exactly $\bar{\mu}$--complete
very local $\lambda$--tree creating pair for $\bH$, and $|\bH(\alpha)|<
\lambda$ for each $\alpha<\lambda$. Then the forcing notion
$\bbQ^\tree_{D^*}(K,\Sigma)$ is strongly fuzzy proper over quasi
$D$--diamonds for $W$.    
\end{theorem}

\begin{proof}
By \ref{getcomp} we know that $\bbQ^\tree_{D^*}(K,\Sigma)$ is
$\lambda$--complete. 

Let $\gR^{\rm tr},\bar{\gY}^{\rm tr}$ be the trivial $\lambda$--base defined
as in the proof \ref{lambdaplus} (but for $\bbP=\bbQ^\tree_{D^*}(K,
\Sigma)$). We are going to show that for this $\lambda$--base and for $c=(
\bar{\lambda},K,\Sigma)$ the condition \ref{1.2}(2)($(\circledast)^+$)
holds. So assume that $N,h,\bar{F}=\langle F_\delta:\delta\in S\rangle$ and
$\bar{q}=\langle q_{\delta,x}:\delta\in S\mbox{ limit }\&\ x\in\cX_\delta
\rangle$ are as there and $p\in\bbQ^\tree_{D^*}(K,\Sigma)\cap N$. Note that
$\cX_\delta=\{0\}$ (for all relevant $\delta$) and thus we may think that
$\bar{q}=\langle q_\delta:\delta\in S\mbox{ limit}\rangle$.  

Let $\bcI=\langle \cI_\xi:\xi<\lambda\rangle$ list all open dense subsets of 
$\bbQ^\tree_{D^*}(K,\Sigma)$ from $N$. For $i<\lambda$ let $\xi_i$ be such
that $\cI_{\xi_i}$ consist of conditions $p\in\bbQ^\tree_{D^*}(K,\Sigma)$
with $\lh(\mrot(p))>i$, and let $E$ be a club of $\lambda$ such that 
\[(\forall \delta\in E)(\forall i<\delta)(\delta\mbox{ is limit and }\xi_i<
\delta).\]
By induction on $\alpha<\lambda$ choose conditions $p_\alpha\in 
\bbQ^\tree_{D^*}(K,\Sigma)\cap N$ and sets $Y_\alpha\in D^*$ such that
\begin{enumerate}
\item[(i)] $p_0=p$, $\mrot(p_\alpha)=\mrot(p)$, and $p_\alpha\leq p_\beta$ and 
$Y_\beta\subseteq Y_\alpha\subseteq S_0$ for $\alpha<\beta<\lambda$, 
\item[(ii)] $Y_\alpha$ witnesses $p_\alpha\in\bbQ^\tree_{D^*}(K,\Sigma)$
  (see \ref{treeforcing}(2)),
\item[(iii)] for every $\alpha<\beta<\lambda$ and $\nu\in (T^{p_\alpha}
)_\alpha$ we have $\nu\in T^{p_\beta}$ and $t^{p_\alpha}_\nu=
t^{p_\beta}_\nu$, 
\item[(iv)] if $\alpha<\lambda$ is a successor, $\xi<\alpha$ and $\eta\in
(T^{p_\alpha})_\alpha$, then for some $\nu\in (T^{p_\alpha})^{[\eta]}$ we
have: $(p_\alpha)^{[\nu]}\in \cI_\xi$ and $(\forall\rho\in T^{p_\alpha})
(\eta\trianglelefteq\rho\vtl\nu\ \Rightarrow\ \nor[
t_\rho^{p_\alpha}]=0)$,
\item[(v)] if $\delta\in\bigcap\limits_{\alpha<\delta} Y_\alpha$ is a limit
ordinal, then $\delta\in Y_\beta$ for every $\beta\geq\delta$, 
\item[(vi)] if $\delta\in S\cap E\setminus S_0$ and $\langle h\circ
F_\delta(i):i<\delta\rangle$ is an increasing $\bcI$--exact sequence of
members of $N\cap \bbQ^\tree_{D^*}(K,\Sigma)$ such that
\[(\forall\alpha<\delta)(\exists i<\delta)(p_\alpha\leq h\circ
F_\delta(i)),\]
and $\eta\in (T^{p_\delta})_\delta$ is such that every $h\circ F_\delta(i)$
is compatible with $(p_\delta)^{[\eta]}$, then $(p_\delta)^{[\eta]}\leq
q_\delta=(p_\delta)^{[\mrot(q_\delta)]}$ and 
\[(\forall\rho\in T^{p_\delta})(\eta\trianglelefteq\rho\vtl\mrot(q_\delta)\
\Rightarrow\ \nor[t_\rho^{p_\delta}]=0).\] 
[Note that there is at most one $\eta$ as above; remember the choice of $E$.]
\end{enumerate}
It should be clear that the inductive construction of the $p_\alpha$'s and  
$Y_\alpha$'s is possible (for (v) remember $\delta<\lambda_\delta$; note
also that there is no collision between (v) and (vi) because
$Y_\alpha\subseteq S_0$). Now letting $\mrot(r)=\mrot(p)$, $T^r=
\bigcup\limits_{\alpha<\lambda}(T^{p_\alpha})_\alpha$, $t^r_\nu=
t^{p_{\lh(\nu)+1}}_\nu$ we get a condition $r\in\bbQ^\tree_{D^*}(K,\Sigma)$
(as witnessed by $\mathop{\triangle}\limits_{\alpha<\lambda}
Y_\alpha$). Also note that $r$ is stronger than all $p_\alpha$'s. 

\begin{claim}
\label{cl4}
The condition $r$ is $(\gR^{\rm tr},\bar{\gY}^{\rm tr})$--fuzzy generic for
$\bar{q}$.
\end{claim}

\begin{proof}[Proof of the Claim]
First note that the condition $r$ is $(N,\bbQ^\tree_{D^*}(K,
\Sigma))$--generic by clause (iv) above. Therefore we may use
\ref{1.2A}(3), and it is enough that we show that the Generic player has
a strategy in the game $\Gfl(r,N,\bcI,h,\bbQ^\tree_{D^*}(K,\Sigma),\bar{F}, 
\bar{q})$ which guarantees that the result $\langle r_i,C_i:i<\lambda
\rangle$ of the play satisfies \ref{pre1.2}(5)($\beta$). Let us describe
such a strategy.

First, for $\alpha<\lambda$ let $\zeta_\alpha<\lambda$ be such that 
\[(\forall q\in\cI_{\zeta_\alpha})(\mbox{ either }p_\alpha\leq q\mbox{ or
$p_\alpha,q$ are incompatible }),\]
and let $E'=\{\delta\in E:(\forall\alpha<\delta)(\zeta_\alpha<\delta)\}$ (it is a
club of $\lambda$).

Now, suppose that during a play of $\Gfl(r,N,\bcI,h,\bbQ^\tree_{D^*}(K,
\Sigma),\bar{F},\bar{q})$ the players have arrived at stage $i\in S$ having 
constructed a sequence $\langle r_j,C_j:j<i\rangle$.

If either $i$ is a successor ordinal or $i\notin\bigcap\limits_{j<i} C_j$,
then the Generic player plays the $<^*_\chi$--first condition $r_i\in
\bbQ^\tree_{D^*}(K,\Sigma)$ such that $(\forall j<i)(r_j \leq r_i)$ and
$\lh(\mrot(r_i))>i$, and the set $C_i=E'\setminus (S_0\cup
\lh(\mrot(r_i)))$. 

If $i\in\bigcap\limits_{j<i} C_j$ is a limit ordinal (so also $i\in
E'\setminus S_0$), then Generic asks
\begin{enumerate}
\item[$(*)$] is $\langle h\circ F_i(\alpha):\alpha<i\rangle$ an increasing 
$\bcI$--exact sequence such that 
\[(\forall j<i)(\exists\alpha<i)(p_j\leq h\circ F_i(\alpha))\mbox{ ?}\]
\end{enumerate}
If the answer to $(*)$ is ``no'', then she plays like at the successor
stage. 

\noindent [Note that if the answer to $(*)$ is ``no'' and $\langle h\circ
F_i(\alpha):\alpha<i\rangle$ is increasing $\bcI$--exact, then for some
$j<i$ and $\alpha<i$ the conditions $p_j$ and $h\circ F_i(\alpha)$ are
incompatible, and hence $r_i$ and $h\circ F_i(\alpha)$ are incompatible.]\\ 
If the answer to $(*)$ is ``yes'', then Generic looks at clause (vi) (of the
choice of $p_\alpha$'s) and $\eta=\bigcup\limits_{j<i}\mrot(r_j)$ (note that
$\lh(\eta)=i$). If $(p_i)^{[\eta]}$ is incompatible with some $h\circ
F_i(\alpha)$, $\alpha<i$, then she plays $C_i,r_i$ as in the successor
case. 

\noindent [Note that then $r_i,h\circ F_i(\alpha)$ are incompatible.]\\
Otherwise $\eta\trianglelefteq\mrot(q_i)\in T^{p_i}$, $q_i=(p_i)^{[\mrot(
q_i)]}$ and $(\forall\rho\in T^{p_i})(\eta\trianglelefteq\rho
\vtl\mrot(q_i)\ \Rightarrow\ |\pos[t_\rho^{p_i}]|=1)$. 
Therefore, $\mrot(q_i)\in T^{r_j}$ and $q_i\leq (r_j)^{[\mrot(q_i)]}$
(for each $j<i$). So the Generic player can play $C_i=E'\setminus i$ and the
$<^*_\chi$--first condition $r_i$ stronger than all $r_j$ (for $j<i$) and
$q_i$. 
\medskip

It follows immediately from the comments stated during the description of
the strategy that every play according to it satisfies
\ref{pre1.2}(5)($\beta$), finishing the proof of the claim and that of the  
theorem. 
\end{proof}
\end{proof}

\begin{theorem}
\label{secondproper}
Let $A,W,D$ be as in \ref{Con3.1} and let $D^*$ be a normal filter on
$\lambda$ such that for some $S_0\in D^*$ we have $\lambda\setminus S_0\in
D$. Assume that $(K,\Sigma)$ is an exactly $\bar{\lambda}$--complete very
local $\lambda$--creating pair for $\bH$, $|\bH(\alpha)|<\lambda$ for each 
$\alpha<\lambda$. Then the forcing notion $\bbQ^*_{D^*}(K,\Sigma)$ is
strongly fuzzy proper over quasi $D$--diamonds for $W$.   
\end{theorem}

\begin{proof}
Like \ref{firstproper}
\end{proof}

\begin{theorem}
\label{thirdproper}
Let $A=\cH_{<\lambda}(\lambda^*)$, $\lambda^*>\lambda$ and $W\subseteq
[A]^{\textstyle \lambda}$ be as in \ref{Con3.1}. Let $\bar{\mu}$ be an
increasing sequence of regular cardinals cofinal in $\lambda$.
Suppose that $(K,\Sigma)$ is an exactly $\bar{\mu}$--complete very local
$\lambda$--tree creating pair for $\bH$, $(\forall\alpha<\lambda)(|\bH(
\alpha)|<\lambda)$, and $D^*$ is a normal filter on $\lambda$. Then the
forcing notion $\bbQ^\tree_{D^*}(K,\Sigma)$ is fuzzy proper for $W$.
\end{theorem}

\begin{proof}
The proof closely follows the lines of that of \ref{firstproper}. Let $D$ be 
a normal filter on $\lambda$ such that there is a $D$--diamond.  

Just only to simplify somewhat the definition of a $\lambda$--base which we
will use, let us assume that $\bigcup\limits_{\delta<\lambda}\prod\limits_{
\alpha<\delta}\bH(\alpha)\subseteq a$ for every $a\in W$. Now we let $\gR=
\gR(K,\Sigma)$ consist of all triples $(p,\delta,\eta)$ such that $\delta<
\lambda$, $\eta\in\prod\limits_{\alpha<\delta}\bH(\alpha)$ and $p\in
\bbQ^\tree_{D^*}(K,\Sigma)$ satisfies $\eta\vtl\mrot(p)$. Next, for $a\in W$
let $\gY_a=\gY_a(K,\Sigma):\lambda\longrightarrow [a]^{\textstyle<\lambda}$
be given by $\gY_a(\delta)=\prod\limits_{\alpha\leq\delta}\bH(\alpha)
\subseteq a$ (for $\delta<\lambda$). It should be clear that $(\gR,
\bar{\gY})$ is a $\lambda$--base for $\bbQ^\tree_{D^*}(K,\Sigma)$ over $W$.    

We claim that $(\gR,\bar{\gY})$ and $c=(\bar{\lambda},\bH,K,\Sigma)$ witness 
the condition $(\circledast)$ of \ref{1.2}(1). To this end, let $N,h,\bar{F}=
\langle F_\delta:\delta\in S\rangle$ and $\bar{q}=\langle q_{\delta,x}:
\delta\in S\mbox{ limit }\&\ x\in\cX_\delta\rangle$ be as in \ref{1.2}(1)$(
\circledast)$, $p\in\bbQ^\tree_{D^*}(K,\Sigma)\cap N$. Let $\bcI=\langle
\cI_\xi:\xi<\lambda\rangle$ list all open dense subsets of $\bbQ^\tree_{D^*}
(K,\Sigma)$ from $N$. For $i<\lambda$ let $\xi_i$ be such that $\cI_{\xi_i}$
consist of conditions $p\in\bbQ^\tree_{D^*}(K,\Sigma)$ with
$\lh(\mrot(p))>i$, and let $E$ be a club of $\lambda$ such that  
\[(\forall \delta\in E)(\forall i<\delta)(\delta\mbox{ is limit and }\xi_i<
\delta).\]
By induction on $\alpha<\lambda$, like in \ref{firstproper} (but note the
change in (vi) below!), we choose conditions $p_\alpha\in\bbQ^\tree_{D^*}(K,
\Sigma)\cap N$ and sets $Y_\alpha\in D^*$ such that 
\begin{enumerate}
\item[(i)] $p_0=p$, $\mrot(p_\alpha)=\mrot(p)$, and $p_\alpha\leq p_\beta$ and 
$Y_\beta\subseteq Y_\alpha$ for $\alpha<\beta<\lambda$, 
\item[(ii)] $Y_\alpha$ witnesses $p_\alpha\in\bbQ^\tree_{D^*}(K,\Sigma)$,
\item[(iii)] for every $\alpha<\beta<\lambda$ and $\nu\in (T^{p_\alpha}
)_\alpha$ we have $\nu\in T^{p_\beta}$ and $t^{p_\alpha}_\nu=
t^{p_\beta}_\nu$, 
\item[(iv)] if $\alpha<\lambda$ is a successor ordinal and $\eta\in (
T^{p_\alpha})_\alpha$, then for some $\nu\in (T^{p_\alpha})^{[\eta]}$ we
have: $(p_\alpha)^{[\nu]}\in \bigcap\limits_{\xi<\alpha}\cI_\xi$ and
$(\forall\rho\in T^{p_\alpha})(\eta\trianglelefteq\rho\vtl\nu\ \Rightarrow\
\nor[t_\rho^{p_\alpha}]=0)$,
\item[(v)] if $\delta\in\bigcap\limits_{\alpha<\delta} Y_\alpha$ is a limit
ordinal, then $\delta\in Y_\beta$ for every $\beta\geq\delta$, 
\item[(vi)] if $\delta\in S\cap E$, $\langle h\circ F_\delta(i):i<\delta
\rangle$ is increasing $\bcI$--exact, $\eta=\bigcup\limits_{i<\delta} 
\mrot(h\circ F_\delta(i))$ and $\lh(\eta)=\delta$, and $(\forall\alpha<
\delta)(\exists i<\delta)(p_\alpha\leq h\circ F_\delta(i))$,\\ 
then ($\eta\in T^{p_\delta}$ and) for every $\nu\in\pos[t^{p_\delta}_\eta]
\cap \bigcap\limits_{i<\delta}\pos[t^{h\circ F_\delta(i)}_\eta]$ we have  
\[(p_\delta)^{[\nu]}\leq q_{\delta,\nu}=(p_\delta)^{[\mrot(q_{\delta,
\nu})]}\ \mbox{ and }\ (\forall\rho\in T^{p_\delta})(\nu\trianglelefteq
\rho\vtl\mrot(q_{\delta,\nu})\ \Rightarrow\ \nor[t_\rho^{p_\delta}]=0).\] 
[Note that, in the situation as in (vi), $\cX_\delta=\bigcap\limits_{i<
\delta}\pos[t^{h\circ F_\delta(i)}_\eta]$.] 
\end{enumerate}
Plainly, the inductive construction of the $p_\alpha$'s and $Y_\alpha$'s is 
possible (for (v) remember $\delta<\lambda_\delta$). Now letting $\mrot(r)=
\mrot(p)$, $T^r=\bigcup\limits_{\alpha<\lambda}(T^{p_\alpha})_\alpha$,
$t^r_\nu=t^{p_{\lh(\nu)+1}}_\nu$ we get a condition $r\in\bbQ^\tree_{D^*}(K,
\Sigma)$ stronger than all $p_\alpha$'s. 

\begin{claim}
\label{cl8}
The condition $r$ is $(\gR,\bar{\gY})$--fuzzy generic for $\bar{q}$.
\end{claim}

\begin{proof}[Proof of the Claim]
It is very much like \ref{cl4}. We note that $r$ is $(N,\bbQ^\tree_{D^*}(K,
\Sigma))$--generic (by clause (iv)), and therefore it is enough to show that
the Generic player has a strategy in the game $\Gfl(r,N,\bcI,h,
\bbQ^\tree_{D^*}(K,\Sigma),\bar{F},\bar{q})$ which guarantees that the
result $\langle r_i,C_i:i<\lambda\rangle$ of the play satisfies
\ref{pre1.2}(5)($\beta$) (remember \ref{1.2A}(3)). Let us describe such a
strategy. First, for $\alpha<\lambda$ let $\zeta_\alpha<\lambda$ be such
that  
\[(\forall q\in\cI_{\zeta_\alpha})(\mbox{ either }p_\alpha\leq q\mbox{ or
$p_\alpha,q$ are incompatible }),\]
and let $E'=\{\delta\in E:(\forall\alpha<\delta)(\zeta_\alpha<\delta)\}$ (it is a
club of $\lambda$).

Now, suppose that during a play of $\Gfl(r,N,\bcI,h,\bbQ^\tree_{D^*}(K,
\Sigma),\bar{F},\bar{q})$ the players have arrived at stage $i\in S$ having 
constructed a sequence $\langle r_j,C_j:j<i\rangle$.

If either $i$ is a successor ordinal or $i\notin\bigcap\limits_{j<i} C_j$,
then Generic plays the $<^*_\chi$--first condition $r_i\in
\bbQ^\tree_{D^*}(K,\Sigma)$ such that $(\forall j<i)(r_j\leq r_i)$ and
$\lh(\mrot(r_i))>i$ and $C_i=E'\setminus \lh(\mrot(r_i))$.

If $i\in\bigcap\limits_{j<i} C_j\subseteq E'$ is a limit ordinal, then
Generic asks 
\begin{enumerate}
\item[$(*)$] is $\langle h\circ F_i(\alpha):\alpha<i\rangle$ an increasing 
$\bcI$--exact sequence such that 
\[(\forall\alpha<i)(\exists j<i)(h\circ F_i(\alpha)\leq r_j)\mbox{ ?}\]
\end{enumerate}
If the answer to $(*)$ is ``no'', then she plays like at the successor
stage.\\
If the answer to $(*)$ is ``yes'', then Generic takes $\eta=
\bigcup\limits_{j<i}\mrot(r_j)$ and she notes that $\lh(\eta)=\delta$ (by 
the choice of $C_j$'s at successor stages) and $\eta=\bigcup\limits_{i<
\delta}\mrot(h\circ F_\delta(i))$. Also, by the exactness and the choice of 
$E'$, we have 
\[(\forall j<i)(\exists\alpha<i)(p_j\leq h\circ F_i(\alpha)).\]
So now Generic looks at clause (vi) of the choice of $p_\alpha$'s.  She
picks (say, the $<^*_\chi$--first) $\nu\in\bigcap\limits_{j<i}\pos[
t^{r_j}_\eta]\subseteq\bigcap\limits_{j<i}\pos[t^{h\circ F_i(\alpha)}_\eta]$ 
and notices that  (by (vi)) $\nu\trianglelefteq\mrot(q_{i,\nu})\in T^{p_i}$,
$q_{i,\nu}= (p_i)^{[\mrot(q_{i,\nu})]}$ and $(\forall\rho\in T^{p_i})(\nu
\trianglelefteq\rho\vtl\mrot(q_{i,\nu})\ \Rightarrow\
|\pos[t_\rho^{p_i}]|=1)$. Therefore, $\mrot(q_{i,\nu})\in T^{r_j}$ and
$q_{i,\nu}\leq (T^{r_j})^{[\mrot(q_{i,\nu})]}$ (for each $j<i$). So the
Generic player can play $C_i=E'\setminus i$ and the $<^*_\chi$--first
condition $r_i$ stronger than all $r_j$ (for $j<i$) and $q_{i,\nu}$. 
\medskip

Easily, the strategy described above has the required property, and the
proof is completed. 
\end{proof}
\end{proof}

\begin{problem}
Unlike that was in the case of \ref{firstproper}, it is not clear how the
proof of \ref{thirdproper} can be modified to get the parallel result for
non-tree case. So, assuming that $A,W,\bH$ and $D^*$ are as in
\ref{thirdproper} and $(K,\Sigma)$ is an exactly $\bar{\lambda}$--complete
very local $\lambda$--creating pair for $\bH$, is the forcing notion
$\bbQ^*_{D^*}(K,\Sigma)$ fuzzy proper for $W$?
\end{problem}

\section{More examples and applications}
Here we are going to present some direct applications of the methods
developed in this paper. Though we do keep our basic assumptions from
\ref{incon}, we are going to introduce more parameters, so let us
fully state the context we are working in now. 

\begin{context}
\begin{enumerate}
\item[(a)] $\lambda$ is a strongly inaccessible cardinal, $2^\lambda=
\lambda^+$, and $2^{\lambda^+}=\lambda^{++}$, and 
\item[(b)] $\bar{\mu}=\langle\mu_\alpha:\alpha<\lambda\rangle$,
$\bar{\lambda}=\langle\lambda_\alpha:\alpha<\lambda\rangle$ and
$\bar{\kappa}=\langle\kappa_\alpha:\alpha<\lambda\rangle$ are strictly
increasing sequences of uncountable regular cardinals, each cofinal in
$\lambda$, 
\item[(c)] for each $\alpha<\lambda$,
\begin{itemize}
\item $\mu_\alpha^+<\lambda_\alpha<\kappa_\alpha$, 
\item $\prod\limits_{\beta<\alpha}\lambda_\beta <\lambda_\alpha$ and
  $(\forall\xi<\lambda_\alpha)(|\xi|^\alpha<\lambda_\alpha)$, 
\end{itemize}
\item[(d)] $A=\cH_{<\lambda}(\lambda^*)$, $\lambda^*>\lambda$ and
  $W\subseteq [A]^{\textstyle \lambda}$ are as in \ref{Con3.1}, 
\item[(e)] $D$ is a normal filter on $\lambda$ such that there is a
  $D$--diamond.  
\end{enumerate}
\end{context}

Let us recall some notions related to cardinal characteristics of
$\lambda$--reals. 

\begin{definition}
\begin{enumerate}
\item Let $\cS_{\bar{\mu}}$ be the family of all sequences
$\bar{a}=\langle a_\alpha;\alpha<\lambda\rangle$ such that $a_\alpha\in
[\lambda]^{\textstyle<\mu_\alpha}$ (for all $\alpha<\lambda$). We
define 
\[\begin{array}{rl}
c(\bar{\mu})=&\min\big\{|\cY|:\cY\subseteq\cS_{\bar{\mu}}\ \&\ (\forall
f\in\bairel)(\exists\bar{a}\in\cY)(\forall\alpha<\lambda)(f(\alpha)
\in a_\alpha)\big\},\\  
c_{\cl}^-(\bar{\mu})=&\min\big\{|\cY|:\cY\subseteq\cS_{\bar{\mu}}\ \&\
(\forall f\in\bairel)(\exists\bar{a}\in\cY)\big(\{\alpha<\lambda:f(
\alpha)\in a_\alpha\}\in (\cD_\lambda)^+\big)\big\},
  \end{array}\]
and also
\[\begin{array}{rl}
e_{\cl}(\bar{\mu})=\min\big\{|\cG|:&\cG\subseteq\prod\limits_{\alpha<
\lambda}\mu_\alpha\ \mbox{ and}\\
&(\forall f\in \prod\limits_{\alpha<\lambda}\mu_\alpha)(\exists g\in
\cG)\big(\{\alpha<\lambda:f(\alpha)\neq g(\alpha)\}\in\cD_\lambda\big)
\big\},  
  \end{array}\]
\item For an ideal $\cJ$ of subsets of a set $\cX$, the {\em covering
  number\/} $\cov(\cJ)$ of $\cJ$ is 
\[\cov(\cJ)=\min\{|\cY|:\cY\subseteq\cJ\ \&\ \bigcup\cY=\cX\}.\]
\end{enumerate}
\end{definition}

\begin{proposition}
It is consistent that $c(\bar{\lambda})<e_\cl(\bar{\mu})$. 
\end{proposition}

\begin{proof}
Let $\bH_0(\alpha)=\mu_\alpha$ (for $\alpha<\lambda$)  and let $K_0$
consist of all $\lambda$--tree creatures $t\in\TCR[\bH_0]$ such that: 
\begin{itemize}
\item $\dis[t]\in\mu_{\lh(\eta[t])}+1$, 
\item if $\dis[t]=\mu_{\lh(\eta[t])}$, then $\pos[t]=\{\eta[t]\conc
\langle\xi\rangle:\xi<\mu_{\lh(\eta[t])}\}$ and $\nor[t]=\mu_{\lh(\eta
[t])}$, 
\item if $\dis[t]<\mu_{\lh(\eta[t])}$, then $\pos[t]=\{\eta[t]\conc
  \langle\dis[t]\rangle\}$ and $\nor[t]=0$.
\end{itemize}
Let $\Sigma_0$ be a local tree-composition operation on $K_0$ (so its
domain consists of singletons only) such that
\begin{itemize}
\item if $\dis[t]<\mu_{\lh(\eta[t])}$, then $\Sigma_0(t)=\{t\}$, 
\item if $\dis[t]=\mu_{\lh(\eta[t])}$, then $\Sigma_0(t)=\{t'\in K_0:
  \eta[t']=\eta[t]\}$. 
\end{itemize}
It should be clear that $(K_0,\Sigma_0)$ is a very local exactly
$\bar{\lambda}$--complete tree creating pair. The forcing notion
$\bbQ^\tree_{\cD_\lambda}(K_0,\Sigma_0)$ has the strong
$\bar{\lambda}$--Sacks property (by \ref{getSacks}). Let $\name{W}$ be
the canonical $\bbQ^\tree_{\cD_\lambda}(K_0,\Sigma_0)$--name for the
generic function in $\prod\limits_{\alpha<\lambda}\mu_\alpha$, so 
\[p\forces_{\bbQ^\tree_{\cD_\lambda}(K_0,\Sigma_0)}\mbox{`` }\mrot(p)
\vtl\name{W}\mbox{ ''.}\]
Then we have 
\[\forces_{\bbQ^\tree_{\cD_\lambda}(K_0,\Sigma_0)}\mbox{`` }
\big(\forall f\in\prod\limits_{\alpha<\lambda}\mu_\alpha\cap\bV\big)
\big(\{\alpha<\lambda:\name{W}(\alpha)=f(\alpha)\}\in(\cD_\lambda)^+
\big)\mbox{ ''.}\] 
Now let $\bbP$ be the limit of a $\lambda$--support iteration,
$\lambda^{++}$ in length, of the forcing notions
$\bbQ^\tree_{\cD_\lambda}(K_0,\Sigma_0)$. Then, by \ref{itSacks} +
\ref{lppcc} + \ref{first} + \ref{1.6}, 
\begin{itemize}
\item $\bbP$ is $\lambda$--complete, $\lambda$--proper and satisfies
  the $\lambda^{++}$--cc, and it has a dense subset of size
  $\lambda^{++}$, thus forcing with $\bbP$ does not collapse
  cardinals, 
\item $\bbP$ has the $\bar{\lambda}$--Sacks property, it is weakly
  fuzzy proper for $W$,  
\item $\forces_{\bbP}$`` $2^\lambda=2^{\lambda^+}=\lambda^{++}=e_\cl(
  \bar{\mu})\ \mbox{ and }\ c(\bar{\lambda})=\lambda^+ $ ''
\end{itemize}
\end{proof}

\begin{remark}
The forcing $\bbQ^\tree_{\cD_\lambda}(K_0,\Sigma_0)$ is a ``bounded
relative'' of $\bbD_\lambda$ from \cite[4.10]{RoSh:655} (remember
\ref{same}). It is also a generalization of the forcing notions
$\bbD_{\cX}$ from \cite{NeRo93}. 
\end{remark}

\begin{proposition}
\label{different}
It is consistent that  $c(\bar{\lambda})<c^-_\cl(\bar{\mu}^+)=c(
\bar{\mu}^+)$, where $\bar{\mu}^+=\langle\mu_\alpha^+:\alpha<\lambda
\rangle$.    
\end{proposition}

\begin{proof}
Let $\bH_1(\alpha)=\mu_\alpha^+$ (for $\alpha<\lambda$)  and let $K_1'$
consist of all $\lambda$--tree creatures $t\in\TCR[\bH_1]$ such that: 
\begin{itemize}
\item $\dis[t]\subseteq\mu_{\lh(\eta[t])}^+$, either  $|\dis[t]|=1$ or
  $\dis[t]$ is a club of $\mu_{\lh(\eta[t])}^+$, 
\item $\pos[t]=\{\eta[t]\conc\langle\xi\rangle:\xi\in \dis[t]\}$, 
 \item if $|\dis[t]|=1$ then $\nor[t]=0$, if $|\dis[t]|>1$ then 
$\nor[t]=\mu_{\lh(\eta[t])}$.  
\end{itemize}
Let $\Sigma_1'$ be a local tree-composition operation on $K_1'$ such
that 
\[\Sigma_1'(t)=\{t'\in K_1':\eta[t']=\eta[t]\ \&\ \dis[t']\subseteq
\dis[t]\}.\]
Then $(K_1',\Sigma_1')$ is a very local $\bar{\mu}$--complete
$\lambda$--tree creating pair. Let $(K_1,\Sigma_1)$ be the
$\bar{\mu}$--exactivity of $(K_1',\Sigma_1')$ (see \ref{exactivity}); thus
$(K_1,\Sigma_1)$ is a very local exactly $\bar{\mu}$--complete
$\lambda$--tree creating pair. The forcing notion $\bbQ^\tree_{\cD_\lambda}
(K_1,\Sigma_1)$ is $\lambda$--complete fuzzy proper for $W$ and it has the
strong $\bar{\lambda}$--Sacks property. Also, letting $\name{W}$ be the
canonical name for the generic function in $\prod\limits_{\alpha<\lambda}
\mu_\alpha^+$ (i.e., $p\forces_{\bbQ^\tree_{\cD_\lambda}(K_1,\Sigma_1)}
\mbox{`` }\mrot(p)\vtl\name{W}$ ''), we have
\[\forces_{\bbQ^\tree_{\cD_\lambda}(K_1,\Sigma_1)}\mbox{`` }
\big(\forall \bar{a}\in \cS_{\bar{\mu}^+}\cap\bV\big)\big(\{\alpha<\lambda:
\name{W}(\alpha)\notin a_\alpha)\}\in\cD_\lambda\big)\mbox{ ''.}\] 
Let $\bbP$ be the limit of a $\lambda$--support iteration, $\lambda^{++}$ in
length, of the forcing notions $\bbQ^\tree_{\cD_\lambda}(K_1,
\Sigma_1)$. Then (by \ref{itSacks} + \ref{lppcc} + \ref{first} + \ref{1.6})
we have: 
\begin{itemize}
\item $\bbP$ is $\lambda$--complete, $\lambda$--proper and satisfies
  the $\lambda^{++}$--cc, and it has a dense subset of size
  $\lambda^{++}$, thus forcing with $\bbP$ does not collapse
  cardinals, 
\item $\bbP$ has the $\bar{\lambda}$--Sacks property, it is weakly
  fuzzy proper for $W$,  
\item $\forces_{\bbP}$`` $2^\lambda=2^{\lambda^+}=\lambda^{++}=c^-_\cl(
  \bar{\mu}^+)\mbox{ and }\  c(\bar{\lambda})=\lambda^+ $ ''
\end{itemize}
\end{proof}

\begin{remark}
The result in \ref{different} is of interest as it shows that the
$\lambda$--versions of cardinal characteristics of the reals may behave
totally differently from their ``ancestors''. Recall that if for an
increasing function $f\in {}^\omega\omega$ we let $\cS^f$ consist of all
sequences $\bar{a}=\langle a_n:n<\omega\rangle$ with $a_n\in
[\omega]^{\textstyle \leq f(n)+1}$ (for $n<\omega$), then 
\[\begin{array}{l}
\min\big\{|\cY|:\cY\subseteq\cS^f\ \&\ (\forall
h\in {}^\omega\omega)(\exists\bar{a}\in\cY)(\forall n<\omega)(h(n)
\in a_n)\big\}=\\
\min\big\{|\cY|:\cY\subseteq\cS^g\ \&\ (\forall
h\in {}^\omega\omega)(\exists\bar{a}\in\cY)(\forall n<\omega)(h(n)
\in a_n)\big\}
  \end{array}\]
for any increasing $f,g\in {}^\omega\omega$

The $\lambda$--tree creating pair $(K_1,\Sigma_1)$ may be treated (in some
sense) as a special case of the $\lambda$--tree creating pairs $(K(
\bar{\cA}),\Sigma(\bar{\cA}))$ from \ref{treeforA} below. 
\end{remark}

\begin{definition}
\label{pre}
Let $\cA$ be a family of subsets of $\kappa$ such that $\kappa\in\cA$.
\begin{enumerate}
\item A game $\Game^*(\cA,\mu)$ of two players, I and II, is defined as
follows. A play lasts $\mu$ moves, in the $\alpha^{\rm th}$ move a set
$A_\alpha\in\cA$ is chosen, and player I chooses $A_\alpha$ for even
$\alpha$'s. In the end player II wins if $\bigcap\limits_{\alpha<\mu}
A_\alpha\neq\emptyset$. 
\item The family $\cA$ is a {\em $\mu$--category prebase on $\kappa$} if
player II has a winning strategy in the game $\Game^*(\cA,\mu)$ and
$(\forall A\in\cA)(\forall\xi<\kappa)(\exists B\in\cA)(B\subseteq
A\setminus\{\xi\})$. 
\item A set $X\subseteq\kappa$ is {\em $\cA$--presmall} if 
\[(\forall A\in\cA)(\exists B\in\cA)(B\subseteq A\setminus X).\]
\end{enumerate}
\end{definition}

Of course, every $\mu^+$--complete uniform filter $D^*$ on $\kappa$ is a
$\mu$--category base on $\kappa$ and then a set is $D$--presmall if and only 
if its complement is in $D^*$.   

\begin{definition}
\label{small}
\begin{enumerate}
\item {\em A $\bar{\lambda}$--smallness base on $\bar{\kappa}$} is a
sequence $\bar{\cA}=\langle\cA_\alpha:\alpha<\lambda\rangle$ such that each 
$\cA_\alpha$ is a $\lambda_\alpha$--category prebase on $\kappa_\alpha$.

Let $\bar{\cA}$ be a $\bar{\lambda}$--smallness base on $\bar{\kappa}$.
\item Let $T\subseteq\bigcup\limits_{\alpha<\lambda}\prod\limits_{\beta<
\lambda}\kappa_\alpha$ be a complete $\lambda$--tree with $\max(T)=
\emptyset$ and $D^*$ be a filter on $\lambda$. We say that  
\begin{itemize}
\item {\em $T$ is $\bar{\cA}$--small\/} if for every $\eta\in (T)_\alpha$,
$\alpha<\lambda$, the set $\{\xi<\kappa_\alpha:\eta\conc\langle\xi\rangle\in 
T\}$ is $\cA_\alpha$--presmall;
\item {\em $T$ is $(D^*,\bar{\cA})$--small\/} if 
\[\Big\{\alpha<\lambda:\mbox{ for every $\eta\in (T)_\alpha$ the set
}\{\xi<\kappa_\alpha:\eta\conc\langle\xi\rangle\in T\} 
\mbox{ is $\cA_\alpha$--presmall }\Big\}\in D^*.\]
\end{itemize}
\item Let $\cJ_{\bar{\kappa}}(\bar{\cA})$ consist of all subsets $X$ of
$\prod\limits_{\alpha<\lambda}\kappa_\alpha$ such that $X\subseteq
\bigcup\limits_{\vare<\lambda}\llim(T_\vare)$ for some
$\bar{\cA}$--small trees $T_\vare\subseteq\bigcup\limits_{\alpha<
\lambda}\prod\limits_{\beta<\lambda}\kappa_\alpha$ (for $\vare<
\lambda$).  

$\cJ_{\bar{\kappa}}(D^*,\bar{\cA})$ is defined similarly, replacing
``$\bar{\cA}$--small'' by ``$(D^*,\bar{\cA})$--small''.
\end{enumerate}
\end{definition}

\begin{proposition}
Let $\bar{\cA}$ be a $\bar{\lambda}$--smallness base on $\bar{\kappa}$. Then
both $\cJ_{\bar{\kappa}}(\bar{\cA})$ and $\cJ_{\bar{\kappa}}(D^*,\bar{\cA})$
are proper $\lambda^+$--complete ideals of subsets of
$\prod\limits_{\alpha<\lambda}\kappa_\alpha$, $\cJ_{\bar{\kappa}}(\bar{\cA})
\subseteq\cJ_{\bar{\kappa}}(D^*,\bar{\cA})$. They contain singletons and
$\lambda<\cov\big(\cJ_{\bar{\kappa}}(D^*,\bar{\cA})\big)\leq\cov\big(\cJ_{
\bar{\kappa}}(\bar{\cA})\big)$.
\end{proposition}

\begin{proposition}
\label{treeforA}
Let $\bar{\cA}$ be a $\bar{\lambda}$--smallness base on $\bar{\kappa}$ and
$D^*$ be a normal filter on $\lambda$. It is consistent that
$\cov\big(\cJ_{\bar{\kappa}}(D^*,\bar{\cA})\big)>\lambda^+$.
\end{proposition}

\begin{proof}
First we define a $\lambda$--tree creating pair $(K(\bar{\cA}),\Sigma(
\bar{\cA}))=(K,\Sigma)$. For $\alpha<\lambda$ let $\bH(\alpha)=
\kappa_\alpha$ and let ${\rm st}_\alpha$ be a winning strategy of player II
in the game $\Game^*(\cA_\alpha,\lambda_\alpha)$.  

$K$ consists of all $\lambda$--tree creatures $t\in\TCR[\bH]$ such that
letting $\alpha=\lh(\eta[t])$: 
\begin{itemize}
\item {\bf either}\ $\dis[t]=(\delta,\langle A^t_i:i<\delta\rangle)$, where
$\delta<\lambda_\alpha$ and $\langle A^t_i:i<\delta\rangle$ is (an initial
segment of) a play of $\Game^*(\cA_\alpha,\lambda_\alpha)$ in which player
II uses strategy ${\rm st}_\alpha$,\\
{\bf or}\ $\dis[t]=\langle \xi\rangle$ for some $\xi<\kappa_\alpha$;
\item if $\dis[t]=\langle\xi\rangle$, then $\pos[t]=\{\eta[t]\conc\langle\xi
\rangle\}$ and $\nor[t]=0$;
\item if $\dis[t]=(\delta,\langle A^t_i:i<\delta\rangle)$, then $\pos[t]=\{
\eta[t]\conc\langle\xi\rangle:\xi\in\bigcap\limits_{i<\delta} A_i\}$ and
$\nor[t]=\alpha+1$. (If $\delta=0$ then we stipulate $\pos[t]=\{\eta[t]\conc
\xi:\xi<\kappa_\alpha\}$.)
\end{itemize}
The domain of the tree composition operation $\Sigma$ consists of singletons
only, and\\  
if $\nor[t]=0$ then $\Sigma(t)=\{t\}$,\\
if $\nor[t]>0$, $\alpha=\lh(\eta[t])$ and $\dis[t]=(\delta,\langle A^t_i:i<
\delta\rangle)$, then $\Sigma(t)$ consists of those $t'\in K\cap\TCR_{\eta[t
]}[\bH]$ for which:
\begin{itemize}
\item either $\nor[t']=0$ and $\pos[t']\subseteq\pos[t]$,
\item or $\nor[t']>0$, $\dis[t']=(\delta',\langle A^{t'}_i:i<\delta'\rangle
)$ and $\langle A^t_i:i<\delta\rangle\trianglelefteq\langle A^{t'}_i:i<
\delta'\rangle$.
\end{itemize}

\begin{claim}
\label{cl5}
$(K,\Sigma)$ is an exactly $\bar{\lambda}$--complete very local tree
creating pair for $\bH$. Hence the forcing notion $\bbQ^\tree_{D^*}(K,
\Sigma)$ is fuzzy proper for $W$.
\end{claim}

\begin{proof}[Proof of the Claim]
Should be clear.
\end{proof}

We finish the proof of the proposition in a standard way: we force with
$\lambda$--support iteration, $\lambda^{++}$ in length, of the forcing
notion $\bbQ^\tree_{D^*}(K(\bar{\cA}),\Sigma(\bar{\cA}))$.
\end{proof}


\def\germ{\frak} \def\scr{\cal} \ifx\documentclass\undefinedcs
  \def\bf{\fam\bffam\tenbf}\def\rm{\fam0\tenrm}\fi 
  \def\defaultdefine#1#2{\expandafter\ifx\csname#1\endcsname\relax
  \expandafter\def\csname#1\endcsname{#2}\fi} \defaultdefine{Bbb}{\bf}
  \defaultdefine{frak}{\bf} \defaultdefine{mathfrak}{\frak}
  \defaultdefine{mathbb}{\bf} \defaultdefine{mathcal}{\cal}
  \defaultdefine{beth}{BETH}\defaultdefine{cal}{\bf} \def\bbfI{{\Bbb I}}
  \def\mbox{\hbox} \def\text{\hbox} \def\om{\omega} \def\Cal#1{{\bf #1}}
  \def\pcf{pcf} \defaultdefine{cf}{cf} \defaultdefine{reals}{{\Bbb R}}
  \defaultdefine{real}{{\Bbb R}} \def\restriction{{|}} \def\club{CLUB}
  \def\w{\omega} \def\exist{\exists} \def\se{{\germ se}} \def\bb{{\bf b}}
  \def\equivalence{\equiv} \let\lt< \let\gt> \def\implies{\Rightarrow}

\end{document}